%% file: zakharov.tex
\documentclass[11pt]{article}
\usepackage{amsmath}
\usepackage{amsfonts}
\usepackage{amsthm}
\usepackage{ifthen}
\usepackage[latin1]{inputenc}
\usepackage{graphicx}
\usepackage{epsfig}

\def\a{\alpha}

\def\b{\beta}
\def\d{\partial}
\def\g{\gamma}
\def\lam{\lambda}
\def\s{\sigma}
\def\om{\omega}

\def\ompi{{\rm op} }
\def\ompie{{\rm op}_\e^\psi }

\def\C{\mathbb{C}}

\def\phi{\varphi}
\def\e{\varepsilon}
\def\S{\mathbb{S}}

\def\R{\mathbb{R}}
\def\N{\mathbb{N}}

\def\Z{\mathbb{Z}}

\def\un{\underline}
\def\id{{\rm Id }\!}

\newtheorem{defi}{Definition}[section]
\newtheorem{theo}[defi]{Theorem}
\newtheorem{prop}[defi]{Proposition}
\newtheorem{lem}[defi]{Lemma}

\newtheorem{ass}[defi]{Assumption}

\numberwithin{equation}{section}

\begin{document}

\renewcommand{\refname}{References}

\title{Derivation of the Zakharov equations }
\author{Benjamin Texier\thanks{
Indiana University, Bloomington, IN 47405;
btexier@indiana.edu.
 This research was partially supported
under NSF grant number DMS-0300487. The author warmly thanks Christophe Cheverry, Thierry Colin, David Lannes, Guy M\'etivier, and Kevin Zumbrun, for the interest they showed for this work and many very interesting discussions.}}
\maketitle

 \vspace{2cm}
 {\footnotesize {\bf Abstract} - This paper continues the study, initiated in \cite{T-AA,CEGT}, of the validity of the Zakharov model describing Langmuir turbulence. We give an existence theorem for a class of singular quasilinear equations. This theorem is valid for well-prepared initial data. We apply this result to the Euler-Maxwell equations describing laser-plasma interactions, to obtain, in a high-frequency limit, an asymptotic estimate that describes solutions of the Euler-Maxwell equations in terms of WKB approximate solutions which leading terms are solutions of the Zakharov equations. Because of transparency properties of the Euler-Maxwell equations put in evidence in \cite{T-AA}, this study is led in a supercritical (highly nonlinear) regime. In such a regime, resonances between plasma waves, electromagnetric waves and acoustic waves could create instabilities in small time. The key of this work is the control of these resonances. The proof involves the techniques of geometric optics of Joly, M\'etivier and Rauch \cite{JMR_Ind,JMR_TMB}, recent results of Lannes on norms of pseudodifferential operators \cite{L0}, and a semiclassical, paradifferential calculus.}

\newpage 

\tableofcontents

\newpage 
\begin{section}{Introduction}

 We describe solutions of initial value problems for quasilinear, hyperbolic systems of the form,
 \begin{equation} \label{em3_0}
    \d_t u + \frac{1}{\e^2} {\cal A}(\e, \e u, \e \d_x) u  =  \frac{1}{\e} {\cal B}(u, u),
\end{equation}
 in the high-frequency limit $\e \to 0.$ 

  In \eqref{em3_0}, ${\cal A}$ is a symmetric hyperbolic, differential, or pseudo-differential operator; 
 the singular source term ${\cal B}$ is bilinear. The unknown $u^\e$ has values in $\R^n.$ It depends on time $t \in \R_+$ and space $x \in \R^d,$ and is subject to the initial condition,
 \begin{equation} \label{em3_1}
  u^\e(0,x) = a^\e(x),
 \end{equation}
where $a^\e$ is a bounded family in $H^s(\R^d),$ for some large $s.$

  In this setting, the existence, uniqueness and regularity of solutions to (\ref{em3_0})-(\ref{em3_1}) for fixed $\e > 0$ is classical. 

 The limit $\e \to 0$ is singular in two ways: first, solutions develop fast oscillations in time, with frequencies of typical size $O(1/\e^2);$ second, the amplitude $O(1)$ of the initial datum is large, hence the singular source term ${\cal B}/\e$ could create instabilities in small time $O(\e).$ 


 Under appropriate assumptions, we prove the existence of solutions to \eqref{em3_0}-\eqref{em3_1} over time intervals independent of $\e,$ and their stability under initial perturbations of the form $\e^{k_0} \varphi^\e,$ where $k_0$ is large enough, and $\varphi^\e$ is bounded in a semiclassical Sobolev space, in particular, may contain fast oscillations of the form $e^{i k x / \e}.$

 We show that our assumptions are satisfied by the Euler-Maxwell equations describing laser-plasma interactions. This implies in particular that, in a high-frequency limit, solutions of the Euler-Maxwell equations are well approximated by WKB approximate solutions which leading terms are solutions of the Zakharov equations. 

 Our assumptions and results are precisely stated in section \ref{AR}.


\begin{subsection}{Discussion: weakly nonlinear versus highly nonlinear geometric optics} \label{wnl-vs-of}

 Consider \eqref{em3_0}, and assume for instance that ${\cal A}$ is a differential operator of the form 
  \begin{equation} \label{dispersif} {\cal A}(\e,\e u, \e \xi) = {\cal A}_0(\e) + \e {\cal A}_1(\e,\xi) + \e^2 {\cal A}_2(\e,u,\xi),\end{equation}
 where ${\cal A}_1$ and ${\cal A}_2$ are linear in $\xi,$ and ${\cal A}_0,$ ${\cal A}_1$ and ${\cal A}_2$ are hermitian, for all $\e, u, \xi.$ Suppose that the  family of initial data has the form
 $ u^\e(0,x) = \e^p a^\e(x),$
 where $a^\e$ is a bounded family in $H^s(\R^d).$

  The \emph{weakly nonlinear} regime corresponds to $p =1.$ In this regime, the classical existence proof provides a maximal existence time $t^*(\e)$ that satisfies $\liminf_{\e \to 0} t^*(\e) > 0.$ Indeed, if one lets $v^\e = \e u^\e,$ then the initial datum for $v^\e$ is $O(1),$ and the equation in $v^\e$ is   $$ \d_t v^\e + \frac{1}{\e^2} {\cal A}(\e, \e v^\e, \e \d_x) v^\e = {\cal B}(v^\e, v^\e).$$
 The classical $H^s$ energy estimate for quasilinear symmetric hyperbolic operators then yields
 $$ |v^\e(t) |_{H^s} \leq | a^\e |_{H^s} + C \int_0^t |v^\e (t') |_{H^s} \, dt',$$
 where $C$ depends on $|v^\e|_{W^{1, \infty}},$ then, with Gronwall's lemma, the uniform bound
 $$  |v^\e(t) |_{H^s} \leq | a^\e |_{H^s} e^{C t}.$$
 This uniform estimate is the key of the proof of the existence of a solution over a time interval independent of $\e.$

 On the contrary, when $p = 0,$ the initial data have a \emph{large amplitude} $O(1).$ In this regime, the maximal existence time a priori satisfies $t^*(\e) = O(\e),$ and in particular $\liminf_{\e \to 0} t^*(\e) = 0.$ The $H^s$ energy estimate gives indeed  
 $$ | u^\e(t) |_{H^s} \leq | a^\e |_{H^s} + \frac{C}{\e} \int_0^t |u^\e (t') |_{H^s} \, dt',$$
 whence,
 $$ |u^\e(t)|_{H^s} \leq | a^\e |_{H^s} e^{C t/\e}.$$
  This shows that the singular term ${\cal B}/\e$ may cause the solution to blow-up in small time. This highly nonlinear (or supercritical) regime is called \emph{oscillations fortes}, or \emph{strong oscillations}, by Cheverry, Gu\`es and M\'etivier in \cite{CGM}.

 \emph{ Our goal is to state conditions on \emph{(\ref{em3_0})-(\ref{em3_1})} that are sufficient to have $\liminf_{\e \to 0} t^*(\e) > 0$ in the regime $p = 0,$ and that are satisfied by the Euler-Maxwell equations.}
\end{subsection}


\begin{subsection}{The Euler-Maxwell equations and the Zakharov equations} \label{EMZ} 

  This paper is in the direct continuation of \cite{T-AA, CEGT}. The underlying physical context is the study of laser-plasma interactions; in particular, the question of the rigorous derivation of the Zakharov model from fundamental equations. 

  We take here as a system of fundamental equations the Maxwell equations coupled with the Euler equations \cite{DB,SS},
 $$ \begin{aligned} \d_{t} B^{\flat} +  c\,  \nabla \times  E^{\flat} & =   0,  \\ \d_{t}  E^{\flat} -  c\, \nabla \times  B^{\flat}  & =  4 \pi e ((n_0 + n_e^{\flat}) v_e^{\flat} - (n_0 + n_i^{\flat}) v_i^{\flat})), \\ 
 m_e (n_0 + n_e^{\flat}) ( \d_{t} v_e^{\flat} + (v_e^{\flat} \cdot  \nabla ) v_e^{\flat})  & =  - \gamma_e T_e \nabla n_e^{\flat} - e (n_0 + n_e^{\flat}) (E^{\flat} + \frac{1}{c} v_e^{\flat} \times B^{\flat}), \\
 m_i  (n_0 + n_i^{\flat}) (\d_{t} v_i^{\flat} + (v_i^{\flat} \cdot  \nabla ) v_i^{\flat})  & =  - \gamma_i T_i \nabla n_i^{\flat} + e (n_0 + n_e^{\flat})(E^{\flat} + \frac{1}{c} v_i^{\flat} \times B^{\flat}), \\
 \d_{t} n_e^{\flat} + \nabla \cdot ( (n_0 + n_e^{\flat}) v_e^{\flat} )  & =  0,  \\
 \d_{t} n_i^{\flat} + \nabla \cdot ( (n_0 + n_e^{\flat}) v_i^{\flat} )  & = 0. \end{aligned}$$
  The variables are $B^{\flat}, E^{\flat}$ the electromagnetic field, $v_e^{\flat}, v_i^{\flat}$ the velocities of the electrons and of the ions, and $n_e^{\flat}$ and $n_i^{\flat}$ the density fluctuations from the equilibrium $n_0.$ The first two equations are Maxwell's equations describing the time evolution of the electromagnetic field, the next two equations are the equations of conservation of momentum for the electrons and the ions, and the last two equations are the equations of conservation of mass for the electrons and the ions. The electric charge of the electrons is $-e;$ to simplify, we assume that the charge of the ions is $+e.$ The parameters are $m_e$ and $m_i$ the masses of both species, $\gamma_e$ and $\gamma_i$ the specific heat ratios of both species, $T_i$ and $T_e$ the temperatures of both species and $n_0$ the (assumed constant and isotropic) density of the plasma at equilibrium.

 In the above system, Maxwell's equations are coupled to Euler's equations by the current density term in the right hand side of the Amp\`ere equation and by the Lorentz force in the right hand side of the equations of conservation of momentum. 

 The additional divergence equations
 \begin{equation} \label{3fev} \nabla \cdot B^{\flat} = 0, \quad \nabla \cdot E^{\flat} = 4 \pi e (n_e^{\flat} - n_i^{\flat}),\end{equation} 
  are satisfied at all times if they are satisfied by the initial data.
 A brief discussion of the relevance of this model is given in \cite{T-AA}. We work in this paper on the non-dimensional form of these equations introduced in \cite{T-AA}:
 $$\mbox{(EM})^\sharp \left\{ \begin{aligned}   \d_t B +   \nabla \times E  =   0 , \nonumber \\ \d_t E -  \nabla \times  B  =  \frac{1}{\e} (1 + n_e^\sharp) v_e - \frac{1}{\e} \frac{\theta_i}{\theta_e} ( 1  + n_i^\sharp) v_i, \nonumber \\
   \d_t v_e + \theta_e (v_e \cdot \nabla ) v_e  =  - \theta_e  \frac{\nabla n_e^\sharp}{1 + n_e^\sharp} -  \frac{1}{\e} ( E + \theta_e v_e \times B), \nonumber \\
  \d_t n_e^\sharp + \theta_e  \nabla \cdot( (1 + n_e^\sharp) v_e)  = 0, \\
\d_t v_i +  \theta_i (v_i \cdot \nabla) v_i  =  - \alpha^2 \theta_i  \frac{\nabla n_i^\sharp}{1 + n_i^\sharp} + \frac{1}{\e} \frac{\theta_i}{\theta_e} ( E +  \theta_i v_i \times B),\nonumber \\
  \d_t n_i^\sharp + \theta_i  \nabla \cdot( (1 + n_i^\sharp) v_i)  = 0.
  \end{aligned} \right. $$
In $\mbox{(EM})^\sharp,$ the variable is 
 $$ u^\sharp = (B, E, v_e, n_e^\sharp, v_i, n_i^\sharp) \in \R^{14}.$$
The change of variables for small amplitudes,
 \begin{equation} \label{chgt} 1 +  n_e^\sharp = e^{n_e},\;\;  1 +  n_i^\sharp = e^{n_i},\end{equation}
  leads to the system,
 \begin{equation} \mbox{(EM)} \left \{\begin{aligned} \d_{t} B + \nabla \times E  & =   0,  \nonumber \\  \d_{t} E -  \nabla \times  B  & =  \frac{1}{\e} ( e^{n_e} v_e - \frac{\theta_i}{\theta_e} e^{n_i} v_i ), \nonumber \\
   \d_{t} v_e + \theta_e (v_e \cdot \nabla) v_e  & =  - \theta_e  \nabla n_e -  \frac{1}{\e} (E + \theta_e v_e \times B), \nonumber \\
  \d_{t} n_e + \theta_e \nabla \cdot v_e   +  \theta_e  (v_e \cdot \nabla) n_e  & =  0,  \nonumber \\
\d_{t} v_i + \theta_i (v_i \cdot \nabla) v_i  & =  - \alpha^2 \theta_i \nabla n_i +  \frac{\theta_i}{\e \theta_e} ( E + \theta_i v_i \times B),\nonumber \\
  \d_{t} n_i +  \theta_i \nabla \cdot  v_i   +  \theta_i  v_i \cdot \nabla  n_i  & = 0. \nonumber 
						\end{aligned}\right.\end{equation} 
						In (EM), the variable is 
  $$ \tilde u = (B, E, v_e, n_e, v_i, \frac{n_i}{\a}) \in \R^{14},$$
 where $(B, E) \in \R^{3 + 3}$ is the electromagnetic field, $(v_e, v_i) \in \R^{3+3}$ are the velocities of the electrons and of the ions and $(n_e, n_i) \in \R^{1+1}$ are the fluctuations of densities of both species. The variable $\tilde u$ depends on time $t \in \R_+$ and space $x\in \R^3.$  The small parameter $\e$ is
 $$ \e := \frac{1}{\om_{pe} t_0},$$
 where $\om_{pe}$ is the electronic plasma frequency:
 \begin{equation} \label{o-p} \om_{pe} :=  \sqrt{\frac{4 \pi e^2 n_0}{m_e}}, \end{equation}
 and $t_0$ is the duration of the laser pulse. A typical value for $\e$ in realistic physical applications is $\e \simeq 10^{-6}.$
 The parameters $\a, \theta_e$ and $\theta_i$ are 
  $$ \theta_e := \frac{1}{c} \sqrt{\frac{\gamma_e T_e}{m_e}}, \quad \theta_i := \frac{1}{c} \sqrt{\frac{\gamma_e T_e}{m_i}}, \quad \alpha := \frac{T_i}{T_e},$$  
 Typically, $\theta_e \simeq 10^{-3}.$ Because the ions are much heavier than the electrons, $\theta_i$ and $\a$ are much smaller than $\theta_e.$ We consider the specific regime
 \begin{equation} \label{theta_i} \theta_i = \e,\end{equation}
 and we look for solutions to (EM) with initial data of size $O(\e)$ defined over diffractive times $O(1/\e),$ that is, we make the ansatz
 \begin{equation} \label{regime} u(t,x) := \e \tilde u(\e t, x). \end{equation}
 Written as a system of equations in the variable $u,$ (EM) takes the form \eqref{em3_0}. 
  
 The Zakharov system is a simplified model for the description of the nonlinear interactions between $\tilde E,$ the envelope of the electric field, and $\bar n,$ the mean mode of the ionic fluctuations of density in the plasma,
 $$ \mbox{(Z)} \left\{ \begin{aligned}  i \d_t  \tilde  E  + \Delta \tilde E  & =   \bar n  \tilde E, \\
   \d_t^2 \bar n - \Delta \bar n & =  \Delta | \tilde E|^2, \end{aligned}\right.$$

 This model was derived by Zakharov and his collaborators in the seventies \cite{Zak}. It describes nonlinear interactions between high-frequency, electromagnetic waves and low-frequency, acoustic waves. In (Z), the Schr\"{o}dinger operator is the classical three-scale approximation of Maxwell's equations \cite{JMR_Ind}; the wave operator is the classical long-wave approximation of the Euler equations. The nonlinear term in the right-hand side of the equation in $\tilde E$ directly comes from the current density term in the Amp\`ere equation. The term $\Delta | \tilde E|^2$ comes from the convective terms and the nonlinear force term in the equations of conservation of momentum. A WKB expansion of the $\mbox{(EM)}^\sharp$ system is performed in section \ref{3-1}. In the nonlinear regime of our interest, the limit system is (Z).  

 

\end{subsection}


\begin{subsection}{Description of the results}

  We extend here the results of \cite{T-AA,CEGT}, as we show that, in the high-frequency limit and in a highly nonlinear regime, solutions of the Euler-Maxwell equations are well-approximated by solutions of the Zakharov equations. In section \ref{deux}, under appropriate assumptions, we prove the following theorem:
 \begin{theo} \label{th-main} The unique solution to {\rm (\ref{em3_0})-(\ref{em3_1})} is defined over a time interval independent of $\e.$
\end{theo}
 Then we apply this result to the Euler-Maxwell equations, to obtain: 
 \begin{theo} \label{th-EM} In the high-frequency limit, solutions of the Euler-Maxwell equations initiating from polarized, large-amplitude initial data, are well approximated by solutions of the Zakharov equations initiating from nearby initial data, in the sense that there exists $t_0 > 0$ and $C > 0,$ independent of $\e,$ such that,
  $$ \sup_{0 \leq t \leq t_0} \sup_{x} ( | E - (\tilde E e^{i \om_{pe} t / \e^2} + c. c.)| + | n^\e - \e \bar n | ) \leq C \e^2,$$
  for $\e$ small enough, where $E^\e$ and $n^\e$ represent the electric field and the (electronic or ionic) fluctuation of density in the solution of the {\rm (EM)} system, and $\tilde E, \bar n,$ is the solution of {\rm (Z)}. 
  \end{theo}
  
  We show in section \ref{application} how Theorem \ref{th-EM} follows from Theorem \ref{th-main}. Precise statements are given in sections \ref{AR} (Therorem \ref{th1}) and \ref{3-2} (Theorem \ref{theo}).

\end{subsection}

\begin{subsection}{Outline of the proof} 

  The proof of Theorem \ref{th-main} (section \ref{pr}) goes along the following lines: the construction of a precise, regular, polarized \emph{approximate solution} defined over a time interval independent of $\e,$
 the \emph{preparation} of the system, and finally, the \emph{control of the resonant interactions of oscillating waves}.

 \begin{subsubsection}{Existence of an approximate solution}

  In section \ref{application}, we consider the initial value problem for the (EM) system, and show that it takes the form \eqref{em3_0}-\eqref{em3_1}. We construct a WKB approximate solution,
  \begin{equation} \label{wkb-0} 
   u_a^\e = u_0 + \e u_1 + \dots,
   \end{equation}
   under an assumption of polarization for the initial data. 
 
   The polarization condition is a well-preparedness condition; it is \emph{necessary} for the construction of a WKB approximate solution that is continuous in time, at $t = 0.$ Consider indeed $u_a^\e$ in the form \eqref{wkb-0}, a solution to \eqref{em3_0}-\eqref{em3_1}, where the operator ${\cal A}$ has the form \eqref{dispersif}, with ${\cal A}_0(0) \neq 0,$ as in the (EM) system. The limits $\e \to 0$ yields
   $$ {\cal A}_0(0) u_0 = 0,$$
   then, if $u_0$ is continuous at $t = 0,$
   $$   {\cal A}_0(0) a^0 = 0.$$ The above condition is the well-preparedness assumption for the initial datum $a^\e.$ For the (EM) system, it takes the form
  $$ a^0 = \big( 0_{\R^3}, \tilde E + (\tilde E)^*, \frac{i}{\om} \tilde E - \frac{i}{\om} \tilde E^*, 0_{\R}, 0_{\R^3}, 0_{\R} \big),$$
  for some fundamental frequency $\om,$ defined in terms of $\om_{pe},$ and some complex amplitude $\tilde E$ (above, $(\tilde E)^*$ denotes the complex conjugate of $\tilde E$).  
  
  WKB solutions to the Euler-Maxwell equations, initating from highly-oscillating, and well-prepared, initial data, are considered in \cite{T-AA}. It is shown in \cite{T-AA} that,
   \begin{itemize}
   \item[1)] the (EM) equations satisfy \emph{transparency} properties, that is, null conditions for coefficients describing constructive interactions of characteristic waves. As a result, the weakly nonlinear (in the sense of section \ref{wnl-vs-of}) approximation of the (EM) system is a \emph{linear} transport equation, and
    \item[2)] WKB solutions of the (EM) equations, initiating from large-amplitude solutions, satisfy, in the high-frequency limit $\e \to 0,$ systems of the form
   $$ \mbox{(Z)}_c \left\{ \begin{aligned} i (\d_t + c \d_z) E + \Delta_\perp E & = n E, \\ (\d_t^2 - \Delta_\perp) n & = \Delta_\perp (|E|^2), \end{aligned} \right. $$
  where $z$ is the direction of propagation of the laser pulse, and $\Delta_\perp$ is the Laplace operator in the transverse directions. 
  \end{itemize}
 
  The approximate solution that is constructed in section \ref{application} satisfies the ansatz,
   \begin{equation} \label{ansatz-intro} u_a^\e(t, x) = U_a^\e (t, x, \frac{\om t}{\e^2}).\end{equation}
  
  In particular, there are \emph{three} times scales. This is consistent with the well-known fact that the Schr\"{o}dinger equation is an approximation of the Maxwell equations in the diffractive limit (that is, $t = O(1)$ and oscillations with frequencies $O(1/\e^2)$). Note, however, that the wave equation, also present in the (Z) system, is an approximation of the Euler equation in the geometric optics limit (that is, $t = O(1)$ and oscillations in $O(1/\e)$). The third scale is actually built in the Euler equations by the "cold ions" assumption $\theta_i = \e.$ 
   
  In \eqref{ansatz-intro}, the profiles are purely time-oscillating. In particular, the initial data do not have fast oscillations.  The limit system is $\mbox{(Z)}_0,$ that is, a Zakharov system with zero group velocity (see the characteristic variety pictured on figure \ref{fig-char1}).  Such waves are called \emph{plasma waves} in the physical literature. The approximate solution has the form
  \begin{equation} \label{10fev} u_a^\e = (u_{0,1} e^{i \om t/\e^2} + u_{0,1}^* e^{- i \om t/\e^2}) + \e(u_{1,0} + \dots) + \e^2(\dots),\end{equation}
   where the components of $u_{0,1}$ and $u_{1,0}$ satisfy (Z). 
  
 For general equations of the form (\ref{em3_0}), studied in section \ref{deux}, the existence of an approximate solution is assumed (Assumption \ref{ass_5}).
 
  \end{subsubsection}
 \begin{subsubsection}{Preparation of the system} \label{142}

  Given a precise approximate solution $u_a^\e$ of the form \eqref{10fev}, we look in section \ref{251} for the exact solution $u^\e$ as a perturbation of $u_a^\e,$
   $$ u^\e = u_a^\e + \e^k \dot u^\e.$$
 The initial condition is $u^\e( t = 0) = a^\e + \e^{k_0} \phi^\e,$ where $\phi^\e$ has a high Sobolev regularity and $k_0$ is large enough. In the definition of $\dot u^\e,$ $k$ is chosen in terms of $k_0.$ One assumes that $u_a^\e$ is accurate at an order $l_0,$ much larger that $k_0.$ The equation in $\dot u^\e$ has the form,  
  \begin{equation} \label{10fev2} \d_t \dot u^\e + \frac{1}{\e^2} {\cal A}(\e (u_a^\e + \e^k \dot u^\e)) \dot u^\e  = \frac{1}{\e} ({\cal B}(u_a^\e, \dot u^\e) + {\cal B}(\dot u^\e, u_a^\e)),\end{equation}
    The propagator ${\cal A}/\e^2$ is hyperbolic (Assumption \ref{ass_1}), thus generates highly oscillating waves, with frequencies of typical size $O(1/\e^2).$ We write the spectral decomposition of ${\cal A}$ as follows:
     $$ {\cal A} = \sum_{\mbox{{\footnotesize kg}}} i \lam_{\mbox{{\footnotesize kg}}} \Pi_{\mbox{{\footnotesize kg}}} + \sum_{\mbox{{\footnotesize ac}}} i\lam_{\mbox{{\footnotesize ac}}} \Pi_{\mbox{{\footnotesize ac}}}.$$
      The real eigenvalues $\lam_{\mbox{{\footnotesize kg}}}$ are called Klein-Gordon modes, while the real eigenvalues $\lam_{\mbox{{\footnotesize ac}}}$ are called acoustic modes. The characteristic variety for the (EM) system (that is, the union of the graphs $\xi \mapsto \lam_j(\xi),$ at $u = 0,$ for $j = \mbox{kg}$ and $j = \mbox{ac}$) is pictured on figure \ref{fig-char1}. In particular, $\lam_{\mbox{{\footnotesize kg}}} \sim O(1) \xi,$ while $\lam_{\mbox{{\footnotesize ac}}} \sim \e \xi,$ a consequence of the cold ions hypothesis \eqref{theta_i}. 
      
   The Klein-Gordon waves generated by ${\cal A}$ interact with the highly oscillating approximate solution, through the convection term, and through the source term ${\cal B}.$ These interactions create low-frequency waves, which can be seen as source terms in the equations for the components of the solutions in the directions of the acoustic modes. Thus low-frequency, and high-frequency waves are propagated. The Zakharov system pretends to describe how these waves interact. 
   
     Equation \eqref{10fev2}, together with an initial datum of size $O(1),$ can be likened to an ordinary differential equation,
  $$ y'  + \frac{i \a}{\e^2} y = \frac{1}{\e} y^2,$$
  with an initial datum $y(0) = y_0.$ The singular source term in the right-hand side may cause the solution to blow-up in small time $O(\e),$ but exponential cancellations are expected to happen because of the rapid oscillations created by the term in $1/\e^2.$ 
  
   To investigate these exponential cancellations, it is natural to project the source term ${\cal B}/\e$ over the eigendirections of ${\cal A}.$ Let the total eigenprojectors, 
   $$\Pi_0 = \sum_{{\rm kg}} \Pi_{{\rm kg}}, \quad \Pi_s = \sum_{{\rm ac}} \Pi_{{\rm ac}}.$$
    We compute
   $$ \left( \begin{array}{cc} \Pi_0 {\cal B} \Pi_0 & \Pi_0 {\cal B} \Pi_s \\ \Pi_s {\cal B} \Pi_0 & \Pi_s {\cal B} \Pi_s \end{array}  \right) = \left(\begin{array}{cc} * & O(1) \\ ** & * \end{array} \right).$$
  Because constructive interaction of waves between low- and high-frequency do occur for the (EM) system (see figure \ref{fig-0-s}), the $O(1)$ term in the right block of the above interaction matrix can be interpreted as an \emph{absence of transparency.} 
  
  We then rescale the solution (section \ref{resc}) as follows, 
   $$ v^\e := (\Pi_0 \, \dot u^\e, \frac{1}{\e} \Pi_s \,  \dot u^\e).$$
  The equation in $v^\e$ has the form
 $$ \d_t v^\e + \frac{1}{\e^2} A v^\e = \frac{1}{\e^2} B v^\e + \frac{1}{\e} D v^\e + O(1) v^\e,$$
  with the notations, 
  $$ B = \left(\begin{array}{cc} 0 & 0 \\ ** & 0 \end{array}\right), \quad D := \left(\begin{array}{cc} * & 0 \\ 0 & * \end{array}\right).$$
  In the equation for $v,$ the propagator $A$ is diagonal and the leading source term $B$ is nilpotent. The system is prepared.    
    
  \end{subsubsection}
  \begin{subsubsection}{Control of the constructive interactions of waves} \label{143}
 
In a third step (sections \ref{red0} to \ref{red3}), the singularity in the right-hand side in $v^\e$ is reduced. Consider a change of variable in the form    
 $$ w^\e := (\id \, + N)^{-1} v^\e, \quad N = \left(\begin{array}{cc} 0 & 0 \\ N & 0 \end{array} \right).$$
  The equation satisfied by $w^\e$ is
  $$ \d_t w^\e +  \frac{1}{\e^2} A w^\e = \frac{1}{\e^2} (B - [\e^2 \d_t + A, N])  w^\e + \frac{1}{\e} D w^\e.$$
  We look for $N$ solution of the \emph{homological equation}, 
  $$  B - [\e^2 \d_t + A, N] = O(\e^2),$$
  This equation takes the form  $ \Phi N = B + O(\e^2),$
   where the phase $\Phi$ is 
$ \Phi = \lam_{{\rm kg}} - \lam_{{\rm ac}} - \om.$
  The equation $\Phi = 0$ is the \emph{resonance} equation. Its solutions are pictured in Figures \ref{fig-0-0}, \ref{fig-0-s} and \ref{fig-0-0-s}. The crucial transparency assumption (Assumption \ref{ass_6}) states that the interaction coefficient $B$ is sufficiently small at the  resonances, that is, 
  \begin{equation} \label{int00} |B| \leq C \e^2 | \Phi|.
  \end{equation}
   The equation in $w^\e$ becomes
  $$ \d_t w^\e + \frac{1}{\e^2} A w^\e = \frac{1}{\e} D w^e.$$
  A symmetrizability assumption for $D$ (Assumption \ref{ass_sym}) eventually allows to perform energy estimates (section \ref{256}), which yield uniform bounds for $w^\e,$ and a continuation argument concludes the proof.   
  
  In section \ref{application}, we describe the resonance equations for the Euler-Maxwell equations, and check that an estimate of the form \eqref{int00} is satisfied. 
 \end{subsubsection}
 
 \begin{subsubsection}{Technical issues} \label{144} 
   
  The symbols in the spectral decompositions of ${\cal A}$ are necessarily pseudo-differential operators, even when ${\cal A}$ is differential. They also depend on the solution $u^\e,$ because the equations are nonlinear. We are naturally led to consider pseudo-differential operators of the form
   $$ q(\e, x, \xi) = p(\e, v(x), \xi),$$
  where $v$ has a Sobolev regularity. The questions of the bounds of the corresponding operators, in Sobolev spaces, and of the existence of a symbolic calculus, naturally arise. Lannes recently gave optimal bounds in \cite{L0}. For an operator of order $m,$ these bounds have the form,
   $$ \| \ompi(p(v, \xi) u \|_{H^s} \leq C(| v|_{L^\infty}) (\| v \|_{H^{s}} \| u \|_{H^{s_0}} + \| u \|_{H^{s + m}}).$$
   A symbolic calculus is available; the operator $\ompi(p_1(v)) \ompi(p_2(v))$ has the symbol
  $$ p_1(v) p_2(v) + \sum_{|\a| =1} \d_\xi^\a (p_1(v)) \d_x^\a (p_2(v)) + \dots$$
  When the symbols depend on $x$ through the solution $u^\e,$ the subprincipal symbol depends on $\d_x u^\e,$ and its operator norm in $H^s$ depends on $ \| u^\e \|_{H^{s + 1}}.$ That is, compositions of such operators lead to losses of derivatives.

   To overcome this difficulty, it is classical to differentiate the equation up to order $s,$ and then perform energy estimates in $L^2.$

  We now explain why we chose a different approach.

  In the perturbation equations \eqref{10fev2}, all the derivatives are $\e\mbox{-derivatives}$ (see below). The equation in $\e \d_x \dot u^\e$ (we let $d = 1$ in this discussion), has a singular source term in $\dot u^\e.$ The variable $U^\e := (\dot u^\e, (\e \d_x) \dot u^\e),$ solves 
 $$ \d_t U^\e + \frac{1}{\e^2} {\cal A}(\e u^\e) U^\e  = \frac{1}{\e} \left( \begin{array}{cc} {\cal B} & 0 \\ \e {\cal B}' & {\cal B} \end{array}\right) U^\e - \frac{1}{\e^2} \left(\begin{array}{c} 0 \\ {} [\e \d_x, {\cal A}(\e u^\e)] \end{array}\right) U^\e,
  $$
  where ${\cal B}$ is short for ${\cal B}(u_a^\e)$ and ${\cal B}'$ is short for ${\cal B}(\d_x u_a^\e).$ Up to $O(\e^2),$ the commutator in the above equation is $\e \d_v {\cal A}(u_a^\e).$ Following the approach of section \ref{142}, we rescale the solution, by letting 
  $$V^\e = ( \Pi_0 \dot u^\e, \frac{1}{\e} \Pi_s \dot u^\e,  \Pi_0 \e \d_x \dot u^\e, \frac{1}{\e} \Pi_s \e \d_x \dot u^\e).$$
   The equation becomes, 
   $$ \d_t V^\e + \frac{1}{\e^2} A(\e u^\e) V^\e  = \frac{1}{\e^2} \un B V^\e + \frac{1}{\e} \un D V^\e,$$
 with the notation,
  $$ \un B := \left(\begin{array}{cccc} 0 & 0 & 0 & 0 \\ \Pi_s {\cal B} \Pi_0 & 0 & 0 & 0 \\ 0 & 0 & 0 & 0 \\ \e( \Pi_s {\cal B}' - \d_v {\cal A}(u_a)) \Pi_0 & 0 & \Pi_s {\cal B} \Pi_0 & 0 \end{array} \right),$$
  Because of the singular rescaling, the commutator $\Pi_s \d_v {\cal A}(u_a)) \Pi_0$ is now a singular source term. It cannot be symmetrized in the matrix $\un B,$ and might well be not transparent. That is, $\un B$ might not satisfy an estimate of the form \eqref{int00}, even if ${\cal B}$ does. In the case of the (EM) system, this term can actually be shown to satisfy a transparency estimate of the form \eqref{int00}, so that the method of differentiating the equation could be applied. To obtain a $H^s_\e$ estimate, one would have to write a system of size $2 n \sum_{s' = 0}^s {\footnotesize\left(\begin{array}{c} s' + d - 1 \\ d-1 \end{array}\right)},$ with a leading singular source term in a bidiagonal form
  $$ \frac{1}{\e^2} \left(\begin{array}{cccc} {\cal B} & 0 & 0 & 0 \\ \e {\cal B}' & {\cal B} & 0 & 0 \\ 0 & \ddots & \ddots & 0 \\ 0 & 0 & \e {\cal B}' &  {\cal B} \end{array}\right),$$
 and notational complications would arise in the normal form reductions of sections \ref{red0} to \ref{red3}. This method would arguably be conceptually simpler than the one we chose, namely paradifferential smoothing.

 The paradifferential smoothing of Bony \cite{Bony} is another classical way to overcome the artificial losses of derivatives that occur in the compositions of operators. We denote the paradifferential operators by $\ompi^\psi,$ where $\psi$ is an admissible cut-off (see section \ref{para}).  It is classical that for an operator of order $m,$ 
  $$ \| \ompi^\psi(p(v), \xi) u \|_{H^s} \leq C(|v|_{L^\infty}) \| u \|_{H^{s + m}},$$
  and the norm of the subprincipal symbol in the composition of two para-differential operators depending on $v,$ is in $| v|_{W^{1,\infty}}.$ 
  
  The setting of our interest is semi-classical, in the sense that the operators depend on $\xi$ through $\e \xi.$ It is easy to check that the above bounds and symbolic calculus can be adapted to this setting (section \ref{para}). In a semiclassical setting, subprincipal symbols arising in the compositions come with a prefactor $\e.$ Because the singularity is in $1/\e^2,$ this implies that we need only consider the principal and the subprincipal symbols. The perturbation of the initial data is accordingly assumed to have a semiclassical Sobolev regularity. In particular, it can take the form $\phi(x) e^{i k x/\e},$ with $\phi \in H^s.$ The final asymptotic estimate \eqref{est-th} is formulated in $H^s_\e.$ It implies an estimate in $L^\infty.$ 
  
  We finally mention a technical point, associated with a lack of regularity of the operators involved in the changes of variables described in sections \ref{142} and \ref{143}, caused by the fact that the spectral decomposition of the hyperbolic operator in the (EM) equations becomes singular for small frequencies.  
  
  The wave equation in (Z), that comes up as a geometric optics approximation of the Euler equations, is associated with symbols in $\pm \e |\xi|.$ In particular, these symbols are only bounded at the origin. Because resonances between Klein-Gordon and acoustic modes occur precisely for small frequencies $|\xi | \sim \e,$ a smoothing procedure, or the introduction of a cut-off, would create large error terms. At this point, we make a crucial use of the fact that all the symbols depend on the solution $u^\e,$ only through $\e u^\e,$ and that we work on \emph{perturbation} equations: $u^\e = u_a^\e + \e^k \dot u^\e,$ where $k$ is large enough. Because we need to handle symbols only up to $O(\e^2),$ we can approximate $p(\e u^\e, \e\xi)$ by $p(0,\e \xi) + \e \d_u p(0, \e \xi) \cdot u_a^\e.$ This approximate symbol is easier to handle, for two reasons. First, it depends on $x$ only through $u_a^\e,$ the approximate solution. Second, it has the simple form $p_1(v) p_2(\e \xi).$ In section \ref{limited}, we describe how operators with non-smooth symbols of this form operate in semi-classical Sobolev spaces. Classically, norms of pseudodifferential operators depend on derivatives in $x$ and in $\xi$ of the symbol, and, because $x$ and $\xi$ play somehow symmetric roles, derivatives in $\xi$ can be shifted to derivatives in $x.$ Here we can afford to lose derivatives in $x.$ The eventual energy estimate in $H^s_\e$ involves $\| u_a^\e \|_{H^{s'}},$ with $s' > s.$ This does not harm the proof if the initial datum is assumed to have enough Sobolev regularity.

\end{subsubsection}
\end{subsection}


\begin{subsection}{Background and references}

 The (Z) system was derived from kinetic models by Vladimir Zakharov and his collaborators in the seventies \cite{Zak}.
 
 The initial value problem for the (Z) equations has received much attention. Global existence of smooth solutions in one space dimension (and of weak solutions in two and three space dimensions, for small initial data) was proved by Sulem and Sulem \cite{SS0}. Global existence of smooth solutions in two space dimensions, for small initial data, was proved by Added and Added \cite{AA}. Schochet and Weinstein \cite{SW} and Ozawa and Tsutsumi \cite{OT} showed existence of local in time, smooth Sobolev solutions. Colliander and Bourgain \cite{CB} and Ginibre, Tsutsumi and Velo \cite{GTV} studied critical regularity issues for local solutions. For large initial data, no evidence of singularity in finite time is known in space dimension greater than one.  

 In their book on the Schr\"{o}dinger equation \cite{SS}, Catherine Sulem and Pierre-Louis Sulem show how the Zakharov equation can be formally derived from the Euler-Maxwell equations; the WKB asymptotics of section \ref{3-1} is based on their description, and on discussions with Vladimir Tikhonchuk and Thierry Colin.

 To our knowledge, the results of \cite{T-AA} and \cite{CEGT} were the first results establishing rigorous links between Euler-Maxwell and Zakharov. 
 
 Formal WKB expansions, carried out in \cite{T-AA}, have shown how the weakly nonlinear limit of (EM) fails to describe nonlinear interactions; such a phenomenon had been observed by Joly, Métivier and Rauch in the context of the Maxwell-Bloch equations. Joly, M\'etivier and Rauch's paper \cite{JMR_TMB}, that describes large-amplitude solutions of semilinear systems of Maxwell-Bloch type by means of normal form reductions, is the main inspiration of the present work.

 In \cite{CEGT}, Klein-Gordon-waves systems were formally derived from Euler-Maxwell, and the Zakharov equations were rigorously derived as a high-frequency limit of these Klein-Gordon-waves systems. 
 The stationnary phase arguments of \cite{CEGT}, where solutions were represented, through Fourier analysis, in the form $\e^{-1} \int_0^t e^{i t \Phi/\e} B(t') dt',$ are analogous to the normal form reductions of the present work. The above integrals are bounded if the ratio $B/\Phi$ is bounded, which echoes the transparency condition \eqref{int00}.

 Highly-oscillating, large-amplitude solutions of quasilinear systems were considered by Serre in \cite{Se}, and by Cheverry, Gu\`es and M\'etivier in \cite{CGM}. These papers deal with conservation laws, in particular, non-dispersive systems, unlike the Euler-Maxwell system. In \cite{CG}, Cheverry studies the parabolic relaxation of the instabilities put in evidence in \cite{CGM} and applies his results to the equations of the large-scale motions in the atmosphere.

 In \cite{EG}, Grenier studies a class of singular equations of the form \eqref{em3_0}, with ${\cal A}$ of the form \eqref{dispersif}, and ${\cal B} = 0.$ He proves existence of solutions over time intervals independent of $\e,$ under the assumption that ${\cal A}$ possesses a 'good' symmetrizer, in the sense that no singular terms are created by subprincipal symbols occurring in the symmetrization process. Grenier is naturally led to study operators depending on $x$ through $v(x),$ where $v$ has a Sobolev regularity. He does not assume that the initial data are well-prepared, and studies the high-frequency behaviour of the solutions.

 Lannes recently gave precise bounds for the norms of pseudodifferential operators depending on $x$ through $v(x),$ where $v$ has a Sobolev regularity, and for the norms of commutators of such operators. These questions had previously been considered by Taylor \cite{Tay}, and by Grenier in the article mentioned above. We use a consequence of Lannes' description of the paradifferential remainder (formulated as Proposition \ref{remainder}; it is used in sections \ref{251} and \ref{red3}). 
 
 In the approximate solution $u_a^\e$ to (EM) that is constructed in section \ref{3-1}, the envelope of the electric field is $O(1),$ while the fluctuation of density has size $O(\e),$ in particular, it vanishes in the high-frequency limit. However, in the (Z) system, the fluctuation of density has a finite effect on the electric field. This means that there is a strong coupling between $E$ and $n$ in the Euler-Maxwell equations, an evidence of which is the non-transparent condition \eqref{non-transp}. Such a phenomenon was called "ghost effect" by the Kyoto school of Sone, Aoki and Takata. These authors extensively studied this phenomenon in the context of small Knudsen number analysis of rarefied gases; a good reference is Sone's book \cite{Sone}, and the paper \cite{TA}. As Sone explains in \cite{Sone}, ghosts effects are characteristic of situations where large temperature variations are recorded.  It would be interesting to show a formal similarity between their formal Hilbert expansions (describing the continuum limit) and the WKB expansions of highly nonlinear geometric optics (describing high-frequency limits).  
 
 M. Colin an T. Colin propose a generalization of the Zakharov system in \cite{CC}. Their system consists in four Schr\"{o}dinger equations coupled with quasilinear terms and a wave equation. It describes three-wave interactions, in particular, the generation of a Raman backscattered field. It is an interesting question, to know whether or not the result of this paper could be generalized to their extension of the (Z) system.
 
 We conclude this introduction by mentioning open questions and directions for future work.
 
 It is natural to ask whether or not the result still holds when the initial condition is assumed to be oscillatory, that is, has the form $a^\e(x) e^{i k /\e^2},$ where $a^\e$ is a bounded family in $H^s.$ It is shown in \cite{T-AA} that, if $k \neq 0,$ the limit system is $\mbox{(Z)}_c,$ where $c = \om'(k),$ and $\om$ is a local parameterization on the characteristic variety (see figure \ref{fig-char1}). Linares, Ponce and Saut prove in \cite{LPS} that this system is well-posed in Sobolev spaces \cite{LPS}; Colin and M\'etivier prove in \cite{CM} that it is ill-posed in $L^\infty.$ 
 
 Another interesting direction for future work is to consider the case of large perturbations of the initial data, of the form $\e^{k_0} \phi^\e,$ with $k_0 < 3 + \frac{d}{2}.$ Our guess is that the strong coupling between the electric field and the mean mode of the fluctuation of density would then create strong instabilities in short time.

\end{subsection}
\end{section}



\begin{section}{A class of singular equations} \label{deux}


\begin{subsection}{Symbols} \label{symbols} 

 We consider profiles $u, v, \dots$ depending on $x \in \R^d,$ with values in $\C^n,$ and symbols $p, q, \dots$ depending on $(\e, v, \xi) \in (0,1] \times \C^n \times \R^d,$ or on $\e, x, \xi,$ with values in the $n \times n$ matrices with complex coefficients. 
  
 For $0 < \e \leq 1$ and $s \in \R,$ we let
 $$ \| v \|_{\e,s} := \| (1 + | \e \xi|^2)^{s/2} \hat v(\xi) \|_{L^2(\R^d_\xi)}.$$
 In particular, $ \| \cdot \|_{1,s}$ denotes the classical Sobolev norm. A profile $v$ is said to belong to $H^{s}_\e(\R^d)$ when $\| v \|_{\e,s}$ is finite. The space $H_1^s(\R^d),$ or simply $H^s(\R^d),$ is the classical Sobolev space. 
  Remark that 
 $$ \| h_\e v\|_{1,s} = \| v \|_{\e,s},$$
where 
 $$ h_\e v(x) := \e^{d/2} v(\e x).$$

 For $k \in \N,$ let
 $$ | v |_{k,\infty} := \sum_{0 \leq |\a| \leq k} \sup_{x \in \R^d} |\d^\a_x v|.$$
  A profile $v$ is said to belong to $W^{k, \infty}(\R^d)$ when $| v |_{k,\infty}$ is finite. Let $d_0 > \frac{d}{2}.$ For all $k \in \N,$ the embedding $H^{k + d_0}_\e \hookrightarrow W^{k,\infty}$ holds:
 $$ | v |_{k,\infty} \leq C \e^{-k - d/2} \| v \|_{\e,k + d_0}.$$

 We now define, and somehow adapt to our context, the class of symbols studied by Lannes in \cite{L0} (see also Taylor \cite{Tay} and Grenier \cite{EG}).

  A symbol $p(\e, v, \xi)$ defined in $(0, 1) \times \C^n \times \R^d,$ is said to belong to the class $C^\infty {\cal M}^m,$ $m \in \N,$ when there exists $\e_0 >0$ such that 
 \begin{itemize}
 \item $ p_{|\{|\xi| \leq 1\}} \in C^\infty((0,\e_0) \times \R^d, L^\infty(\{|\xi| \leq 1\}),$ and 
 \item for all $\a, \b,$ there exists a non-decreasing function $C_{\a,\b}$ such that for all $v,$ 
 \begin{equation} \label{s1}
   \sup_{\e \in (0,\e_0)} \; \; \sup_{|\xi| \geq 1/4} \; \langle \xi \rangle^{|\b| - m} | \d_{\e,v}^\a \d_{\xi}^\b p(\e, v, \xi) | \leq C_{\a,\b}(|v|).
 \end{equation}
 \end{itemize}
 In particular, if $p \in C^\infty {\cal M}^m,$ then for all $\a,$ $\d_{\e,v}^\a p \in C^\infty {\cal M}^m.$

 A symbol $p \in C^\infty {\cal M}^m$ is said to be \emph{k-regular at the origin} when 
 $$p_{|\{|\xi| \leq 1\}} \in C^\infty((0,\e_0) \times \R^d, W^{k,\infty}(\{|\xi| \leq 1\}).$$ A symbol is said to be smooth at the origin if it is $k\mbox{-regular}$ for all $k.$ Symbols in $C^\infty {\cal M}^m$ that depend only on $\e, \xi,$ are called Fourier multipliers. 

 If $p \in C^\infty {\cal M}^m$ is evaluated at $v \in H^s_\e,$ $s > \frac{d}{2},$ then Moser's inequality implies that
 $$ \sup_{|\xi| \leq 1} \| p(v(\cdot), \xi) - p(0,\xi) \|_{\e,s} \leq C_s(|v|_{0,\infty}) \| v \|_{\e,s},$$ 
   and for all $\b,$
  
 $$ \sup_{|\xi| \geq 1/4} (1 + |\xi|^2)^{|\b| - m}) \| \d_\xi^\b(p(v(\cdot), \xi) - p(0,\xi)) \|_{\e,s} \leq C_{\b,s}(|v|_{0,\infty}) \| v \|_{\e,s},$$
for some nondecreasing functions $C_s, C_{\b,s}.$ \end{subsection}


\begin{subsection}{Para-differential operators} \label{para}

  The class ${\cal S}^m_{k},$ $m \in \R, k \in \N,$ is defined as the space of symbols $q(\e, x, \xi)$ such that there exists $0 < \e_0 < 1,$ such that, for all $0 < \e < \e_0,$ for all $\a, \b,$ with $|\a| \leq k,$ there exists $C_{\e, \a, \b}(t)$ such that for all $x, \xi,$
 $$ \langle \xi \rangle^{|\a| - m} | \d_x^\a \d_\xi^\b q(\e, x, \xi)| \leq C_{\e,\a,\b}.$$

 With this definition, if $p \in C^\infty {\cal M}^m$ is smooth at the origin, and if $v$ is a profile in $W^{k,\infty}(\R^d)),$ then $q := p(v)$ belongs to ${\cal S}^m_{k}.$

 To $q \in {\cal S}^m_{k},$ one associates the pseudo-differential operator $\ompi_{\e'} (q)$, $0 < \e' \leq 1,$ formally defined by its action as
 $$ (\ompi_{\e'}(q) z)(x) := (2 \pi)^{-d/2} \int_{\R^d} e^{i x \cdot \xi} q(\e, x , \e' \xi) \hat z(t, \xi) d\xi.$$
 With this definition,
  $$
 \ompi_\e(q) := (h_\e)^{-1} \ompi_1(\tilde q) h_\e,
$$
 where $\tilde q(\e,t,x,\xi) := q(\e,t,\e x,\xi).$

 Symbols in ${\cal S}^m_{k}$ are smoothed into paradifferential symbols as follows. Let $\chi_0: \R_+ \to \R$ be a smooth function, such that $0 \leq \chi \leq 1,$ and
 $$ \chi_0(\lam) = 1, \quad \mbox{for $\lam \leq 1.1$}; \quad \chi_0(\lam) = 0, \quad \mbox{for $\lam \geq 1.9$}.$$ 
 For $k \geq 0,$ define $\phi_k: \R^d \to \R$ by
 $$ \phi_k(\xi) := \chi_0(2^{-k} (1  + |\xi|^2)^{1/2}) - \chi_0(2^{-(k-1)} (1  + |\xi|^2)^{1/2} ).$$ With these notations, for all $\xi \in \R^d,$
 $$ 1 = \sum_{k \geq 0} \phi_k(\xi).$$
  Let $\chi: \R^d \to \R$ be a smooth function, $0 \leq \chi \leq 1,$ and such that
 $$ \chi(\eta) = 1, \quad \mbox{for $|\eta| \leq 1.1$}; \quad \chi(\eta) = 0, \quad \mbox{for $|\eta| \geq 1.9$}.$$ 
Let $\psi: \R^d \times \R^d \to \R$ be defined by
 $$\psi(\eta,\xi) = \sum_{k \geq 0} \chi(2^{-k + 3} \eta) \phi_k(\xi).$$
 Then $\psi$ satisfies, 
 $$ \psi(\eta, \xi) = \left\{\begin{array}{cc} 1, & |\eta| \leq  2^{-5} (1  + |\xi|^2)^{1/2},\\  0, & |\eta| \geq 2^{-1} (1  + |\xi|^2)^{1/2}.\end{array}\right.$$
 In particular,
 \begin{equation} \label{tri}
  \psi(\eta,\xi) = 1, \quad \mbox{for all $\xi,$ for all $|\eta| \leq 2^{-5}.$}\end{equation}
  One lets
 $$ \widehat{q^\psi} (\eta, \xi) := \psi(\eta, \xi) \hat q(\eta, \xi).$$ 
 The paradifferential operator associated with $q$ is $$\ompi_{\e'}^\psi(q) := \ompi_{\e'}(q^\psi).$$

 The following proposition describes how well the action of a pseudo-differential operator is approximated by the action of its associated para-differential operator, a classical result in the case of differential symbols and when $\e =1,$ of which Lannes gave an extension to pseudo-differential symbols in \cite{L0}. We check below that the result of Lannes extends to $0 < \e <1;$ in particular, that the action of the para-differential remainder, in $H^s_\e,$ is very small with respect to $\e,$ when $s$ is large.

 $C$ denotes non-decreasing functions, and $d_0$ is a real number such that $[\frac{d}{2}] < d_0 \leq [\frac{d}{2}] + 1.$


 \begin{prop}[remainder] \label{remainder} Let $v \in H^{s}_\e,$ $s > \frac{d}{2},$ and $p \in C^\infty {\cal M}^m$ be smooth at the origin. Then  $$ \ompi_\e(p(v)) =  \ompi_\e^\psi(p(v)) + \e^{s - d/2} \ompi_\e(R_{p(v)}),$$
 and, for all $u \in H^{m + d_0}_\e,$
 \begin{equation} \label{est-r}
 \| \ompi_\e(R_{p(v)}) u \|_{\e,s} \leq C(|v|_{0,\infty}) \| v\|_{\e,s} \| u \|_{\e, m + d_0}.
 \end{equation}
 \end{prop}
 
 Above, $C$ denotes a nondecreasing function, and $d_0$ is a real number such that $[\frac{d}{2}] < d_0 \leq [\frac{d}{2}] + 1.$

\begin{proof} We indicate how \eqref{est-r} follows from Propositions 3.3 and 3.4 of \cite{L0}. Let $q := p(v) - p(0).$ The operation of para-differential smoothing is a convolution in $x,$ and thus $R_q = R_{p(v)}.$ The Fourier transform of the symbol of $R_q$ is 
 $$ \widehat{R_q} = (1 - \psi(\eta,\xi)) \hat q(\eta,\xi).$$
 The point is that because of \eqref{tri}, the above symbol is identically zero for small $\eta.$ Let $\tilde \phi_0 := \chi,$ and for $k \geq 1,$
 $$ \tilde \phi_k(\eta) := \chi(2^{-k} \eta) - \chi(2^{-(k-1)} \eta).$$
Then, for all $\eta \in \R^d,$
 $$ 1 = \sum_{k \geq 0} \tilde \phi_k(\eta).$$
One can write
 $$
 \hat q(\eta,\xi) = \sum_{| k - k' | \geq 3} \tilde \phi_{k'}(\eta) \hat q(\eta,\xi) \phi_k(\xi) + \sum_{| k - k' | < 3} \tilde \phi_{k'}(\eta) \hat q(\eta,\xi) \phi_k(\xi),$$ where the sums run over integers $k, k' \geq 0.$ 
 The first sum in $\hat q$ is further decomposed into $\hat q_1 + \hat q_2,$ where 
 $$ \hat q_1 := \sum_{k \geq 3} \sum_{k' \leq k - 3} \tilde \phi_{k'}(\eta) \hat q(\eta,\xi) \phi_k(\xi), \quad \hat q_2 := \sum_{k' \geq 3} \sum_{k \leq k' - 3} \tilde \phi_{k'}(\eta) \hat q(\eta,\xi) \phi_k(\xi).$$
 Remark that  
 $$ (\psi \hat q)(\eta, \xi) = \hat q_1(\eta, \xi)  + \sum_{k < 3} \chi(2^{-k + 3}\eta) \hat q(\eta, \xi) \phi_k(\xi).$$
 Thus we have
 $$ R_q = q_2 + q_{r,1} + q_{r,2},$$
 where
 \begin{eqnarray} \hat q_{r,1} & := & \sum_{k' \geq 0, |k - k'| < 3} \tilde \phi_{k'}(\eta) \hat q(\eta,\xi) \phi_k(\xi),\nonumber\\
 \hat q_{r,2} &:=  & \chi(\eta) \hat q(\eta,\xi) \chi_0(2^{-3} (1 + |\xi|^2)^{1/2}) - \sum_{k < 3} \chi(2^{-k + 3}\eta) \hat q(\eta,\xi) \phi_k(\xi).\nonumber \end{eqnarray}
The symbols $\hat q_2,$ $\hat q_{r,1},$ $\hat q_{r,2}$ correspond to the symbols $\s_{II},$ $\s_{R,1}$ and $\s_{R,2}$ in \cite{L0}.  We want to bound
 $$ \| \ompi_\e(R_{q}) u \|_{\e,s} = \| \ompi_1 (\tilde q_2 + \tilde q_{r,1} + \tilde q_{r,2}) h_\e u \|_{1,s},$$
 where $\tilde q(x) := q(\e x).$ Propositions 3.3 and 3.4 of \cite{L0} imply that
 \begin{equation} \label{m2} \| \ompi_1(\tilde q_2 + \tilde q_{r,1}) h_\e u \|_{1,s} \leq M(\tilde q) \| u \|_{\e,m + d_0},\end{equation}
 where  $$ M(q) := \sup_{|\g| \leq \g_0} \sup_{\xi \in \R^d} (1 + |\xi|^2)^{|\g| - m} \|\d_\xi^\g q(\cdot, \xi)\|_{1,s},$$
 where $\g_0$ depends only on $d.$ Now, owing to \eqref{tri}, $q$ can be replaced by $(1 - \ompi_1(\tilde \chi)) q$ in the symbol of $R_q,$ the function $\tilde \chi$ being smooth, identically equal to one for $|\eta| \leq 2^{-6},$ and identically equal to zero for $|\eta| \geq 2^{-5}.$ For all $w \in H^s,$ for all $|k| \leq s,$ there holds
 \begin{equation} \label{trick} \| (1 - \ompi_1(\tilde \chi)) w \|_{1,s} \leq C \| \d_x^k w \|_{1,s-k}.\end{equation}
 Applying \eqref{trick} to \eqref{m2}, one finds that the contribution of $\tilde q_2 + \tilde q_{r,1}$ to the operator norm of $R_{q}$ is bounded by
 \begin{equation} \label{bound} C \e^{s - d/2} \sup_{|\g| \leq \g_0} \sup_{\xi \in \R^d} (1 + |\xi|^2)^{|\g| - m} \| \d_\xi^\g(p(v)  - p(0)) \|_{\e,s},
 \end{equation}
 which, in turn, is bounded by $\e^{s- d/2} C(|v|_{0,\infty}) \| v\|_{\e,s}.$ The support of $\hat q_{r,2}$ is included in a ball $|\eta| + |\xi| \leq A.$ Thus we can use Lemma 4.4 of \cite{L0} to estimate the contribution of $\ompi_\e(q_{r,2}).$ Using \eqref{trick} again, one sees that, up to a multiplicative constant, it is also bounded by \eqref{bound} 
 \end{proof}

 The next propositions are classical results of para-differential calculus (see for instance Appendix 2 of \cite{MZ}), formulated in our $C^\infty {\cal M}^m$ framework. In the following statements, $C$ denotes non-decreasing functions, and $d_0$ is a real number such that $[\frac{d}{2}] < d_0 \leq [\frac{d}{2}] + 1.$

 \begin{prop}[action] \label{action} Let $v \in L^\infty$ and $p \in C^\infty {\cal M}^m$ be smooth at the origin. The operator $\ompi_\e^\psi(p(v))$ maps $H^{s+m}_\e$ to $H^{s}_\e,$ for all $s,$ and for all $u \in H^{s+m}_\e,$
  $$
   \| \ompi_\e^\psi(p(v)) u \|_{\e,s} \leq C (|v|_{0,\infty}) \| u \|_{\e,s+m}.
$$ \end{prop}


 \begin{prop}[composition] \label{composition} Let $p_1 \in C^\infty {\cal M}^m_1$ and $p_2 \in C^\infty {\cal M}^m_2$ be smooth at the origin, and let  
 $$ p_1 \sharp p_2  := \sum_{|\a| = 1} \frac{(-i)^\a}{ \a !} \d_\xi^\a p_1 \d_x^\a p_2.$$
 If $p_1, p_2$ are evaluated at $v_1, v_2 \in W^{1,\infty},$ for all $u \in H_\e^{s + 1 - m_1 - m_2},$
 $$
 \| \;  (\ompi_\e^\psi(p_1) \ompi_\e^\psi(p_2) - \ompi_\e^\psi(p_1 p_2)) u \; \|_{\e,s} \leq  \e C (| v_1, v_2|_{1,\infty}) | u \|_{\e, s +1  - m_1 - m_2}.
  $$
 If $p_1, p_2$ are evaluated at $v_1, v_2 \in W^{2,\infty},$ then for all $u \in H_\e^{s + 2 - m_1 - m_2},$
$$\begin{aligned} \| \;  (\ompi_\e^\psi(p_1) \ompi_\e^\psi(p_2) & - \ompi_\e^\psi(p_1 p_2 - \e p_1 \sharp p_2)) u \; \|_{\e,s} \\ & \leq \e^2 C (| v_1, v_2|_{2,\infty}) | u \|_{\e, s + 2 - m_1 - m_2}.
 \end{aligned}$$
 \end{prop}


 \begin{prop}[adjoint] \label{adjoint} Let $v \in W^{2,\infty}$ and $p \in C^\infty {\cal M}^m$ be smooth at the origin. Let $p(v)^*$ denote the complex adjoint of the matrix $p(v),$ and let $\ompie(p(v))^*$ denote the adjoint of $\ompie(p(v))$ in $L^2.$ Let 
 $$ r_*(p) := \sum_{|\a| = 1} \frac{(-i)^\a}{\a!} \d_\xi^\a \d_x^\a p(v)^*.$$
 Then, for all $u \in H^{s + m-2}_\e,$ 
  $$ \begin{aligned}
 \| \; (\ompi_\e^\psi(p(v))^*  - \ompi_\e^\psi(p(v)^*) & - \e \ompie(r_*(p(v))) ) u \; \|_{\e,s}  \\ & \leq  \e^2 C(| v |_{2,\infty}) \| u \|_{\e,s + m - 2}.\end{aligned}
 $$
 \end{prop}

\end{subsection}


\begin{subsection}{Pseudo-differential operators with limited regularity}\label{limited}
  
  We now consider non-smooth symbols that have the simple product structure:
 \begin{equation} \label{structure-produit} p(\e, v, \xi) = p_1(v) p_2(\e, \xi).\end{equation}
If $p$ is a symbol in $C^\infty {\cal M}^m,$ with the structure \eqref{structure-produit}, then, in particular, $p_1$ is a smooth map, and $p_2 \in C^\infty {\cal M}^m.$ Matrix-valued symbols are said to have the structure \eqref{structure-produit} when every entry can be written as a sum of terms of the form \eqref{structure-produit}. 

We will use the following lemma:
 \begin{lem} \label{dav-trick} Let $u, v \in H^s(\R^d),$ $s > \frac{d}{2},$ and assume that $\hat v$ has compact support. Then
$$\| u v \|_{1,s} \leq C \| u \|_{1,s} | v|_{0,\infty},$$
 where $C$ depends only on $s$ and on the space dimension.
 \end{lem}

\begin{proof} This estimate follows easily from a dyadic decomposition. More details can be found in Lemma 3.1 of \cite{L0}, for instance.
\end{proof} 

 In the following propositions, $C$ denotes nondecreasing functions, and $[\frac{d}{2}] < d_0 \leq [\frac{d}{2}] + 1.$

\begin{prop} \label{action-pdo-trick} Let $p \in C^\infty {\cal M}^m,$ of the form \eqref{structure-produit}, and $v\in H^s_\e,$ $s > 1 + \frac{d}{2}.$ For all $u \in H^{s + m}_\e(\R^d),$
   $$
   \| \ompi_\e(p(v)) u \|_{\e,s} \leq C (\| u\|_{\e,s + m} + \| v \|_{1,d_0} \| u \|_{L^2} + \e^{s - d/2} \| v \|_{\e,s} \| u \|_{\e,m + d_0}),
 $$
 where $C$ depends on $|v|_{1,\infty}.$
\end{prop}

\begin{proof} If $p$ does not depend on $v,$ the result obviously holds, and so, changing $p_1$ to $p_1 - p_1(0)$ if necessary, we are reduced to the case $p_1(0) = 0.$

  We use the smooth truncation $\chi$ introduced at the beginning of section \ref{para}.  As $p_2$ is smooth for $|\xi| \geq 1/4,$ the action of $p_1(v) (1 - \ompi_\e(\chi)) \ompi_\e(p_2)$ can be estimated with Propositions \ref{remainder} and \ref{action}. It remains to bound
 \begin{equation} \label{re}  \| \ompi_\e(p_1(v) \chi p_2) u \|_{\e,s} = \| p_1(\tilde v) \ompi_1(\chi p_2) h_\e u \|_{1,s},\end{equation}
 where $\tilde v(x) = v(\e x).$ We can write $p_1( \tilde v)$ as the sum of $\ompi_1(\chi) p_1(\tilde v)$ and $(1 - \ompi_1(\chi)) p_1(\tilde v),$
 and apply the above Lemma. Thus we need to bound
 \begin{equation} \label{fi} |\ompi_1(\chi) p_1(\tilde v)|_{0,\infty} \| \ompi_1(\chi p_2) h_\e u\|_{1,s},\end{equation}
and
 \begin{equation} \label{se} \| (1 - \ompi_1(\chi)) p_1(\tilde v) \|_{1,s} | \ompi(\chi p_2) h_\e u|_{0,\infty}.\end{equation}
The first factor in \eqref{fi} is bounded by 
 \begin{eqnarray} \sup_x \big| \int e^{i \e x \xi} \chi(\e \xi) \widehat{p_1(v)} (\xi) d \xi \big| & \leq & C \| \widehat{p_1(v)} \|_{L^1} \nonumber \\ & \leq & C(|v|_{0,\infty}) \| v \|_{1,d_0}.\nonumber \end{eqnarray}
 With \eqref{trick}, the first factor in \eqref{se} is bounded by $\e^{s - d/2} C(|v|_{0,\infty}) \| v \|_{\e,s}.$
 Finally, the bounds
 $$ \| \ompi_1(\chi p_2) h_\e u\|_{1,s} \leq C \| u \|_{L^2},\quad | \ompi_1(\chi p_2) h_\e u|_{0,\infty} \leq C \| u \|_{\e,d_0},$$
 yield the desired estimate.
\end{proof}

 Next we describe the composition of two operators of the form \eqref{structure-produit}.  Because the composition of two such operators involve remainder terms that do not have the form \eqref{structure-produit}, we need the following notations and lemmas. 
 
 Let $q$ be defined on $\R^d_\xi,$ and $p$ be smooth on $\R_x^d.$ Introduce the notations, 
  $$ Q^\e(\xi,\xi') := q(\e \xi + \xi') - q(\xi'),$$
  and, formally,
  $$ \rho(q,p)_{|(\e,x,\xi)} := (\check Q^\e(\cdot, \xi) * p)(x),$$
  where $\check Q^\e$ denotes the inverse Fourier transform of $Q^\e$ in its first variable, and $*$ is a convolution in $x \in \R^d.$ 
  
  With these notations, $\ompi_\e(\rho(q,p)) u$ is, formally, the inverse Fourier transform of 
  $$ \int (q(\e \xi) - q(\e \xi')) \widehat{p(v)}(\xi - \xi') \hat u(\xi') d\xi'.$$
  
 \begin{lem} \label{le-1} Let $p$ be smooth on $\R^n,$ such that $p(0) = 0,$ and $\chi$ be the smooth truncation introduced at the beginning of section \ref{para}.
  \begin{itemize}
  \item[{\rm (i)}] If $q_{\{ |\xi| \leq 1 \}} \in L^\infty$ and $v \in H^{s + d_0},$ 
  $$ \| \ompi_\e(\rho(\chi q, p(v)) u \|_{\e,s} \leq C(|v|_{0,\infty}) \| v \|_{1, s + d_0} \| u \|_{L^2}$$
  \item[{\rm (ii)}] If $q_{\{ |\xi| \leq 1 \}} \in W^{1,\infty}$ and $v \in H^{s+ d_0 + 1},$  then
  $$ \| \ompi_\e(\rho(\chi q, p(v)) u \|_{\e,s} \leq \e C(|v|_{0,\infty}) \| v \|_{1, s + d_0 + 1} \| u \|_{L^2}. $$
     \item[{\rm (iii)}] If $q \in C^\infty {\cal M}^m$ is a Fourier multiplier, if $v \in H^{s+ d_0 + m},$ then  
   $$ \| \ompi_\e(\rho((1 - \chi) q, p(v))) u \|_{\e,s} \leq \e C(|v|_{0,\infty}) \| v \|_{1, s + m +  d_0} \| u \|_{\e, s + m -1}. $$
  \end{itemize}
   \end{lem}
   
 The third estimate is not tame, but it will be sufficient for our purposes.

 \begin{proof} (i) Let $q_0 := \chi q,$ and $w_0 := \ompi_\e(\rho^\e(q_0, p(v)) u.$
  The $H^s_\e$ norm of $w_0$ is equal to the $L^2_\xi$ norm of
 \begin{equation} \label{convo} (1 + |\e \xi|^2)^{s/2} \int (q_0(\e \xi') - q_0(\e \xi)) \widehat{p(v)}(\xi - \xi') \hat u(\xi') d\xi',\end{equation}
 With Peetre's inequality and because $q_0$ is compactly supported, the $L^2_\xi$ norm of \eqref{convo} is bounded by
 $$ |q_0|_{0,\infty} \; \big\| \int  (1 + |\e ( \xi - \xi')|^2)^{s/2} | \widehat{p(v)}(\xi - \xi')| |\hat u(\xi')| d\xi' \big\|_{L^2_\xi}.$$
 The above integral is a convolution. Thus 
 \begin{eqnarray}
  \| w_0 \|_{\e,s} & \leq & C \| \, {\cal F}^{-1} | (1 + |\e \xi|^2)^{s/2} \widehat{p(v)}| \,\, {\cal F}^{-1} | \hat u | \, \|_{L^2_x} \nonumber \\  \nonumber & \leq & C | \,{\cal F}^{-1} | (1 + |\e \xi|^2)^{s/2} \widehat{p(v)}| \, |_{0,\infty} \, \| u \|_{L^2_x}.\end{eqnarray}
 Now the $L^\infty$ norm of ${\cal F}^{-1} (1 + |\e \xi|^2)^{s/2} | \widehat{p(v)}|$ is bounded by the $L^1$ norm of $(1 + |\e \xi|^2)^{s/2} \widehat{p(v)},$ which in turn is bounded by $C(|v |_{0,\infty}) \| v\|_{1,s + d_0}.$ This yields the desried estimate. 
 
 (ii) If $\d_x^\g q_0$ is bounded, then $w$ is the inverse Fourier transform of 
 $$ \e \sum_{|\gamma| = 1} \int \int_0^1 \d_\xi^\gamma q_0(\e \xi + \e t(\xi' - \xi)) \, dt \; \widehat{\d_x^\gamma p(v)}(\xi - \xi')  \hat u(\xi')  d\xi',$$
 and the bounds that led to (i) are easily adapted to obtain (ii). 

 (iii) Let $q_1 := (1 - \chi) q,$ and $ w_1 = \ompi_\e(\rho((1 - \chi) q, p(v))) u.$ The $H^s_\e$ norm of $w_1$ is equal to the $L^2_\xi$ norm of 
  $$ \e \sum_{|\g| = 1} (1 + |\e \xi|^2)^{s/2} \int \int_0^1 \d_\xi^\gamma q_1(\e \xi + \e t(\xi' - \xi)) \, dt \; \widehat{p(v)}(\xi - \xi') \hat u(\xi') d\xi'.$$
  There exists $C$ such that, for all $\e, \xi, \xi',$
   $$ | \d_\xi^\gamma q_1(\e \xi + \e t(\xi' - \xi)) | \leq C ((1 + |\e \xi'|^2)^{(m-1)/2} + (1 + |\e (\xi - \xi')|^2)^{(m-1)/2}).$$
  The rest of the proof of (iii) is similar to the proof of (i).

\end{proof}

 Given two symbols of the form \eqref{structure-produit}:
 \begin{equation} \label{8fev} 
  p(\e,v, \xi) = p_1(v) p_2(\e,\xi), \quad q(\e,v,\xi) = q_1(v) q_2(\e,\xi),
  \end{equation}
  let
   \begin{equation} \label{8fev2} {\tt r}(p,q) :=  q_1 \rho(\chi q_2, p_1) p_2.\end{equation}
    
 \begin{lem} \label{le-2} Let $p, q$ of the form \eqref{8fev}, where $p_1, q_1$ are smooth and vanish at $v =0,$ $\chi q_2 \in L^\infty,$ and $p_2 \in C^\infty {\cal M}^{m}.$ Let ${\tt r} := {\tt r}(p, q).$ If $p_1$ and $q_1$ are estimated at $v \in H^{s + d_0 + 1},$ for all $u \in H^{m_{}}_{\e},$
 $$ \| \ompi_\e({\tt r}) u \|_{\e,s} \leq C( \| v \|_{{1,s + d_{0}}}) \| u \|_{\e,m}.$$
  and, if $a \in C^\infty {\cal M}^1$ is a Fourier multiplier, for all $u \in H^{m}_{\e},$
$$\| (\ompi_\e(a) \ompi_\e({\tt r}) - \ompi_\e(a {\tt r})) u \|_{\e,s} \leq  \e C( \| v \|_{1, s + d_0 + 1}) \| u \|_{\e,m}. $$   \end{lem}

 \begin{proof} Let $z := \ompi_\e(p_2) u,$ $p_1 = p_1(v),$ and $q_1 = q_1(v).$ The Fourier transform of $\ompi_\e({\tt r}) u$ is a convolution in $\xi,$
  $$ \int \widehat{q_1}(\xi - \xi') ((\chi q_2)(\e \xi'') - (\chi q_2)(\e \xi'))  \widehat{p_1}(\xi' - \xi'') \hat z(\xi'') d\xi'' d\xi',$$
   and, because $\chi q_2 \in L^\infty,$ the first estimate in the lemma is obtained in the same way as Lemma \ref{le-1} (i). 
   
   Let $w:= (\ompi_\e(a) \ompi_\e({\tt r}) - \ompi_\e(a {\tt r})) u.$ The Fourier transform of $w$ is the sum, over $|\g| = 1,$ of 
     $$ \e \int a_\g(\e,\e \xi,\e \xi') \widehat{\d_x^\g q_1}(\xi - \xi') ((\chi q_2)(\e \xi'') - (\chi q_2)(\e \xi'))  \widehat{p_1}(\xi' - \xi'') \hat z(\xi'') d\xi'' d\xi',$$ 
     and
   $$ \e \int a_\g(\e,\e \xi',\e \xi'') \widehat{q_1}(\xi - \xi') ((\chi q_2)(\e \xi'') - (\chi q_2)(\e \xi'))  \widehat{\d_x^\g p_1}(\xi' - \xi'') \hat z(\xi'') d\xi'' d\xi',$$   
 where $a_\g(\e,\eta,\eta') := \int_0^1 \d_\xi^\g a(\e,\eta + t(\eta' - \eta)) dt.$ Again, because $a_\g, \chi q_2 \in L^\infty,$ these convolutions can be bounded in the same way as in the proof of Lemma \ref{le-1}, to yield the second estimate. 
  \end{proof}

  We can now state a proposition that describes the composition of two symbols of the form \eqref{structure-produit}. 
  
  \begin{prop} \label{comp-trick-pdo-BIS} Given $p \in C^\infty {\cal M}^{m_1}$ and $q \in C^\infty {\cal M}^{m_2},$
  \begin{itemize}
  \item[{\rm (i)}] if $p$ and $q$ have the form \eqref{8fev}, if $p$ is 1-regular at the origin, if $p$ and $q$ are estimated at $v \in H^{s + d_0 + 1},$ for all $u \in H^{s + m_1 + m_2 -1}_\e,$
$$ \begin{aligned}
  \|  \; [ \ompi_\e(p), \ompi_\e(q) ] u  & - \ompi_\e([p, q] + {\tt r}(p,q)) u \; \|_{\e,s} \\ & \leq \e C( \| v \|_{1,s + d_{0} +1}) \| u \|_{\e,s + m_1 + m_2 -1}; \end{aligned}
  $$
   \item[{\rm (ii)}] if $p$ is a Fourier multiplier and is 2-regular at the origin, if $q$ has the form \eqref{structure-produit} and is estimated at $v \in H^{s + d_0 + 2},$ for all $u \in H^{s + m_1 + m_2 - 2}_\e,$
  $$\begin{aligned}  
 \| \;  [\ompi_\e(p), \ompi_\e(q)] u  & - \ompi_\e( [p, q] + \e (p \sharp q)) u \; \|_{\e,s} \\ & \leq \e^2 C( \| v \|_{1,s + d_{0} + 2 }) \| u \|_{\e,s + m_1 + m_2 - 2}. \end{aligned}
 $$
     \end{itemize}
 \end{prop}
 
 \begin{proof} (i) We compute 
  \begin{eqnarray} \ompi_\e(p) \ompi_\e(q) u & = &  p_1 \ompi_\e(p_2) (q_1 \ompi_\e(p_2) u) \nonumber \\
   & = & p_1 ( \ompi_\e(p_2 q_1) - \ompi_\e(\rho(p_2, q_1))) \ompi_\e(p_2) u. \nonumber  \end{eqnarray}
  Thus 
  $$ \begin{aligned} {} [\ompi_\e(p), \ompi_\e(q)] & - \ompi_\e([p, q] + {\tt r}(p,q)) \\ & = p_1 \ompi_\e(\rho(p_2, q_1)) \ompi_\e(q_2) - q_1 \ompi_\e(\rho((1 - \chi) q_2, p_1) \ompi_\e(p_2), \end{aligned}$$
  and the action of the right-hand side on $u$ is estimated with Lemma \ref{le-1} (ii) and (iii). 
  
  (ii) Let $w := ([\ompi_\e(p), \ompi_\e(q)]  - \ompi_\e[p, q] + \e \ompi_\e(p \sharp q))u.$ The Fourier transform of $w$ is 
   $$ \e^2 \sum_{|\g|  = 2} \int p_\g(\e \xi, \e \xi') \widehat{\d_x^\g q_1(v)}(\xi - \xi') q_2(\xi') \hat u(\xi') d\xi' ,$$   where $p_\g(\e \xi, \e \xi') := \int_0^1 (1 - t) \d_\xi^\g p(\e \xi' + \e t (\xi - \xi')) dt,$ and we obtain the second estimate as above.     
   \end{proof}

  \end{subsection}


\begin{subsection}{Assumptions and results} \label{AR}

 We use the profile and symbol spaces introduced in section \ref{symbols}. Let ${\cal A}, {\cal B}, {\cal G}^\e$ such that
 \begin{itemize}
 \item ${\cal A}$ is a smooth symbol in $C^\infty {\cal M}^1;$
 \item ${\cal B}$ is a bilinear map $\R^n \times \R^n \to \R^n;$
 \item ${\cal G}^\e$ is a family of smooth maps: $\R^n \to \R^n,$ such that ${\cal G}^\e (0) = 0.$
 \end{itemize}

Using the notations of section \ref{symbols}, we denote by $\ompi_\e({\cal A}(\e, \e v))$ the semiclassical pseudo-differential operator with symbol ${\cal A}(\e, \e v, \xi).$ We write the Taylor expansion of ${\cal A}$ in $v$ as follows:
 $$ {\cal A}(\e, \e v) = {\cal A}^{(0)}(\e,0) + \e {\cal A}^{(1)}(v) + \e^2 {\cal A}^{(2)}(\e, v),$$
 where ${\cal A}^{(j)} \in C^\infty {\cal M}^1,$ for $j= 0, 1, 2,$ and where ${\cal A}^{(1)} (v) := \d_v {\cal A}(0,0) \cdot v$ is linear in $v$ and has the form \eqref{structure-produit}. 

 We study the initial value problem
 \begin{equation} \label{j_0} \left\{ \begin{aligned}\d_t u^\e + \frac{1}{\e^2 } \ompi_\e({\cal A}(\e,\e u^\e)) u^\e  & = \frac{1}{\e} {\cal B}(u^\e, u^\e) + {\cal G}^\e(u^\e), \\
  u^\e(0, x) & = a^{\e}(x) + \e^{k_0} \varphi^\e(x), \end{aligned}\right.
 \end{equation}
 where the initial datum $a^{\e}$ belongs to $H^\s(\R^d),$ for some Sobolev index $\s,$ much larger than $\frac{d}{2}$ (we will actually need $\s > 6 + d$), uniformly with respect to $\e:$
 $$
  \sup_{0 < \e <\e_0} \| a^{\e} \|_{1,\s} < \infty.
$$
  The perturbation $\e^{k_0} \varphi^\e$ is such that 
 \begin{equation} \label{k_0}  k_0 > 3 + \frac{d}{2},\end{equation}
 and $\varphi^\e$ belongs to $H^s_\e(\R^d),$ uniformly with respect to $\e:$
 $$
  \sup_{0 < \e < \e_0} \| \varphi^\e \|_{\e,s} < \infty,
 $$ for some large Sobolev index $s,$ smaller than $\s.$ 

Our first assumption is a hyperbolicity assumption that implies in particular the local well-posedness of the initial value problem (\ref{j_0}).

 \begin{ass}[hyperbolicity]  \label{ass_1} For all $\e, v, \xi,$ the matrix ${\cal A}(\e, v, \xi)$ is hermitian. Let
 \begin{equation} \label{spectral-dec} {\cal A} =  \sum_{1 \leq j \leq n_0} i \lam_j \Pi_j + \sum_{n_0 + 1 \leq k \leq n} i \lam_k \Pi_k,\end{equation}
 be its spectral decomposition, where the eigenvalues $\lam_j, \lam_k$ are real and the eigenprojectors $\Pi_j, \Pi_k$ are orthogonal. We assume:
 \begin{itemize}
 \item[{\rm (i)}] the eigenvalues can be ordered as follows: for all $\e,v,\xi,$
  $$ \sup_{n_0 + 1 \leq k} | \lam_k(\e,v,\xi) | < \sup_{j \leq n_0}  | \lam_j(\e,v,\xi)|;$$
 \item[{\rm (ii)}] for all $1 \leq m \leq n,$ $\lam_m(\e,0,\xi) \in C^\infty {\cal M}^1,$ and $\Pi_m(\e,0,\xi) \in C^\infty {\cal M}^0;$ 
 \item[{\rm (iii)}] for all $n_0 + 1 \leq m \leq n,$ $\lam_m(0,0,\xi) = 0,$ for all $\xi.$ 
\end{itemize}
 \end{ass}

 In Assumption \ref{ass_1}, we do not assume that the eigenvector decomposition is not singular, as we want to handle the case of the (EM) equations, whose eigenvector decomposition does become singular for small frequencies (see figure \ref{fig-char1} and section \ref{3-2-2}).

 Let
 \begin{equation} \label{def-proj}
 \Pi_0 := \sum_{1 \leq k \leq n_0} \Pi_k, \qquad \Pi_s := \sum_{n_0 + 1 \leq j \leq n} \Pi_j.
 \end{equation}
In addition to Assumption \ref{ass_1}, we will assume that  
 \begin{equation} \label{reg-proj-totaux} \Pi_0, \Pi_s (\e, v ,\xi) \in C^\infty {\cal M}^0.
 \end{equation}
 When ${\cal A}$ depends analytically on $\e, v, \xi,$ \eqref{reg-proj-totaux} follows from standard considerations, as detailed in section \ref{application}.

 In reference to the Euler-Maxwell equations (see section \ref{application}), the eigenvalues $\lam_j$ for $1 \leq j \leq n_0$ are called Klein-Gordon modes, while the eigenvalues $\lam_k,$ for $n_0 + 1 \leq k \leq n,$ are called acoustic modes. Condition (i) in Assumption \ref{ass_1} states that the acoustic modes do not cross the Klein-Gordon modes. Condition (ii) is a regularity assumption, and condition (iii) amounts to say that the acoustic velocities are $O(\e),$ a consequence of \eqref{theta_i}. 

  \begin{ass}[approximate solution] \label{ass_5} For all $l_0 <\s - 2 - \frac{d}{2},$ there exists $t^{*}(l_0) > 0,$ independent of $\e,$ and a family of profiles $u_a^\e,$ such that
  \begin{equation} \label{j_0-app} \left\{ \begin{aligned}\d_t u_a^\e + \frac{1}{\e^2 } \ompi_\e({\cal A}(\e,\e u_a^\e)) u_a^\e  & = \frac{1}{\e} {\cal B}(u_a^\e, u_a^\e) + {\cal G}^\e(u_a^\e) + \e^{l_0} R_a^\e, \\
  u_a^\e(0, x) & = a^\e(x), \end{aligned}\right.\end{equation}
 where $R_a^\e$ is uniformly bounded with respect to $\e,$
 $$ \sup_{0 < \e < \e_0} \sup_{0 \leq t \leq t^*(l_0)} \| R_a^\e (t) \|_{1,\s - l_0 - 2} < \infty.$$
 There exists a characteristic frequency $\om \neq 0,$ a finite set ${\cal R}^* \subset \Z$ and profiles $\{u^\e_{a,p} \}_{p \in {\cal R}^*}$ and $v_a^\e,$ such that $u_a^\e$ decomposes as
 $$ u_a^\e(t, x) = \sum_{p \in {\cal R}^*} e^{i p \om t/ \e^2} u_{a,p}^\e(t, x) + \e v_a^\e(t,x),$$
 with the uniform bounds
 \begin{eqnarray}\nonumber 
  \sup_{0 < \e < \e_0} \sup_{0 \leq t \leq t^*(l_0)} ( \| u_{a,p}^\e (t) \|_{1,\s} + \| \d_t u^\e_{a,p} (t) \|_{1,\s-2}) & <  & \infty,\\
  \nonumber 
  \sup_{0 < \e < \e_0} \sup_{0 \leq t \leq t^*(l_0)} (\| v_a^\e (t) \|_{1,\s} + \| \e^2 \d_t v_a^\e (t) \|_{1,\s-1}) & < &  \infty.
\end{eqnarray}
 \end{ass}

 For the Euler-Maxwell system of equations, such an approximate solution is explicitly contructed in section \ref{application}, under an assumption of well-preparedness for the initial datum. In the following assumption, ${\cal R}^*$ refers to the set of characteristic harmonics introduced in Assumption \ref{ass_5}. 
 
\begin{ass}[resonances] \label{ass_3} There exists $0 < c_l < c_m < C_m,$ such that the resonance equations in $\xi \in \R^d$ and $p, p' \in {\cal R}^*,$
 \begin{eqnarray} \label{res-eq} \Phi_{j,k,p}(\e) & := & \lam_j(\e,0,\xi) - \lam_k(\e,0,\xi) + p \om = 0, \\  \label{res-eq-2} \Psi_{j,p,p'} & := &  \lam_j(0,0,\xi) - (p + p') \om = 0,
 \end{eqnarray}
  are such that,
 \begin{itemize} \item[{\rm (0-0)}] $1 \leq j, k \leq n_0:$ the solutions $\xi, p$ of \eqref{res-eq} at $\e = 0,$ are located in the interval $c_m \leq |\xi| \leq C_m;$ outside this interval, $\Phi_{j,k,p}(0)$ is bounded away from $0,$ uniformly in $\xi,$ for all $p \in {\cal R}^*;$
 \item[{\rm (0-s)}] $1 \leq j \leq n_0,$ $n_0 + 1 \leq k \leq n:$ for $\e$ small enough, the solutions $\xi, p$ of \eqref{res-eq} are located in the interval $|\xi | \leq c_l;$ outside this interval, $\Phi_{j,k,p}(\e)$ is bounded away from $0,$ uniformly in $\xi$ and $\e,$ for all $p \in {\cal R}^*;$ 
 \item[{\rm (0-0-s)}] $1 \leq j \leq n_0:$ for all $p, p' \in {\cal R}^*,$ $\Psi_{j,p,p'}$ is bounded away from $0,$ uniformly in $|\xi| \in [0, c_m].$
 \end{itemize}
\end{ass}

 Next we state the assumptions that describe the interaction coefficients at the resonances. 

Introduce first the notations:
 \begin{eqnarray} 
   {\cal B}(u_a^\e)  & := & {\cal B}(u_a^\e, \cdot) + {\cal B}( \cdot ,u_a^\e) - {\cal A}^{(1)}(\cdot) u_a^\e, \label{b} \\  {\cal D}(u_a^\e)  & := & ({\cal G}^0)'(u_a^\e) - \d_v {\cal A}^{(2)}(0, u_a^\e) (u_a^\e, \cdot),\label{d} \end{eqnarray}
 where $u_a^\e$ is the approximate solution given by Assumption \ref{ass_5}, and, for all $z,$ $\d_v {\cal A}^{(2)}(0, u_a^\e) (u_a^\e, z)$ is the linear term in $z$ in ${\cal A}^{(2)}(0, u_a^\e + z) u_a^\e - {\cal A}^{(2)}(0,u_a^\e) u_a^\e.$ The symbols ${\cal B}$ and ${\cal D}$ depend on $u_a^\e$ and of its $\e\mbox{-derivatives};$ both belong to $C^\infty {\cal M}^0.$  
 
  We assume that for some $j \leq n_0 < k,$ there exists $p,$ $\xi_0$ and $\eta > 0,$ such that $\Phi_{j,k,p}(\xi_0) = 0,$
 and \begin{equation} \label{non-transp}
 | \Pi_j(0,0) {\cal B}(u_{a,p}^\e) \Pi_k(0,0) | > \eta, \quad \mbox{uniformly in $\xi \sim \xi_0,$}\end{equation}
 Inequality \eqref{non-transp} means that the interaction coefficient $\Pi_j {\cal B} \Pi_k$ is \emph{not transparent} for  resonances between Klein-Gordon and acoustic modes.
 
 Let
  \begin{equation} \label{rho}
 \rho(\e, u_a^\e) := (\Pi_s {\cal A}^{(0)}) \sharp \Pi_0 + (\Pi_s \sharp {\cal A}^{(1)}) \Pi_0 + \Pi_s {\cal A}^{(0)} (\Pi_0 \sharp \Pi_0) - (\Pi_s \sharp {\cal B}) \Pi_0,
 \end{equation}
where the projectors $\Pi_0, \Pi_s$ are evaluated at $(\e, \e u_a^\e),$ and ${\cal A}^{(1)}$ is evaluated at $u_a^\e.$ Let 
\begin{equation} 
 \label{tt-b}
  B^{\tt r}(\e, u_a^\e) :=  \Pi_s ({\cal B} + \d_u \Pi_s \cdot (\e^2 \d_t u_a^\e) + \e \rho ) \Pi_0, \end{equation}
  where $\rho$ is evaluated at $(\e,u_a^\e),$ ${\cal B}$ is evaluated at $u_a^\e,$ and the projectors and their derivatives are evaluated at $(\e, \e u_a^\e).$ 
 Remark that $B^{\tt r}(\e,0) = 0,$ and 
   $$ B^{\tt r}(\e, u_a^\e) = \d_u B^{\tt r}(\e,0) \cdot u_a^\e + \e \d_u^2 B^{\tt r}(\e,0) \cdot (u_a^\e, u_a^\e) + O(\e^2),$$
 The linear term $\d_u B^{\tt r}(\e,0) \cdot u_a^\e$ is the crucial interaction coefficient.  The following transparency assumption states that it is sufficiently small at the resonances between Klein-Gordon and acoustic modes.   

 \begin{ass}[transparency] \label{ass_6} There exists $\e_0 > 0$ and $C > 0$ such that, for all $0 < \e < \e_0,$ for $j \leq n_0 < k,$ for all $p \in {\cal R}^*,$
 \begin{equation} \label{interaction-d}
  | \Pi_k(\e,0) {\cal D}(u_a^\e) \Pi_j(\e,0) | \leq C \e,
 \end{equation}
 and 
  \begin{equation} \label{interaction0} | \Pi_k(\e,0) (\d_u B^{\tt r}(\e, 0) \cdot u_a^\e) \Pi_j(\e,0) | \leq C (\e^2 + |\Phi_{j,k,p}(\e)|),
 \end{equation}
 uniformly in $|\xi| \leq c_l,$ $x \in \R^d$ and $t \in [0, t^*(l_0)).$
  \end{ass}

\medskip

 Introduce finally 
 $$ E := \left(\begin{array}{cc} \Pi_0 {\cal B} \Pi_0 & 0 \\ 0 &  \Pi_s {\cal B} \Pi_s\end{array}\right), \quad i A := \left(\begin{array}{cc} \Pi_0 {\cal A} & 0 \\ 0 &  \Pi_s {\cal A}\end{array}\right),$$
 where ${\cal A}$ and the projectors are seen as symbols depending on $(\e, \e u),$ and ${\cal B}$ as a symbol depending on $(\e, u).$

\begin{ass}[symmetrizability] \label{ass_sym} There exists $S,$ a smooth Fourier multiplier in $C^\infty {\cal M}^0,$ 
such that
 $$ \frac{1}{\gamma} \| u \|^2_{\e,s} \leq (\ompi_\e(S) u, u)_{\e,s} \leq \gamma \| u \|^2_{\e,s},$$
 for all $u \in H^s_\e$ and for some $\gamma > 0,$ and 
 $$
 \frac{1}{\e} (S E + (S E)^*) \in C^\infty {\cal M}^0, \quad  \frac{1}{\e^2} (i S A   + (i S A)^*) \in C^\infty {\cal M}^0.
 $$ \end{ass}

\medskip

In the following theorem, $u_a^\e$ is the approximate solution at order $l_0,$ for some $3 + k_0 \leq l_0 \leq \s - 2 - \frac{d}{2},$ whose existence is guaranteed by Assumption \ref{ass_5}, $t^* = t^*(l_0)$ is its maximal existence time, independent of $\e,$ and $s$ is a Sobolev index, such that $1 + \frac{d}{2} < s < \s - l_0 - 2.$  

\begin{theo} \label{th1}
 Under Assumptions \emph{\ref{ass_1}, \ref{ass_5}, \ref{ass_3}, \ref{ass_6}} and \emph{\ref{ass_sym}}, there exists a unique solution $u^\e \in C^0([0, t_0], H^s_\e(\R^d))$ to the initial value problem \emph{(\ref{j_0})}, for all $0 \leq t_0 < t^{*};$ there exists $C > 0$ and $\e_0 > 0,$ such that, for all $0 < \e < \e_0,$ for all $0 \leq t_0 < t^{*},$
 \begin{equation} \label{est-th} \sup_{0 \leq t \leq t_0} \| (u^\e - u_a^\e)(t) \|_{\e,s} \leq C \e^{k_0-1}.\end{equation}
 \end{theo}
 In particular,
 $$
  \sup_{0 \leq t \leq t_0} | (u^\e - u_a^\e)(t) |_{0,\infty} \leq C \e^{k_0-1-d/2}.$$ In the error estimate (\ref{est-th}), $C$ depends on a Sobolev norm of the initial data and on $t_0.$ 

\end{subsection}


\begin{subsection}{Proof of theorem \ref{th1}} \label{pr}
 
  In the proof below, we often drop the epsilons as we write $u$ for $u^\e,$ $u_a$ for $u_a^\e,$ etc. We use the notations and results of section \ref{symbols} to describe symbols and operators.

  We start with (\ref{j_0}). Let $l_0 \geq k_0 + 3,$ and let $u_a$ be the approximate solution at order $l_0$ given by Assumption \ref{ass_5} An existence time for $u_a$ is $t^* >0,$ independent of $\e.$  
  
\begin{subsubsection}{The perturbation equations} \label{251}

  The exact solution $u$ is sought as a perturbation of $u_a:$
 \begin{equation} \label{first-move} u = u_a + \e^{k_0 -1} \dot u,\end{equation}
The symbol ${\cal A}^{(0)}$ does not depend on $v.$ This implies, with Proposition \ref{remainder},
 $$ \ompi_\e({\cal A}(\e, \e u)) = \ompie({\cal A}(\e, \e u)) + \e^{s + 1 - d/2} \ompi_\e(R_{{\cal A}^{(1)}(u)} + \e R_{{\cal A}^{(2)}(\e,u)}).$$
The perturbation equations are 
  \begin{equation} \label{perturb} \left\{ \begin{aligned} \d_t \dot u + \frac{1}{\e^2} \ompie({\cal A}(\e, \e u )) \dot u & = \frac{1}{\e} \ompie({\cal B}) \dot u +  \ompi_\e({\cal D}) \dot u  + \e R^\e, \\ \dot u(0, x) & = \e \phi^\e(x). \end{aligned} \right. \end{equation}
where ${\cal B}$ and ${\cal D}$ are given by \eqref{b} and \eqref{d}, 
 and where 
 $$\begin{aligned} R^\e & :=  \e^{k_0 - 2} \big( {\cal B}(\dot u, \dot u) + \int_0^1 (1-t) ({\cal G}^\e)''(u_a + \e^{k_0 -1} t \dot u) \cdot (\dot u, \dot u) dt \big) \\ & - \e^{k_0-1} \int_0^1 (1-t) \d_v^2 {\cal A}^{(2)}(u_a + \e^k t \dot u) \cdot (u_a, \dot u, \dot u) dt \\ & - \e^{s-1-d/2} \ompi_\e(R_{{\cal A}^{(1)}(u)} + \e R_{{\cal A}^{(2)}(\e,u)}) - \e^{l_0 - k_0 - 2} R_a^\e. \end{aligned}$$
 Under Assumption \ref{ass_1}, standard hyperbolic theory provides the existence of a unique solution $\dot u$ to \eqref{perturb} over a small time interval $[0, t_*(\e)],$ with the uniform estimate
 \begin{equation} \label{borne-u} \sup_{0 < \e < \e_0} \sup_{0 \leq t \leq t_*(\e)} \| \dot u (t) \|_{\e,s} \leq \delta.\end{equation}

The term $R^\e$ is a remainder, in the sense that its $H^s_\e$ norm can be bounded in terms of $\delta,$ uniformly in $\e.$ The $H^s_\e$ norm of the terms in the first line in the definition of $R^\e$ is indeed bounded by   
  $\e^{k_0 - 2 - d/2} C \| \dot u \|_{\e,d_0} \| \dot u \|_{\e,s};$ the terms in the second line $R^\e$ are bounded by $$\e^{k_0 - 1 - d} C ( \| u_a \|_{\e,s+1} \| \dot u \|_{\e,d_0}^2 + |u_a|_{1,\infty} \| \dot u \|_{\e,d_0} \| \dot u \|_{\e,s}),$$ and the last terms are bounded by 
   $$ \e^{s - 1 - d/2} C \| \dot u \|_{\e,1 + d_0} \| u \|_{\e,s}  + \e^{l_0 - k_0 - 2} \| R_a^\e \|_{\e,s}.$$
 In all these estimates, $C$ depends on $|u|_{0,\infty}.$
 
  With the estimates for $u_a$ given in Assumption \ref{ass_5}, the bound for $\dot u$ given in \eqref{borne-u}, and the form of the equation \eqref{perturb}, for $\a + |\b| \leq 2,$ 
  \begin{equation} \label{est-bar-u1}
 \sup_{0 \leq t \leq t_*(\e)} | (\e^2 \d_t)^\a \d_x^\b u |_{0, \infty} \leq c_a + \e^{k_0 - 3 - d/2} \delta,\end{equation}
 where $c_a$ does not depend on $\e.$ The size of the perturbation of the initial data in \eqref{j_0}, namely $O(\e^{k_0})$ in $H^s_\e,$ where $k_0$ satisfies \eqref{k_0}, was chosen in order that the estimate \eqref{est-bar-u1} be uniform in $\e.$

 With these notations, the above estimates give, for $\e_0$ small enough, and $s > 1 + \frac{d}{2},$
 $$ \| R^\e \|_{\e,s} \leq \tilde C \| \dot u \|_{\e,s} + \e^{l_0 - k_0 - 2} \| R_a^\e \|_{\e,s},$$
 where $\tilde C$ is a nondecreasing function of $\e_0^{k_0 - 3 - d/2} \delta,$ $\| u_a\|_{\e,s + 1},$ $s$ and $d.$

 We generically denote by $R_{(0)}$ any pseudo- or para-differential operator, possibly depending on the solution $u,$ such that, for all $z \in H^{s}_\e(\R^d),$
 \begin{equation} \label{r-generic}
 \| R_{(0)} z \|_{\e,s} \leq \tilde C \| z \|_{\e,s},
 \end{equation}
uniformly in $t \in [0, t_*(\e)],$
where $\tilde C$ is nondecreasing, and 
 \begin{equation} \label{tilde-c} \tilde C = \tilde C(\e_0^{k_0 - 3 - d/2} \delta, \| u_a \|_{1, s + d_0 + 2}, s, d),\end{equation}
 and where $[\frac{d}{2}] < d_0 \leq [\frac{d}{2}] + 1.$

  
  We denote by $O(\e^k)$ symbols associated with pseudo- or para-differential operators of the form $\e^k R_{(0)}.$ 
 
 In the next section, we are led to study compositions of para-differential symbols of the form $p_j(u_{\a,\b}),$ where $p_j \in C^\infty {\cal M}^{m_j}$ is smooth, and $u_{\a,\b} = (\e^2 \d_t)^\a \d_x^\b u.$ It follows from Propositions \ref{action} and \ref{composition} that
 \begin{itemize}
  \item If $m_j \leq 0,$ and $\a + |\b| \leq 2,$ then $\ompie(p_j(u_{\a,\b})) = R_{(0)}.$ 
  \item If $m_1 + m_2 \leq 1,$ if $p_1 = p_1(u_{\a,\b}),$ and $p_2 = p_2(u_{\a',\b'}),$ with $\a + \a' + |\b| + |\b'| \leq 1,$ then
     $$ \ompie(p_1) \ompie(p_2) - \ompie(p_1 p_2) = \e R_{(0)},$$
  \item if $m_1 + m_2 \leq 2,$ if $p_1 = p_1(u)$ and $p_2 = p_2(u),$ then
   $$ \ompie(p_1)\ompie(p_2) - \ompie(p_1 p_2) - \e \ompie(p_1 \sharp p_2) = \e^2 R_{(0)}.$$
  \end{itemize}

\end{subsubsection}

\begin{subsubsection}{Projection and rescaling} \label{resc}
  
 Let $\chi_\e: \R^d \to \R$ be a smooth function such that $\e \leq \chi_\e \leq 1,$ and such that $\chi_\e$ is identically equal to 1 for $|\xi| \leq c_0,$ and identically equal to $\e$ for $|\xi| \geq c_1,$ where $0 < c_0 < c_l < c_1 < c_m,$ with the notations of Assumption \ref{ass_3}.  


\begin{figure}[t]
\begin{center}
\scalebox{0.5}{\input{troncs}}
\caption{Truncation in frequency.}
\label{fig-troncs}
\end{center}
\end{figure}
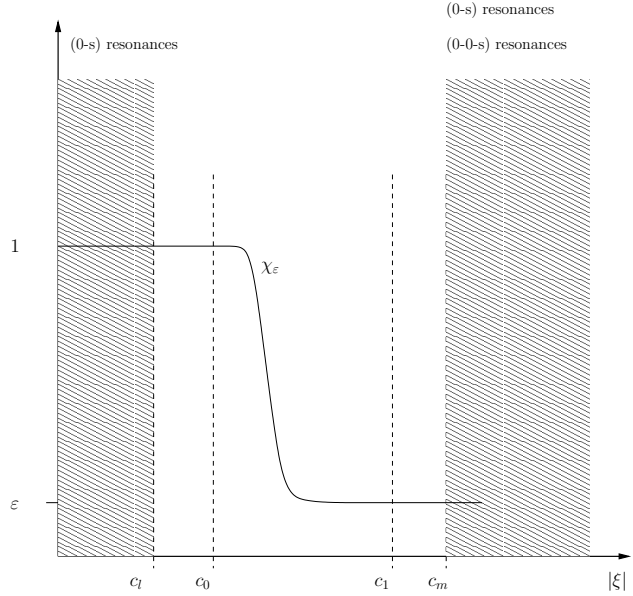

 Introduce the change of variables 
 $$  v_0 := \ompi_\e^\psi (\Pi_0) \dot u, \quad  v_s := \frac{1}{\e} \ompi_\e^\psi(\chi_\e) \ompie(\Pi_s) \dot u, \quad v := (v_0, v_s),
 $$
 where $\Pi_0, \Pi_s$ are evaluated at $(\e, \e(u_a + \e^l \dot u)).$ Then
 \begin{equation} \label{inv-diag}
  \dot u = v_0 + \e \ompie(\chi_\e^{-1}) v_s.
 \end{equation}
 With Proposition \ref{composition} and the orthogonality of $\Pi_0$ and $\Pi_s,$
 \begin{eqnarray}
  \ompie(\Pi_0) v_0 & =& v_0 + \e^2 \ompie(\Pi_0 \sharp \Pi_0) \dot u + \e^3 R_{(0)} \dot u, \nonumber \\
  \ompie(\Pi_s) v_s & = & v_s + \e \ompie(\chi_\e(\Pi_s \sharp \Pi_s) + (\chi_\e \sharp \Pi_s) \Pi_s) \dot u + \e^2 R_{(0)} \dot u. \nonumber\end{eqnarray}
We multiply (\ref{perturb}) by $\ompi_\e^\psi (\Pi_0)$ (resp. $\ompie(\chi_\e) \ompie(\Pi_s)$) to the left to find the equation satisfied by $v_0$ (resp. $v_s$). We use Proposition \ref{composition} to spell out the compositions.

 The terms in $\d_t$ are
 \begin{eqnarray} \ompi_\e^\psi(\Pi_0) \d_t \dot u & = & \d_t v_0 - \ompie(\d_t \Pi_0) \dot u \nonumber \\ 
 & = & \d_t v_0 - \frac{1}{\e} \ompie((\e \d_t \Pi_0) \Pi_0) v_0 \nonumber \\ & - & \frac{1}{\e} \ompie( \e \chi_\e^{-1} (\e \d_t \Pi_0) \Pi_s) v_s +  R_{(0)} v; \nonumber\end{eqnarray}
and
\begin{eqnarray}
  \frac{1}{\e} \ompie(\chi_\e) \ompi_\e^\psi(\Pi_s) \d_t u & = & \d_t v_s - \frac{1}{\e} \ompie(\chi_\e) \ompie( \d_t \Pi_s) \dot u \nonumber \\
 & = &  \d_t v_s - \frac{1}{\e} \ompie( (\e \d_t \Pi_s) \Pi_s) v_s \nonumber \\
 & - & \frac{1}{\e^2} \ompie(\chi_\e (\e \d_t \Pi_s) \Pi_0) v_0 \nonumber \\
 & - & \frac{1}{\e} \ompie( (\chi_\e \sharp \e \d_t \Pi_s) \Pi_0) v_0 +  R_{(0)} v. \nonumber 
  \end{eqnarray}

 
 The terms of order one are, with Assumption \ref{ass_1},
 \begin{eqnarray} \ompi_\e^\psi(\Pi_0) \ompi_\e^\psi({\cal A}) u & = & (\ompi_\e^\psi(\Pi_0 {\cal A}) + \e^2 \ompi_\e^\psi(\Pi_0 \sharp {\cal A}) + \e^3 R_{(0)}) u \nonumber \\
  & = & \ompie(\Pi_0 {\cal A}) v_0 +  \e \ompie( (\e \chi_\e^{-1})^2 \Pi_0 {\cal A} (\chi_\e \sharp \Pi_s) \Pi_s) v_s \nonumber \\ & + & \e^2 R_{(0)} v; \nonumber \end{eqnarray}
and
\begin{eqnarray}
 \frac{1}{\e} \ompie(\chi_\e) \ompi_\e^\psi(\Pi_s) \ompi_\e^\psi({\cal A}) u & = & \ompie(\Pi_s {\cal A}) v_s + \e \ompie( \e \chi_\e^{-1} (\chi_\e \sharp (\Pi_s {\cal A}))) v_s \nonumber \\ & + & \e \ompie(\rho_{s0}) v_0 + \e \ompie(\chi_\e \sharp (\Pi_s {\cal A}) \Pi_0) v_0 + \nonumber \\ & + & \e^2 R_{(0)} v; \nonumber 
 \end{eqnarray} 
where
 $$ \rho_{s0} := \chi_\e ( (\Pi_s {\cal A}) \sharp \Pi_0 + (\Pi_s \sharp {\cal A}^{(1)}) \Pi_0 + \Pi_s {\cal A} (\Pi_0 \sharp \Pi_0).$$
In the above symbolic computations, the projectors $\Pi_0$ and $\Pi_s,$ as well as ${\cal A},$ are evaluated at $(\e, \e u_a + \e^{k_0} \dot u),$ and ${\cal A}^{(1)}$ is evaluated at $(\e, u_a + \e^{k_0 - 1} \dot u).$ 


\bigskip

 
The singular terms in the right-hand sides are
 \begin{eqnarray} \ompie(\Pi_0) \ompi_\e({\cal B}) \dot u & =& \ompie(\Pi_0 {\cal B} \Pi_0) v_0 + \ompie( \e \chi_\e^{-1} \Pi_0 {\cal B} \Pi_s) v_s \nonumber \\ & + & \e R_{(0)} v \nonumber;\\
  \frac{1}{\e} \ompie(\chi_\e) \ompie(\Pi_s) \ompie({\cal B}) \dot u & = & \frac{1}{\e} \ompie(\chi_\e \Pi_s {\cal B} \Pi_0) v_0 \nonumber \\ & + & \ompie(\chi_\e (\Pi_s \sharp {\cal B}) \Pi_0 + \chi_\e \sharp (\Pi_s {\cal B}) \Pi_0) v_0 \nonumber \\ & + & \ompie(\Pi_s {\cal B} \Pi_s) v_s + \ompie( \e \chi_\e^{-1} (\chi_\e \sharp \Pi_s {\cal B}) \Pi_s)  v_s \nonumber \\ & + & \e R_{(0)} v.\nonumber  
 \end{eqnarray}
 In the above symbolic computations involving ${\cal B},$ the projectors $\Pi_0, \Pi_s$ are evaluated at $(\e, \e u_a),$ and ${\cal B}$ is evaluated at $u_a.$ 

 The other singular source term is 
 $$ \frac{1}{\e} \ompie(\chi_\e) \ompie(\Pi_s) \ompie({\cal D}) u = \frac{1}{\e} \ompie(\chi_\e \Pi_s {\cal D} \Pi_0) v_0 +  R_{(0)} v,$$ 
 where the projectors are evaluated at $(0,0),$ and ${\cal D}$ is evaluated at $u_a.$

\bigskip
 The equation in $v$ thus takes the form 
 $$ \d_t v + \frac{1}{\e^2} \ompi_\e^\psi(i A) v  = \frac{1}{\e^2} \ompie(\un B) v + \frac{1}{\e} \ompie(\un D) v + R_{(0)} v + r_a^\e,$$
 where 
 \begin{itemize}
 \item $A = A(\e, \e(u_a + \e^{k_0 - 1} \dot u)) \in C^\infty {\cal M}^1$ is defined as 
 $$i A := \left(\begin{array}{cc} {\cal A} \Pi_0 & 0 \\ 0 & {\cal A} \Pi_s \end{array}\right);$$
 \item $\un B = \un B(\e, u_a) \in C^\infty {\cal M}^0$ is defined as
 $$\un B := \left(\begin{array}{cc} 0 & 0 \\ \un B_{s0} & 0 \end{array}\right),$$
 where
 $$\begin{aligned} \un B_{s0} := & \chi_\e (\Pi_s ({\cal B} + \e {\cal D}) \Pi_0 + \e (\d_t \Pi_s) \Pi_0 + \e (\rho_{s0} - (\Pi_s \sharp {\cal B}) \Pi_0) \\ & - \e \chi_\e \sharp ( \e \d_t \Pi_s + \Pi_s {\cal A} - \Pi_s {\cal B} \Pi_0). \end{aligned}
 $$
 where $\Pi_0, \Pi_s, {\cal A}$ are evaluated at $(\e, \e u_a);$ ${\cal A}^{(1)}$ is evaluated at $(\e,u_a),$ and ${\cal B}$  and ${\cal D}$ are evaluated at $u_a.$ 
 \item $\un D = \un D(\e,u_a) = \un B' + \un B'' + E + \un F \in C^\infty {\cal M}^0$ is defined as 
  $$ \un B' := \left(\begin{array}{cc} 0 & 0 \\ 0 & \un B_s \end{array}\right), \quad \un B'' := \left(\begin{array}{cc} 0 & \un B_{0s} \\ 0 & 0 \end{array}\right),$$ where
 \begin{eqnarray}
  \un B_s & := &   - \e \chi_\e^{-1} (\chi_\e \sharp (\Pi_s ({\cal A}- {\cal B})) \Pi_s,\nonumber \\
  \un B_{0s} & := &  \e \chi_\e^{-1} ( \Pi_0 {\cal B} \Pi_s + \e \d_t \Pi_0 - (\e \chi_\e^{-1}) \Pi_0 {\cal A} (\chi_\e \sharp \Pi_s)) \Pi_s,\nonumber \end{eqnarray}
  $$E := \left(\begin{array}{cc} \Pi_0 {\cal B} \Pi_0 \\ 0 & \Pi_s {\cal B} \Pi_s \end{array}\right),$$
 where $\Pi_0, \Pi_s, {\cal A}$ are evaluated at $(0,0),$ ${\cal B}, {\cal D}$ are evaluated at $u_a,$ and
 $$\un F := \left(\begin{array}{cc} (\e \d_t \Pi_0) \Pi_0 \\ 0 & (\e \d_t \Pi_s) \Pi_s \end{array}\right),$$
 where $\e \d_t \Pi_j$ is short for $\d_v \Pi_j (0,0) \cdot (\e^2 \d_t u_a),$ $j = 0, s,$ and where $\Pi_0, \Pi_s$ are evaluated at $(0,0);$
 \item $r_a^\e := \e^{l_0 - k_0 - 2} (\ompie(\Pi_0) R_a^\e, \frac{1}{\e} \ompie(\Pi_s) R_a^\e).$
 \end{itemize}

 Next we polarize the source terms, by letting
 $$ B_{s0} := \Pi_s \un B_{s0} \Pi_0, \quad B_s := \Pi_s \un B_s \Pi_s, \quad B_{0s} := \Pi_0 \un B_{0s} \Pi_s,$$
 and
 $$ B := \left(\begin{array}{cc} 0 & 0 \\ B_{s0} & 0 \end{array}\right), \quad B' := \left(\begin{array}{cc} 0 & 0 \\ 0 & B_s \end{array}\right), \quad B'' := \left(\begin{array}{cc} 0 & B_{0s} \\ 0 & 0 \end{array}\right).$$
 We let also
 $$ F := \un F + \left(\begin{array}{cc} 0 & \un B_{0s} - B_{0s} \\ \un B_{s0} - B_{s0} & \un B_s - B_s \end{array}\right),$$
 and
 $$ D := B' + B'' + E + F.$$
 The equation is now 
\begin{equation} \label{prepared} 
  \d_t v + \frac{1}{\e^2} \ompi_\e^\psi(i A) v  = \frac{1}{\e^2} \ompie(B) v + \frac{1}{\e} \ompie(D) v + R_{(0)} v + r_a^\e.
 \end{equation}
 In \eqref{prepared}, the variables $v_0$ and $v_s$ are coupled only by order-zero terms, the leading singular term is polarized and has a nilpotent structure. The system is prepared. All the symbols in \eqref{prepared} are smooth.
 
 In the next sections, the terms $B'',$ $B'$ and $B$ will be eliminated by normal form reductions. In the subsequent $H^s_\e$ energy estimate, the term $E$ will be symmetrized, while the non-polarized term $F/\e$ will be seen to contribute to $O(1).$ 

\end{subsubsection}


\begin{subsubsection}{First reduction} \label{red0} 

 In this section, the non-resonant term $B''$ is eliminated from \eqref{prepared}.

\begin{prop}\label{p-red1} Under Assumptions {\rm \ref{ass_3}} and {\rm \ref{ass_6}}, there exists a smooth symbol $L \in C^\infty {\cal M}^{-1}$ such that
  \begin{equation} \label{first-H}
 [ \e^2 \d_t + \ompie(i A), \ompie(L(u_a)] = \ompie(B'') + \e R_{(0)}.
 \end{equation}
\end{prop}

\begin{proof} The leading term in the symbol of the source term in \eqref{first-H} is linear in $u_a:$
 $$ B_{0s} = \d_u B_{0s}(0,0) \cdot u_a + O(\e),$$ 
  and satisfies $\Pi_0 B_{0s} \Pi_s = B_{0s}.$ We look a solution $L \in C^\infty {\cal M}^{-1}$ to \eqref{first-H} in the same form: $ L = \left(\begin{array}{cc} 0 & L_{0s} \\ 0 & 0 \end{array}\right),$ where $\Pi_0 L_{0s} \Pi_s = L_{0s},$ and
  $$L_{0s} = \sum_{p\in {\cal R}^*} e^{ i p \om t/\e^2} L_{0s,p}(u_{a,p}),$$ 
 where $L_{0s,p}$ is linear in $u_{a,p}.$ Then
 $$ [ \ompie(i A), \ompie(L)] =  [ \ompie(i A(0,0), \ompie(L)],$$
  up to the commutator
  $$[ \ompie(A - A(0,0)), \ompie(L)],$$
 which, because $A$ depends on $u$ through $\e u,$ and because $L$ is assumed to belong to $C^\infty {\cal M}^{-1},$ has the form $\e R_{(0)}.$ If we suppose in addition that $L$ is smooth at the origin, we can use Proposition \ref{composition} to obtain
 $$  [ \ompie(i A(0,0), \ompie(L)] =  \ompie [i A(0,0), L] + \e R_{(0)}.$$
 Thus to solve \eqref{first-H}, it suffices to solve the equation
 \begin{equation} \label{1h} 
  \e^2 \d_t L_{0s} + [ i A(0,0), L]_{0s} = \d_u B_{0s}(0,0) \cdot u_a,\end{equation}
 up to $O(\e).$ We compute
 $$ \e^2 \d_t L_{0s} + [ i A(0,0), L]_{0s} = \sum_{j \leq n_0 < k} \sum_{p\in {\cal R}^*} e^{i p \om t/ \e^2} \Phi_{j,k,p}(0) \Pi_j L_{0s,p} \Pi_k,$$
 where the projectors $\Pi_j, \Pi_k$ are evaluated at $(0,0),$ and where $\Phi_{j,k,p}(0)$ stands the evaluation at $\e = 0$ of the phase defined in \eqref{res-eq}. Let $\chi_L$ be a smooth function on $\R^d,$ identically equal to 0 for $|\xi| \leq c_l,$ and identically equal to 1 for $|\xi| \geq c_0.$ Let then  
 $$  L_{0s,p} :=  \chi_L \sum_{j \leq n_0 < k} \Phi_{j,k,p}^{-1}(0) \Pi_j (\d_u B_{0s}(0,0) \cdot u_{a,p}) \Pi_k.$$ Because $\Phi_{j,k,p}^{-1},$ for $j \leq n_0 < k,$ is uniformly bounded for $|\xi| \geq c_l,$  the above defines a symbol $L \in C^\infty {\cal M}^{-1}.$ Besides, this symbol is smooth, and solves \eqref{1h}, up to the error term $(1 - \chi_L) B_{0s}.$ Because $|\e \chi_\e^{-1} \chi_L| \leq \e,$ this error is $O(\e).$


\end{proof}

Proposition \ref{action} allows to evaluate the action of $\ompie(L).$ Its norm, as an operator from $H^{s}_\e$ to $H^{s+1}_\e,$ is bounded by $C(|u_a|_{0,\infty}).$ Consider now the change of variables
 \begin{equation} \label{tilde} \tilde v := (\id + \e \ompie(L))^{-1} v,
 \end{equation} 
 Then
 $$ \d_t v = (\id + \e \ompie(L)) \d_t \tilde v +  \ompie(\e \d_t L) \tilde v.$$
 Because $L$ is order $-1$ and is smooth at the origin,
 \begin{eqnarray} (\id + \e \ompie(L))^{-1} \ompie(A)  (\id + \e \ompie(L)) & = &  \ompie(A) + \e [ \ompie(A), \ompie(L) ] \nonumber \\
 & = &  \ompie(A) + \e \ompie [ A, L ],\nonumber \end{eqnarray}
 up to error terms of the form $\e^2 R_{(0)}.$
 Similarly,
 $$ (\id + \e \ompie(L))^{-1} \ompie(B)  (\id + \e \ompie(L)) = \ompie(B) + \e \ompie[ B, L] +  \e^2 R_{(0)},$$
 and
 $$  (\id + \e \ompie(L))^{-1} \ompie(D)  (\id + \e \ompie(L)) = \ompie(D) +  \e R_{(0)}.$$
 The leading term in $\e$ in the commutator $$[B, L] = \left(\begin{array}{cc} - L B_{s0} & 0 \\ 0 & B_{s0} L\end{array}\right),$$
 is $O(\chi_\e) O(\e \chi_\e^{-1}) = O(\e).$ Thus the equation satisfied by $\tilde v$ is 
 \begin{equation} \label{tilde2}
  \d_t \tilde v + \frac{1}{\e^2} \ompi_\e^\psi(i A) \tilde v = \frac{1}{\e^2} \ompie(B) \tilde v + \frac{1}{\e} \ompie(\tilde D) \tilde v + R_{(0)} \tilde v + \tilde r_a^\e, \end{equation}
 where
 $\tilde r_a^\e := (\id + \e \ompie(L))^{-1} r_a^\e,$ 
  and
 $$ \tilde D:= D - [\e^2 \d_t + i A - B, L],$$
 With the above proposition,
 $$ \tilde D = B' + E + F.$$

\end{subsubsection}

\begin{subsubsection}{Second reduction}\label{red1}

 In this section, the non-resonant term $B'$ is eliminated from \eqref{tilde2}.

 \begin{prop}\label{p-red2} Under Assumptions {\rm \ref{ass_3}} and {\rm \ref{ass_6}}, there exists a smooth symbol $M \in C^\infty {\cal M}^{-1}$ such that
  \begin{equation} \label{second-H}
 [ \e^2 \d_t + \ompie(i A), \ompie(M(u_a))] = \ompie(B') + \e R_{(0)}.
 \end{equation}
\end{prop}

\begin{proof} The leading term in the symbol of the source term is linear in $u_a:$ 
 $$ B_s = \d_u B_s(0,0) \cdot u_a + O(\e),$$
 and satisfies $\Pi_s B_s \Pi_s = B_s.$ We look for a solution of \eqref{second-H} in the form of a smooth symbol $M \in C^\infty {\cal M}^{-1},$ such that $M = \left(\begin{array}{cc} 0 & 0 \\ 0 & M_s \end{array}\right),$ where $\Pi_s M_{s} \Pi_s = M_{s},$ and
  $$M_{s} = \sum_{p\in {\cal R}^*} e^{ (i p \om t)/\e^2} M_{s,p}(u_{a,p}),$$ 
  where $M_{s,p}$ is linear in $u_{a,p}.$ We check as in the proof of Proposition \ref{p-red1} that in order to solve \eqref{second-H}, it suffices to solve up to $O(\e)$ the equation
 $$ 
  \e^2 \d_t M_s + [ i A(0,0), M]_{s} = \d_u B_s(0,0) \cdot u_a. 
  $$  Let 
 $$ M_{s,p} :=  \sum_{n_0 < j,k} \Phi_{j,k,p}^{-1}(0) \Pi_j (\d_u B_s(0,0) \cdot u_{a,p}) \Pi_k,$$
 where the projectors are evaluated at $(0,0).$ Condition (iii) in Assumption \ref{ass_1} implies that for $j, k > n_0,$ $\Phi_{j,k,p}^{-1}(0) = i p \om.$ The support of $B_s$ is included $[c_0,c_1].$ Thus $M_{s,p}$ is a smooth symbol with compact support, $M$ solves \eqref{second-H}. 
\end{proof}

Proposition \ref{action} allows to evaluate the action of $\ompie(M).$ Its norm, as an operator from $H^{s}_\e$ to $H^{s+1}_\e,$ is bounded by $C(|u_a|_{0,\infty}).$ Consider now the change of variables
 \begin{equation} \label{check} \check v := (\id + \e \ompie(M))^{-1} \tilde v,
 \end{equation} 
 With the above proposition, the equation satisfied by $\check v$ is 
 \begin{equation} \label{check2}
  \d_t \check v + \frac{1}{\e^2} \ompi_\e^\psi(i A) \check v = \frac{1}{\e^2} \ompie(\check B) \check v + \frac{1}{\e} \ompie(E + F) \check v + R_{(0)} \check v + \check r_a^\e, \end{equation}
 where $\check r_a^\e := (\id + \e \ompie(M))^{-1} \tilde r_a^\e,$ and 
 $$ \check B  :=  B + \e [B, M].$$

\end{subsubsection}

\begin{subsubsection}{Third reduction} \label{red3}

In this section, the resonant term $\check B$ is eliminated from \eqref{check2}, under the transparency condition \eqref{interaction0}. 

\begin{prop} \label{prop-red-3} Under Assumptions {\rm \ref{ass_3}} and {\rm \ref{ass_6}}, there exists $N \in C^\infty {\cal M}^{-1},$ such that
 \begin{equation} \label{third-H} [ \e^2 \d_t + \ompie(i A) - \e \ompie(E), \ompie(N(u_a))] =  \ompie(\check B) + \e^2 R_{(0)},
 \end{equation}
\end{prop}

\begin{proof} 
The source term is $\check B = \left(\begin{array}{cc} 0 & 0 \\ \check B_{s0} & 0 \end{array}\right),$
 where
 $$ \check B_{s0} = \chi_\e (B^{\tt r} + \e \Pi_s {\cal D} \Pi_0) +  \e (B^{\tt nr} - \chi_\e M_s B^{\tt r}) + O(\e^2).$$
 The symbol $B^{\tt r}$ is introduced in \eqref{tt-b}, and
 $$ B^{\tt nr}(\e, u_a) := - \Pi_s (\chi_\e \sharp (\e \d_t \Pi_s + \Pi_s {\cal A}^{(0)} - \Pi_s {\cal B} \Pi_0)) \Pi_0.$$
 With the transparency assumption \eqref{interaction-d}, the term $\e \chi_\e \Pi_s {\cal D} \Pi_0$ is $O(\e^2).$
  The Taylor expansions of $B^{\tt r}$ and $B^{\tt nr}$ in their second variables are
  $$ B^{\tt r}(\e, u_a) = \d_u B^{\tt r}(\e,0) \cdot u_a + \e \d_u^2 B^{\tt r}(\e,0) \cdot (u_a, u_a) + O(\e^2),$$
  and
  $$ B^{\tt nr}(\e, u_a) = \d_u B^{\tt nr}(\e,0) \cdot u_a + O(\e).$$
 Thus, up to $O(\e^2),$ the source $\check B_{s0}$ is the sum of a linear term in $u_a,$
  \begin{equation} \label{26.1-1} \chi_\e \d_u B^{\tt r}(\e,0) \cdot u_a + \e \d_u B^{\tt nr}(\e,0) \cdot u_a,\end{equation}
   and of a bilinear term in $u_a,$
 \begin{equation} \label{26.1-2} \e \chi_\e (\d_u^2 B^{\tt r}(0,0) \cdot (u_a, u_a) - M_s(u_a) \d_u B^{\tt r}(0,0) \cdot u_a).\end{equation}
 All the terms in \eqref{26.1-1}-\eqref{26.1-2} have the product structure \eqref{structure-produit}. Accordingly, we look for $N$ in the form,
 \begin{equation} \label{forme-n} N = N^{(0)} + \e (N^{(1)} +  N^{(2)} + N^{(3)}),\end{equation}
  with 
 $$N^{(j)} = \left(\begin{array}{cc} 0 & 0 \\ N^{(j)}_{s0} & 0 \end{array}\right) \in C^\infty {\cal M}^{-1}, \quad \mbox{ for all $j,$}$$ and where
 \begin{itemize}
  \item[(a)] $N^{(0)},$ $N^{(1)}$ and $N^{(2)}$ have the structure \eqref{structure-produit};
  \item[(b)] all the entries of $N^{(3)}$ have the form $a(\xi) {\tt r}(p, q) b(\xi),$ for some Fourier multiplier $a$ and $b,$ and some symbols $p, q$ in the form \eqref{structure-produit} (the notation ${\tt r}$ is introduced in \eqref{8fev2});
  \item[(c)]  $N^{(0)}$ and $N^{(1)}$ are linear in $u_a,$ and $N^{(2)}$ and $N^{(3)}$ are bilinear in $u_a:$
 \begin{eqnarray} N^{(j)} & = & \sum_{p \in {\cal R}^*} e^{i p \om t /\e^2} N^{(j)}_p(u_{a,p}), \quad j= 0, 1,\nonumber \\
   N^{(j)} & = & \sum_{p,p' \in {\cal R}^*} e^{i (p + p') \om t /\e^2} N^{(2)}_{p,p'}(u_{a,p}, u_{a,p'}), \quad j = 2,3.\nonumber
 \end{eqnarray}
\end{itemize}
We now describe the symbols of the commutators 
  $$ [\ompie(i A - E), \ompi_\e(N^{(j)})], \quad j = 0, 1, 2, 3,$$
  using the results of section \ref{limited}. 
   
  Because the symbol $N^{(0)}$ is assumed to belong to $C^\infty {\cal M}^{-1}$ and to depend on $x$ only through $u_a,$ 
 $$ [\ompie(i A), \ompi_\e(N^{(0)})] = [\ompi_\e(i A(\e, \e u_a)), \ompi_\e(N^{(0)})] +  \e^2 R_{(0)},$$
  and the commutator in the right-hand side of the above equation is
  $$ [\ompi_\e(i A(\e,0)), \ompi_\e(N^{(0)})]  + \e [\ompi_\e(i \d_v A(0,0) \cdot u_a), \ompi_\e(N^{(0)})],$$
 up to $\e^2 R_{(0)}.$ Because $N^{(0)}$ is assumed to have the form \eqref{structure-produit} and $A(\e,0)$ is smooth, Proposition \ref{comp-trick-pdo-BIS} (ii) implies that  
 $$ [\ompi_\e(i A(\e,0), \ompi_\e(N^{(0)})] = \ompi_\e[ i A(\e,0), N^{(0)}] + \e \ompi_\e(i A(0,0) \sharp N^{(0)}),$$
 up to a remainder in $\e^2.$ Because $N^{(0)}$ depends only on $u_a,$ Proposition \ref{comp-trick-pdo-BIS} (ii) implies that this remainder has the form $\e^2 R_{(0)}.$ Remark that 
 $$E (\e, u_a) = \d_u E(0,0) \cdot u_a + O(\e).$$
  The symbol $i \d_u A - E$ is smooth, and thus Proposition \ref{comp-trick-pdo-BIS} (i) implies that the commutator
 $$  [\ompi_\e( \d_u i A(0,0) \cdot u_a - E), \ompi_\e(N^{(0)})],$$
  is equal to 
  $$ \ompi_\e[  \d_v (i A - E)(0,0) \cdot u_a, N^{(0)}] + \ompi_{\e}{\tt r}(i \d_v (A - E)(0,0) \cdot u_a, N^{(0)}),$$
 up to a remainder in $O(\e).$ With Proposition \ref{comp-trick-pdo-BIS} (i), and because $N^{(0)}$ depends only on $u_a,$ this remainder has the form $\e^2 R_{(0)}.$ 
  
  Similarly, for $j= 1,2,3$
  $$ [\ompi_\e(i A), \ompi_\e(N^{(j)})] = [\ompi_\e(i A(0,0)), \ompi_\e(N^{(j)})] + \e R_{(0)}.$$
 For $j = 1,2,$ Proposition \ref{comp-trick-pdo-BIS} (ii) implies that 
 $$ [\ompi_\e(i A(0,0)), \ompi_\e(N^{(j)})]  = \ompi_\e[ i A(0,0), N^{(j)}] + \e R_{(0)}.$$
 Lemma \ref{le-2} implies that 
 $$ [\ompi_\e(i A(0,0)), \ompi_\e(N^{(3)})]  = \ompi_\e[ i A(0,0), N^{(3)}] + \e R_{(0)}.$$

 These symbolic computations show that in order to solve \eqref{third-H}, it is sufficient to solve the system,  
  \begin{eqnarray}
  \e^2 \d_t N^{(0)}_{s0} + [i A(\e,0), N^{(0)}]_{s0} & = &\chi_{\e} \d_u B^{\tt r}(\e,0) \! \cdot \! u_a, \label{n0} \\
  \e^2 \d_t N^{(1)}_{s0} + [i A(0,0), N^{(1)}]_{s0} & = & \! \d_u B^{\tt nr}(0,0) \! \cdot \! u_a - i A_s(0,0) \sharp N^{(0)}_{s0},  \label{n1}\\
  \e^2 \d_t N^{(2)}_{s0} + [i A(0,0), N^{(2)}]_{s0} & = &  B^{\tt bl}(u_a, u_a),  \label{n2} \\
   \e^2 \d_t N^{(3)}_{s0} + [i A(0,0), N^{(3)}]_{s0} & = & {\tt r}( \d_v (i A - E)(0,0) \cdot u_a, N^{(0)})_{s0}, \label{n3} 
 \end{eqnarray}  
 with the notation,
 $$ \begin{aligned} B^{\tt bl}(u_a, u_a) & := \chi_\e (\d_u^2 B^{\tt r}(0,0)\cdot(u_a,u_a) - M_s(u_a) \d_u B^{\tt r}(0,0) \cdot u_a) \\ & + [\d_v (i A - E)(0,0) \cdot u_a, N^{(0)} ]_{s0}. \end{aligned} $$


 
 We now solve \eqref{n0}. For $j \leq n_0 < k,$ for $0 < \e <\e_0,$ let
  $$ \overline \Phi_{j,k,p}(\e) := \left\{ \begin{array}{cc} \Phi_{j,k,p}(\e), & \mbox{if } |\Phi_{j,k,p}(\e)| \geq \e^2/2, \\ \e^2/2, & \mbox{otherwise}. \end{array} \right. $$ 
  These new phases are not continuous, but they are bounded in $\xi,$ uniformly in $\e,$ for $|\xi| \leq c_l.$  
  Let then
 \begin{equation} \label{def-n-0}
  N^{(0)}_{s0,p} := \chi_\e \sum_{j \leq n_0 < k} \overline \Phi_{j,k,p}(\e)^{-1} \Pi_k(\e,0) (\d_u B^{\tt r}_0(\e, 0) \cdot u_{a,p})  \Pi_j(\e,0).
  \end{equation}
 Then every entry of $N_{s0}^{(0)}$ is a product $p(u_a) q(\e,\xi),$ where $p$ is smooth, and the transparency condition \eqref{interaction0} in Assumption \ref{ass_6} ensures that $q(\e,\xi)$ is bounded in $\xi,$ uniformly in $\e.$ Besides, by definition, $q$ is compactly supported. Thus $N^{(0)} \in C^\infty {\cal M}^{-1},$ and it solves \eqref{n0}. 
    
 In \eqref{n1} the source term $A_s(0,0) \sharp N^{(0)}_{s0}$ is identically zero (Assumption \ref{ass_3} (iii)). Let 
  $$
  N^{(1)}_{s0,p} := \sum_{j \leq n_0 < k} \Phi_{j,k,p}(0)^{-1} \Pi_k(0,0) (\d_u B^{\tt nr}(0, 0) \cdot u_{a,p}) \Pi_j(0,0),
  $$
   Because $B^{\tt nr}$ is supported in $|\xi| \in [c_0, c_1],$ far from the Klein-Gordon/acoustic resonances, the above defines a symbol in $C^\infty {\cal M}^{-1}.$ Then $N_{s0}^{(1)}$ has the structure \eqref{structure-produit}, is linear in $u_a,$ and solves \eqref{n1}. 

 Let $\chi_N$ be a smooth truncation function on $\R^d,$ identically equal to 1 for $|\xi| \leq c_1,$ and identically equal to 0 for $|\xi| \geq c_m.$ Let then 
 $$ N^{(2)}_{s0,p,p'} := \chi_N \sum_{j \leq n_0} \Psi_{k,p,p'}^{-1} \Pi_s(0,0) B^{\tt bl}(u_{a,p}, u_{a,p'}) \Pi_j(0,0).$$
 The phases $\Psi_{k,p,p'}$ are uniformly bounded away from 0 for $|\xi| \leq c_m$ (Assumption \ref{ass_3} (0-0-s)). Thus the above defines a symbol in $C^\infty {\cal M}^{-1}.$ Then $N^{(2)}$ has the structure \eqref{structure-produit}, just like the source term $B^{\tt bl},$ is bilinear in $u_a,$ and solves \eqref{n2}, up to the error $(1 - \chi_N) B^{\tt bl}.$ Because $|(1 - \chi_L) \chi_\e| \leq \e,$ the error is $O(\e).$ 
 
 Finally, let 
 $$ N^{(3)}_{s0,p,p'} := \chi_N \sum_{j \leq n_0} \Psi_{k,p,p'}^{-1} \Pi_s(0,0) {\tt r}(u_{a,p}, u_{a,p'})\Pi_j(0,0),$$
 with the notation, 
 $$ {\tt r}(u_{a,p}, u_{a,p'}) := {\tt r}( \d_v (i A - E)(0,0) \cdot u_{a,p}, N_p^{(0)}(u_{a,p'}))_{s0}.$$
 All the entries of $N^{(3)}$ have the form $a(\xi) {\tt r}(p,q) b(\xi),$ where $a$ is smooth and $p, q$ have the form \eqref{structure-produit}. The symbol $N^{(3)}$ solves \eqref{n3}, up to the error $(1 - \chi_N) {\tt r},$ which, because $N^{(0)}$ is $O(\e)$ for $|\xi| \geq c_1,$ is $O(\e)$ as well. 
 
 Finally $N$ defined by \eqref{forme-n} satisfies the assumptions (a), (b), (c), on which the symbolic computations were based, and solves \eqref{third-H}.  
 
  
 \end{proof}
 
 Proposition \ref{action-pdo-trick} and Lemma \ref{le-2} imply that the norm of $\ompi_\e(N),$ as an operator from $H^s_\e$ to $H^{s+1}_\e,$ is bounded by $C(\| u_a \|_{1,s + d_0 + 2}).$ Consider now the change of variables
 \begin{equation} \label{dernier} w := (\id + \ompi_\e(N))^{-1} \check v,
 \end{equation} 
 Then
 $$ \d_t \check v = (\id + \e \ompi_\e(N)) \d_t w +  \ompi_\e(\e \d_t N) w.$$
 Because $N$ is block triangular, 
 $$ (\id + \ompi_\e(N))^{-1} \ompie(A) (\id + \ompi_\e(N)) =  \ompie(A) + [\ompie(A), \ompi_\e(N)];$$
 $$ (\id + \ompi_\e(N))^{-1} \ompie(B) (\id + \ompi_\e(N)) =  \ompie(B),$$
 and
  $$\ompi_\e(N) \ompi_\e(\e \d_t N) = 0, \quad \ompi_\e(N) \ompie(E) \ompi_\e(N) = 0.$$
 Thus, with the above proposition, the equation satisfied by $w$ is 
 \begin{equation} \label{tilde3}
  \d_t w + \frac{1}{\e^2} \ompi_\e^\psi(i A) w  = \frac{1}{\e} \ompie(E + \check F) w + R_{(0)} w + \bar r_a^\e, \end{equation}
 where $\bar r_a^\e := (\id + \ompi_\e(N))^{-1} \tilde r_a^\e,$ and 
 $$ \ompie(\check F) := (\id + \ompi_\e(N))^{-1} \ompie(F) (\id + \ompi_\e(N)).$$

\end{subsubsection}

\begin{subsubsection}{Uniform Sobolev estimates} \label{256}

 We perform energy estimates on \eqref{tilde3}, using the symmetrizer $S$ whose existence is granted by Assumption \ref{ass_sym}. 
 
 We evaluate 
 \begin{eqnarray} \d_t (\ompi_\e(S) w, w)_{\e,s} & = & (\ompi_\e(S) \d_t w, w)_{\e,s} + (\ompi_\e(S)  w, \d_t w)_{\e,s} \nonumber  \\
 & = &  2 \Re (\ompi_\e(S) \d_t w, w)_{\e,s}.\nonumber \end{eqnarray}
 We can apply Proposition \ref{adjoint} to find,
  $$ \ompie(i S A)^* = \ompie((i S A)^*) + \e^2 R_{(0)},$$
  because $A$ depends on $u$ through $\e u.$ 
  Similarly,
   $$ \ompie(i S E)^* = \ompie((i S E)^*) + \e R_{(0)},$$
  and this implies,
  \begin{eqnarray} \frac{1}{\e^2}\Re (\ompi_\e(S) \ompie(i A) w, w)_{\e,s} & = &  (\ompie( \frac{1}{\e^2} (i S A + (i S A)^*) w, w)_{\e,s} \nonumber \\ & + &  \e^2 (R_{(0)} w, w)_{\e,s}; \nonumber\\
 \frac{1}{\e} \Re (\ompi_\e(S) \ompie(E) w, w)_{\e,s} & = &  (\ompi_\e( \frac{1}{\e}(i S E + (i S E)^*) w, w)_{\e,s} \nonumber \\  & + &  \e^2 (R_{(0)} w, w)_{\e,s}.\nonumber
 \end{eqnarray}
 Now Assumption \ref{ass_sym} implies that the symbols in the right-hand sides of the above equations all belong to $C^\infty {\cal M}^0,$ and Proposition \ref{action} implies that these symbols have the form $R_{(0)}.$
 The other source term contributes to
 \begin{equation} \label{f} \frac{1}{\e} \Re (\ompi_\e(S) \ompie(\check F) w, w)_{\e,s}.\end{equation}
Introduce the notation,
  $$ \Pi := \left(\begin{array}{cc} \Pi_0 & 0 \\ 0 & \Pi_s \end{array}\right).$$
  It follows from the definitions of the above changes of variables that
 $$ w = \ompie(\Pi) w + \e R_{(0)} w,$$ 
 and that, up to a term of the form $\e R_{(0)},$ 
 $$ \ompie(\Pi) (\id - \ompi_\e(N)) =  (\id - \ompi_\e(N))  \ompie(\Pi).$$
 Besides, up to a term of the form $\e R_{(0)},$ 
 $$ \ompie(\Pi) \ompi_\e(S) =  \ompi_\e(S) \ompie(\Pi).$$  
 Thus, up to a term of the form $\e (R_{(0)} w,w)_{\e,s},$ \eqref{f} is equal to
 $$ \Re \big(\ompi_\e(S) (\id - \ompi_\e(N)) \ompie(\tilde F) (\id - \ompi_\e(N)) w, w \big)_{\e,s},$$
 where $\tilde F := \Pi F \Pi.$
  Because $\Pi_0$ and $\Pi_s$ are projectors,
 $$ \Pi_0 \d_t \Pi_0 \Pi_0 = 0 , \quad \Pi_s  \d_t \Pi_s \Pi_s = 0,$$
 and it follows from the definition of $F$ (given in section \ref{resc}), that $\tilde F = O(\e).$
Gathering the above estimates, we find that
  \begin{equation} \label{pipi} \d_t (\ompie(S) w, w)_{\e,s} = (R_{(0)} w,w)_{\e,s} + 2 \Re (\ompie(S) \bar r_a^\e, w)_{\e,s}.\end{equation}
 The changes of variables of sections \ref{red0} to \ref{red3} define a normal form $\Psi^\e(u_a),$ such that, 
  $$ w = (\id \, + \Psi^\e(u_a)) v,$$
  and, for $\e_0$ small enough and $0 <\e < \e_0,$ both $\id \, + \Psi^\e(u_a)$ and $(\id \, + \Psi^\e(u_a))^{-1}$ are uniformly bounded as operators $H^s_\e \to H^s_\e,$ with norms depending on $\| u_a \|_{1, s + d_0 + 2}.$ From \eqref{pipi}, we obtain,
 $$ \| v (t) \|^2_{\e,s} \leq C \| \phi^\e \|^2_{\e,s} + C \int_0^t  (\| v (t')\|_{\e,s} + \e^{l_0 - k_0 - 2} \| R_a^\e (t')\|_{\e,s}) \| v (t')\|_{\e,s} dt',$$
 where the constant $C$ depends on $\e_0^{k_0 - 3 - d/2} \delta,$ $\| u_a \|_{1, s + d_0 + 2},$ $s, d,$ and on $\g.$ Because \eqref{inv-diag} implies that $\dot u = R_{(0)} v,$ Gronwall's lemma finally yields, 
 \begin{equation} \label{est-finale} \| \dot u (t) \|_{\e,s} \leq C_1 \| \phi^\e \|_{\e,s} e^{C_1 t}.\end{equation}
 
  A classical continuation argument shows that, for $\e_0$ small enough, the bound \eqref{est-finale} is valid over a time interval $[0, t_0],$ independent of $\e.$ Then \eqref{first-move} yields the asymptotic estimate \eqref{est-th}. 
\end{subsubsection}


\end{subsection}

\end{section}


\begin{section}{Application to the Euler-Maxwell equations} \label{application}

 We show in this section that the Euler-Maxwell equations satisfy the assumptions of Theorem \ref{th1}. The system we consider is (EM), introduced in section \ref{EMZ}, in the specific regime \eqref{theta_i}-\eqref{regime}. 
For the unknown
 \begin{equation} \label{inconnue}
 u^\e(t,x) = (B, E, v_e, n_e, v_i, \frac{n_i}{\a}),
 \end{equation}
 the system takes the form 
 \begin{equation} \label{le-em} \d_t u^\e + \frac{1}{\e^2} {\cal A}(\e, \e u^\e, \e \d_x) u^\e  =  \frac{1}{\e} {\cal B}(u^\e, u^\e) + {\cal G}^\e(u^\e),\end{equation}
where
 $$ {\cal A}(\e, \e u, \e \d_x) =  {\cal A}_0(\e,\e \d_x) + \e {\cal A}_1(\e, u, \e \d_x),$$
 with the notations,
  \begin{eqnarray} {\cal A}_0(\e, \xi) & :=  & \left(\begin{array}{cccccc} 0 & \xi \times & 0 & 0 & 0 & 0 \\
  - \xi \times & 0 & i & 0 & - i \frac{\e}{\theta_e} & 0 \\ 
 0 & -i & 0 & \theta_e \xi & 0 & 0 \\
  0 & 0 & \theta_e \xi \cdot & 0 & 0 & 0 \\
  0 & i \frac{\e}{\theta_e} & 0 & 0 & 0 & \e \a \xi \\
  0 & 0 & 0 & 0 & \e \a \xi \cdot & 0 \end{array}\right),\nonumber \\ 
  {\cal A}_1(\e, u, \xi) & := & \mbox{diag}(0_{\C^3}, 0_{\C^3}, \theta_e (v_e \cdot \xi), \theta_e (v_e \cdot \xi), \e (v_i \cdot \xi),\e (v_i \cdot \xi)), \nonumber \end{eqnarray}
and
$$ \begin{aligned}
  {\cal B}(u,u') & :=  (0_{\C^3}, n_e v'_e, - \theta_e v'_e \times B, 0_\C, 0_{\C^3}, 0_\C),\nonumber \\ 
 {\cal G}^\e(u) &  :=  \frac{1}{\theta_e} (0,  f^\e(n_e) v_e - \frac{1}{\theta_e}(\a n_i + \e f^\e(\a n_i)) v_i, 0, 0, \frac{1}{\theta_e} v_i \times B,0),\end{aligned}$$
 with the notation, $f^\e(x) := \e^{-2}(e^{\e x} - 1 - \e x).$

\begin{subsection}{WKB approximate solution} \label{3-1} 

 Let the initial datum 
 \begin{equation} \label{init} a = (0,  E^{0}, v_e^{0}, 0, 0, 0) \in H^\s,
 \end{equation}
 for some large $\s,$ where the electric field satisfies 
  $$ \nabla \cdot E^0 = 0,$$
 in accordance with \eqref{3fev}.  We assume that $a$ is polarized (or well-prepared), in the sense that
 $$ E^0 = \tilde E + (\tilde E)^*, \quad v_e^0 = \frac{i}{\om} \tilde E - \frac{i}{\om} \tilde E^*,$$
  for some fundamental frequency $\om,$ defined in terms of $\om_{pe},$ and some complex amplitude $\tilde E$ (above, $(\tilde E)^*$ denotes the complex conjugate of $\tilde E$).

 Because the conservative form of the convective terms in the equations of conservations of mass in $\mbox{(EM)}^\sharp$ allows simple formal computations, we carry out the WKB expansion on $\mbox{(EM)}^\sharp$ rather than on (EM). 
 
 Consider $\mbox{(EM)}^\sharp$ in the regime \eqref{theta_i}. We look for an approximate solution $u_{app}^\e$ in the form of a profile
 \begin{equation} \label{ansatz-sharp}  u_{app}^\e(t,x) = \e [{\bf u}_{app}^\e (\e t, x, \theta)]_{\theta = \om t / \e},\end{equation}
 where ${\bf u}_{app}^\e$ has a WKB expansion, 
 $$u_{app}^\e = \sum_{m = 0}^M \e^m {\bf u}_m^\e,$$
  such that for all $m,$ ${\bf u}_m$ is a trigonometric polynomial in $\theta,$
 $$ {\bf u}_m(t, x, \theta) = \sum_{p \in {\cal R}_m} e^{i p \theta} u_{m,p}(t,x).$$
  The sets ${\cal R}_m \subset \Z$ are finite and the $u_{m,p}$ are assumed to satisfy
   $$ u_{m,p}  \in W^{1,\infty}([0,t^*(\e)], W^{k(m),\infty}(\R^3)),$$
 for some $t^*(\e) > 0$ and some large $k(m).$ We plug this ansatz in $\mbox{(EM)}^\sharp$ and find a cascade of WKB equations, which we now describe. We sometimes use below the notation $(v)_p$ to denote the $p^{\mbox{{\footnotesize th}}}$ harmonic in $\theta$ of a trigonometric polynomial ${\bf v}(t,x,\theta).$
 
 \medskip
\emph{Equations for the terms in $O(1/\e^2).$} 
 
  $$\begin{aligned}   \om\d_\theta {\bf E}_0 & =  {\bf v}_{e 0},\nonumber \\
 \om\d_\theta {\bf v}_{e 0} &=  - {\bf E}_0, \nonumber \\
  \om\d_\theta ({\bf B}_0, {\bf n}_{e 0}^\sharp, {\bf v}_{i 0}, {\bf n}_{i 0}^\sharp) &= 0.\nonumber
 \end{aligned}$$
 The dispersion relation is $$ \om^2 - 1 = 0.$$ We choose $\om =1.$ With this choice, and because $\e$ was set to be equal to $(\om_{pe} t_0)^{-1},$ we find that 
  $ e^{i \om t/\e^2} = e^{i \om_{pe} T},$
 where $T$ represents the physical time, $t = \e t_0 T,$ and $\om_{pe}$ is the electronic plasma frequency \eqref{o-p}. Thus the waves we consider are oscillating at the electronic plasma frequency. These waves are called \emph{plasma waves} in the physical literature. The set of characteristic harmonics is $${\cal R}^* = \{-1,1\}.$$ The first term of the expansion satisfies
 $$ E_{0,0} = 0,\; v_{e 0,0} = 0, \; {\bf B}_0 = B_{0,0}, \; {\bf n}_{e 0}^\sharp = n_{e 0,0}^\sharp, \; {\bf v}_{i 0} = v_{i 0,0},\; {\bf n}_{i 0}^\sharp = n_{i 0,0}^\sharp,$$
 and 
 \begin{equation}
  {\bf E}_0 = E_{0,1} e^{i \theta} + E_{0,-1} e^{- i \theta}, \quad  {\bf v}_{e 0} = \frac{i}{\om} E_{0,1} e^{i \theta} - \frac{i}{\om} E_{0,-1} e^{- i \theta}. \label{pola_leading}\end{equation}

\bigskip


\emph{Equations for the terms in $O(1/\e)$.}

 $$\begin{aligned}  \om \d_\theta {\bf B}_1 + \nabla \times {\bf E}_0 &=  0,\nonumber \\
 \om\d_\theta {\bf E}_1 - \nabla \times {\bf B}_0 &= {\bf v}_{e 1} + {\bf n}_{e 0}^\sharp {\bf v}_{e 0} - \frac{1}{\theta_e} {\bf v}_{i 0},\nonumber \\
 \om\d_\theta {\bf v}_{e 1} &= - \theta_e \nabla {\bf n}_{e 0}^\sharp - ({\bf E}_1 + \theta_e {\bf v}_{e 0} \times {\bf B}_0),\nonumber \\ 
\om\d_\theta {\bf n}_{e 1}^\sharp + \theta_e \nabla \cdot {\bf v}_{e 0} &= 0, \nonumber \\
\om\d_\theta {\bf v}_{i 1}  &= \frac{1}{\theta_e} {\bf E}_0,\nonumber \\
 \om\d_\theta {\bf n}_{i 1}^\sharp &= 0.\nonumber 
 \end{aligned}$$
  Let $$ L_0 := \left( \begin{array}{cc} 0 & -1 \\ 1  &  0 \end{array} \right) \in {\cal L}(\C^6).$$ 
 The frequency $\om$ was chosen so that
  $\mbox{det }(i p \om + L_0) = 0.$ In $\C^6,$
  $$ (i p \om + L_0) a = \left(\begin{array}{c} b_1 \\ b_2 \end{array}\right),$$
  implies the compatibility condition
 $$ i p \om b_1 + b_2 = 0.$$
The oscillating terms satisfy
  $$ ( i p\om  + L_0) \left( \begin{array}{c} {\bf E}_{1} \\ {\bf v}_{e 1} \end{array} \right)_p = \left( \begin{array}{c}  {\bf n}_{e 0}^\sharp {\bf v}_{e 0} \\ - \theta_e {\bf v}_{e 0} \times {\bf B}_0  \end{array} \right)_p, \quad p = -1, 1,$$
 and this implies, 
  \begin{equation} \label{compa} i p \om ( {\bf n}_{e 0}^\sharp {\bf v}_{e 0})_p - \theta_e ({\bf v}_{e 0} \times {\bf B}_0)_p  = 0.\end{equation}
 Besides, we find the polarization conditions,
 $$ i p \om B_{1, p} = - \nabla \times E_{0, p}, \;\;  i p \om v_{i 1,p} = \frac{1}{\theta_e} E_{0,p}, \;\;  n_{i 1,p}^\sharp = 0, \;\; n_{e 1 p}^\sharp = \frac{- \theta_e}{(p \om)^2} \nabla \cdot E_{0 p}.$$
 The nonoscillating terms satisfy, 
  \begin{equation} \label{QN1} -  \nabla \times B_{0,0} =  v_{e 1,0}  - \frac{1}{\theta_e} v_{i 0,0},\;\; E_{1, 0} =  \theta_e \nabla n_{e 0 ,0}^\sharp. \end{equation}
  
  \medskip


\emph{Equations for the terms in $O(1).$}

\vspace{-0.4cm}

 \begin{eqnarray} \om \d_\theta {\bf B}_2 + \d_t  {\bf B}_0 & + &  \nabla \times {\bf E}_1 = 0,\nonumber \\
 \om \d_\theta {\bf E}_2 + \d_t {\bf E}_0 & - &  \nabla \times {\bf B}_1 = {\bf v}_{e 2} + {\bf n}_{e 1}^\sharp {\bf v}_{e 0} + {\bf n}_{e 0}^\sharp {\bf v}_{e 1} - \frac{1}{\theta_e} ( {\bf v}_{i 1}  + {\bf n}_{i 0}^\sharp {\bf v}_{i 0}) ,\nonumber \\
 \om\d_\theta {\bf v}_{e 2} + \d_t {\bf v}_{e 0} & + &  \theta_e ( {\bf v}_{e 0} \cdot \nabla)  {\bf v}_{e 0} = - \theta_e \nabla {\bf n}_{e 1}^\sharp + \theta_e {\bf n}_{e 0}^\sharp \nabla {\bf n}_{e 0}^\sharp \nonumber \\ & - &  ( {\bf E}_2 + \theta_e ( {\bf v}_{e 1} \times {\bf B}_0 + {\bf v}_{e 0} \times {\bf B}_1)),\nonumber \\
  \om\d_\theta {\bf n}_{e 2}^\sharp & + &  \d_t  {\bf n}_{e 0}^\sharp + \theta_e \nabla \cdot ( {\bf v}_{e 1} + {\bf n}_{e 0}^\sharp {\bf v}_{e 0}) = 0, \nonumber \\
\om\d_\theta {\bf v}_{i 2} & + &  \d_t {\bf v}_{ i 0} = - \a^2 \nabla {\bf n}_{ i 0}^\sharp + \frac{1}{\theta_e} ( {\bf E}_1 +  {\bf v}_{i 0} \times {\bf B}_0),\nonumber \\
 \om\d_\theta {\bf n}_{i 2}^\sharp & + & \d_t  {\bf n}_{i 0}^\sharp + \nabla \cdot {\bf v}_{i 0} = 0.\nonumber 
 \end{eqnarray}
  Because \eqref{QN1} implies that $E_{1,0}$ is a gradient, the first equation implies that $\d_t B_{0, 0} = 0.$ This yields ${\bf B}_0 = 0,$ and, with \eqref{compa}, we find that 
   $({\bf n}_{e 0}^\sharp {\bf v}_{e 0})_p = 0.$
 Because ${\bf v}_{e 0, p}$ is assumed to be non identically zero (see the form of the initial condition \eqref{init}), this implies $ {n}_{e 0, 0}^\sharp = 0,$ and finally ${\bf n}_{e0}^\sharp = 0.$ That is, the electronic fluctuation of density is $O(\e),$ in accordance with the rescaling of section \ref{resc}. For $p \in \{-1, 1\},$ 
$$ (i p  + L_0) \left( \begin{array}{c} {\bf E}_{2} \\ {\bf v}_{e 2} \end{array} \right)_p \!=\! \left( \begin{array}{c} - \d_t {\bf E}_{0} + \nabla \times {\bf B}_{1} +  ({\bf n}_{e 1}^\sharp {\bf v}_{e 0} + {\bf n}_{e 0}^\sharp {\bf v}_{e 1}) - \frac{1}{\theta_e} {\bf v}_{i 1} \\  - \d_t {\bf v}_{e 0} - \theta_e (\nabla {\bf n}_{e 1}^\sharp - ({\bf v}_{e 0} \cdot \nabla)  {\bf v}_{e 0}) - ( {\bf v}_{e 0} \times {\bf B}_1)) \end{array} \right)_p\!.$$
 In the above right-hand side, the nonlinear terms are  
$$  ({\bf n}_{e 1}^\sharp {\bf v}_{e 0})_p = n_{e 1,0}^\sharp v_{e 0,p}, \quad ({\bf v}_{e 0} \times {\bf B}_1)_p = v_{e 0,p} \times B_{1,0}, \quad p = -1, 1,$$ 
  and 
  $$({\bf v}_{e 0} \cdot \nabla ) {\bf v}_{e 0})_p = 0, \quad p = -1, 1,$$ a transparency relation for the convective term. The compatibilty relation is the Schr\"{o}dinger equation for the electric field,  
 \begin{equation} \label{Schro3} - 2i p \om \d_t E_{0,p} + \Delta_e E_{0,p}  - \frac{1}{\theta_e^2} E_{0,p} -  n_{e 1,0}^\sharp  E_{0,p} + \frac{\theta_e}{i p \om} E_{0 p} \times B_{1,0} = 0,
 \end{equation}
 where
 \begin{equation} \label{delta-e} \Delta_e z := \theta_e^2 \nabla ( \nabla \cdot z ) -   \nabla \times( \nabla \times z).\end{equation}
 The nonoscillating terms satisfy
 \begin{equation} \label{couplage}
  \theta_e \big( ( {\bf v}_{e 0} \cdot \nabla ) {\bf v}_{e 0} \big)_0  = - \theta_e  \nabla n_{e 1,0}^\sharp - ( E_{2,0} + \theta_e ({\bf v}_{e 0} \times {\bf B}_1)_0).
 \end{equation}
 The above equality (\ref{couplage}) is the crucial equation that couples the Schr\"{o}dinger equation (\ref{Schro3}) and the evolution equation for $n_{e 1,0}^\sharp,$ which is made explicit below.  The equations for the terms in $O(1)$ also contain the relation:
 \begin{equation} \label{quasineutralite}
   \nabla \times  B_{1,0}   =   v_{e 2,0}  - \frac{1}{\theta_e}   v_{i 1,0}  +  ({\bf n}_{e 1}^\sharp {\bf v}_{ e 0})_0,
 \end{equation}
 and a linear wave equation for $v_{ i 0,0}$ and $n_{i 1,0}^\sharp.$ Because the initial data for $v_{ i 0,0}$ and $n_{i 1,0}^\sharp$ are null, $n_{i 0,0}^\sharp = 0$ and $v_{i 0,0} = 0.$ With (\ref{QN1}), and because $B_{0,0} = 0,$ this implies $v_{e 1,0} = 0.$ 

\bigskip


\emph{Equations for the terms in $O(\e).$}

\vspace{-.4cm}
\begin{eqnarray}
  \om \d_\theta {\bf B}_3 + \d_t {\bf B}_1 & + &  \nabla \times {\bf E}_2 = 0,\nonumber \\
 \om \d_\theta {\bf E}_3 + \d_t {\bf E}_1 & - &  \nabla \times {\bf B}_2 = {\bf v}_{e 3} + {\bf n}_{e 1}^\sharp {\bf v}_{e 1} + {\bf n}_{e 0}^\sharp {\bf v}_{e 2}  + {\bf n}_{e 2}^\sharp {\bf v}_{e 0} \nonumber \\ & - & \frac{1}{\theta_e} ( {\bf v}_{i 1}  + {\bf n}_{i 1}^\sharp {\bf v}_{i 0} + {\bf n}_{i 0}^\sharp {\bf v}_{i 1}) ,\nonumber \\
 \om\d_\theta {\bf v}_{e 3} + \d_t {\bf v}_{e 1} & + &  \theta_e(  ({\bf v}_{e 1} \cdot \nabla)  {\bf v}_{e 0}  +  ({\bf v}_{e 0} \cdot \nabla)  {\bf v}_{e 1} ) = - \theta_e \nabla {\bf n}_{e 2}^\sharp \nonumber \\ & + &  \theta_e ( {\bf n}_{e 1}^\sharp\nabla {\bf n}_{e 0}^\sharp + {\bf n}_{e 1}^\sharp \nabla {\bf n}_{e 1}^\sharp ) \nonumber \\ & - &  ( {\bf E}_3 + \theta_e ({\bf v}_{e 2} \times {\bf B}_0 + {\bf v}_{e 1} \times {\bf B}_1 + {\bf v}_{e 0} \times {\bf B}_2)),\nonumber \\
\om\d_\theta {\bf v}_{i 3} + \d_t {\bf v}_{ i 1} & + &  ({\bf v}_{i 0} \cdot \nabla) {\bf v}_{i 0} = - \a^2 \nabla {\bf n}_{ i 1}^\sharp + \a^2 {\bf n}_{i 0}^\sharp \nabla {\bf n}_{i 0}^\sharp \nonumber \\ & + &  \frac{1}{\theta_e} ( {\bf E}_2 + {\bf v}_{i 1} \times {\bf B}_0 + {\bf v}_{i 0} \times {\bf B}_1),\nonumber \\
 \om\d_\theta {\bf n}_{e 3}^\sharp + \d_t  {\bf n}_{e 1}^\sharp & + &  \theta_e \nabla \cdot ( {\bf v}_{e 2} +  {\bf n}_{e 1} ^\sharp{\bf v}_{e 0} +  {\bf n}_{e 0}^\sharp {\bf v}_{e 1} ) = 0, \nonumber \\
 \om\d_\theta {\bf n}_{i 3}^\sharp + \d_t  {\bf n}_{i 1}^\sharp & + &  \nabla \cdot ( {\bf v}_{i 1} +   {\bf n}_{ i 0}^\sharp {\bf v}_{ i 0}) = 0.\nonumber 
 \end{eqnarray}
 In particular, 
 $$ \begin{aligned} \d_t v_{i 1,0} + \a^2 \nabla n_{i 1,0}^\sharp & =  \frac{1}{\theta_e} E_{2,0},  \\
   \d_t n_{i 1,0}^\sharp + \nabla \cdot v_{i 1,0} & =   0,\\
  \d_t  n_{e 1,0}^\sharp + \theta_e \nabla \cdot ( v_{e 2,0} + ({\bf n}_{e 1}^\sharp {\bf v}_{e 0})_0 ) & =   0.  
 \end{aligned}$$
 The last two equations in the above system, together with (\ref{quasineutralite}), imply the quasineutrality relation,
 \begin{equation} \label{qn-2006} n_{e 1,0}^\sharp = n_{i 1, 0}^\sharp.\end{equation}
 The first two equations in the above system, together with \eqref{couplage} and \eqref{qn-2006}, give
  $$ \left\{ \begin{aligned} \d_t v_{i 1,0} + (\a^2 + 1 ) \nabla n_{i 1,0}^\sharp  & = -  ({\bf v}_{e 0} \cdot  \nabla) {\bf v}_{e 0} + {\bf v}_{e 0} \times {\bf B}_1)_0,  \\
   \d_t n_{i 1,0}^\sharp  + \nabla \cdot v_{i 1,0}  & = 0,\end{aligned}\right. $$
 where the nonlinear term can be computed with the above polarization conditions,
 $$ \big(({\bf v}_{e 0} \cdot  \nabla) {\bf v}_{e 0} + {\bf v}_{e 0} \times {\bf B}_1 \big)_0 = \nabla |E_{0 ,p}|^2.$$
  The equations at order $O(1)$ also yield
 $$ \d_t B_{1, 0} +  \nabla \times E_{2, 0} = 0,$$ 
 which, together with (\ref{couplage}), implies that $E_{2,0}$ is a gradient. Hence the term $B_{1,0}$ in \eqref{Schro3} vanishes.
 
   Finally, $E_{0, p}$ and $n_{i 1,0}$ satisfy the vector Zakharov system,
 \begin{equation} \label{Z}  \left\{ \begin{aligned} - 2 i p \d_t E_{0,p} + \Delta_e E_{0,p}  - \frac{1}{\theta_e^2} E_{0,p} & = &   n_{i 1,0}^\sharp  E_{0,p} \\
   (\d_t^2 - (\a^2 + 1) \Delta) n_{i 1,0}^\sharp & = &  - \Delta | E_{0,p}|^2. \end{aligned}\right. \end{equation}
 In \eqref{Z}, $ p =1$ or $p = -1.$ The Laplace-type operator $\Delta_e$ was introduced in \eqref{delta-e}. With the initial condition $E_{0,1} = E,$ $n_{i 1,0} = 0,$ and $\d_t n_{i 1,0} = 0,$ Ozawa and Tsutsumi's result \cite{OT} guarantees the existence and uniqueness of a solution $(E_{0,1}, n_{i 1,0})$ to \eqref{Z}, over a time interval $[0, t^*),$ with the same Sobolev regularity as the initial condition. Note that the crucial coupling term $\Delta | E_{0,p}|^2$ comes from the convective term and from the Lorentz force term. In the above Schr\"{o}dinger equation, the term $\frac{1}{\theta_e^2} E_{0,p}$ means to a small shift in frequency. This term was not present in the (Z) system given in the introduction; it accounts for the contribution of the ions to the fundamental frequency. Indeed, by letting $\e = (\om_{pe} t_0)^{-1},$ we took as a reference the electronic plasma frequency $\om_{pe},$ which is only an approximation of the plasma frequency,
  $$ \om_p = \sqrt{4 \pi e^2 n_0 (\frac{1}{m_e} + \frac{1}{m_i})}.$$ 
 Because the ration $m_e/m_i$ is equal to $\e^2/\theta_e^2$ (a consequence of \eqref{theta_i}), at first order in $\e^2,$ $$\om_p = \om_{pe}(1 + \frac{1}{2} \frac{\e^2}{\theta_e^2} ),$$
 and thus the shift in frequency in \eqref{Z} means that the electric field actually oscillates at the plasma frequency $\om_p.$

 \bigskip

\emph{Higher-order terms.} The WKB expansion can be carried out up to any order. For $m \geq 2,$ the terms $E_{m,p}, n_{i, m+1,0},$ $p = -1, 1$ are seen to satisfy a linearized Zakharov system of the form 
  \begin{equation} \label{Z_m}  \left\{ \begin{aligned} (- 2 i p   \d_t + \Delta_e  - \frac{1}{\theta_e^2}) E_{m,p} =    n_{i m+1,0}^\sharp  E_{0,p} + n_{i 1, 0}^\sharp E_{m, p} + r_{m,p},\\
   (\d_t^2 - (\a^2 + 1) \Delta) n_{i m+1,0}^\sharp  =  - \Delta ( E_{m, p} E_{0, -p} + E_{0, p} E_{m, -p}) + r_{m,p}, \end{aligned}\right. \end{equation}
 where $r_{m,p}$ represent the $p^{\mbox{\footnotesize th}}$ harmonics of smooth functions of the profiles $\{ \d_x^{k} {\bf u}_{k'}, k, k' \leq m-1\}.$ The system (\ref{Z_m}) is the linearization of (\ref{Z}) around $(E_{0, p}, n_{i 1,0}^\sharp).$ The initial data for $E_{m,p}$ and $n_{i m+1,0}^\sharp$ are null. Ozawa and Tsutsumi's method for the Zakharov equations \cite{OT} allows to solve the initial value problem for \eqref{Z_m}. It provides an existence time which is a priori smaller than the existence time of the data $r_{m,p},\, E_{0,p}$ and $n_{i 1,0}^\sharp.$ The solution has the same Sobolev regularity as the data. 
 
 Thus, by induction, we can construct a family of profiles ${\bf u}_m$ that determines an approximate solution $u_{app}^\e,$ as follows. If we assume the ${\bf u}_k,$ for $k \leq m-1,$ to be known, with enough Sobolev regularity, then
 \begin{itemize}
  \item
  
   $E_{m,p}$ and $n_{i m+1,0}$ are defined as the unique solution of $\mbox{(Z)}_m,$
  \item the terms $u_{m,p},$ for $p \in {\cal R}^*,$ are deduced from $E_{m,p}$ and $n_{i m+1,0}$ by polarization conditions similar to the ones found in the first terms of the expansion;
  \item the terms $u_{m,p},$ for $p \neq {\cal R}^*,$ are computed by elliptic inversions; that is, in terms of $(i p \om + L_0)^{-1} r_{m-1,p}.$
  \end{itemize}
 We obtain a profile ${\bf u}_m,$ whose Sobolev regularity is smaller than the regularity of ${\bf u}_{m-1}$ by one, and the induction is complete. The profile 
 $${\bf u}_{app}^\e = ({\bf B}^\e, {\bf E}^\e, {\bf v}_e^\e, {\bf n}_e^{\e,\sharp}, {\bf v}_i^\e, {\bf n}_i^{\e,\sharp}),$$ provides, via \eqref{ansatz-sharp}, an approximate solution to $\mbox{(EM)}^\sharp.$ Then
  $$ u_a^\e = [ {\bf B}^\e, {\bf E}^\e, {\bf v}_e^\e, \frac{1}{\e} \log(1 + \e {\bf n}_e^{\e,\sharp}),  {\bf v}_i^\e, \frac{1}{\e} \frac{1}{\a} \log(1 + \e {\bf n}_i^{\e,\sharp}) ]_{\theta = \om t /\e^2},$$
  is an approximate solution to \eqref{le-em}, in the sense of \eqref{j_0-app}. The approximate solution can be made arbitrarily precise, in a Sobolev norm, provided the initial data have enough Sobolev regularity. Then, by construction, $u_a^\e$ satisfies all the conditions stated in Assumption \ref{ass_5}.

\end{subsection}


\begin{subsection}{Stability of the approximate solution} \label{3-2} 

 Consider now a perturbation of the initial condition $a$ introduced in \eqref{init}:
 \begin{equation} \label{1-2} a^\e := a + \e^{k_0} \varphi^\e,\end{equation}
 where $\phi^\e$ is a bounded family in $H^{\s - l_0 - 2}_\e,$ and
   $ k_0 + 3 \leq l_0 < \s - 2 - \frac{d}{2}.$
 We need to assume that \eqref{3fev} is satisfied so that $a^\e$ is a proper initial datum for the Euler-Maxwell system. This amounts to the conditions,
  $$
   \nabla \cdot B_{\phi^\e} = 0, \quad \nabla \cdot E_{\phi^\e} = \frac{1}{\e \theta_e} (n_{e,\phi^\e} - n_{i,\phi^\e}),
   $$  
  for the coordinates of the perturbation $\phi^\e.$

 Consider the approximate solution $u_a^\e$ to \eqref{le-em}, associated with the initial datum $a,$ at order $l_0.$

\begin{theo} \label{theo} If $k_0 > 3 + \frac{3}{2},$ the system {\rm \eqref{le-em}}, together with the initial datum \eqref{1-2}, has a unique solution $u^\e$ defined over $[0, t^*),$ independent of $\e.$ Moreover, for all $0 < t_0 < t^*,$ there holds, for $\e$ small enough,
 \begin{equation} \label{est-EM} \sup_{0 \leq t \leq t_0} \| {u}^\e - {u}_a^\e \|_{\e,s} \leq C \e^{k_0-1},\end{equation}
 with a constant $C$ independent of $\e,$ and a Sobolev index $s > \frac{d}{2}.$   
\end{theo}

 The proof shows that the existence time $t^*$ is bounded from below by the existence time of the approximate solution $u_a^\e.$ Note that the estimate \eqref{est-EM}, the condition $k_0 > 3 + \frac{3}{2},$ and the description of the approximate solution given in the above section, imply
 \begin{equation} \label{justif} \sup_{0 \leq t \leq t_0} \sup_{x \in \R^3} ( | E^\e - (E_{0,1} e^{i \om t / \e^2} + c. c.)| + | n^\e - \e n_{i 1,0} | ) \leq C \e^2,\end{equation}
 where $E^\e, n^\e$ are coordinates of the solution $u^\e$ of \eqref{le-em}, as in \eqref{inconnue}, where $n_{i 1,0} = \frac{1}{\e} \log(1 + \e n_{i 1,0}^\sharp),$ and $E_{0,1}$ and $n_{i 1,0}^\sharp$ solve the Zakharov system \eqref{Z}, with the initial condition $E_{0,1}(t = 0) = E^0,$ $n_{i 1,0}^\sharp(t = 0) = 0,$ $\d_t n_{i 1,0}^\sharp(t = 0) = 0.$  The asymptotic estimate \eqref{justif} is thus the estimate that validates the Zakharov model, as it actually gives a description of the electric field $E^\e$ and the fluctuation of density $n^\e$ in (EM) by means of the solution $(E_{0,1}, n_{i 1,0})$ of $\mbox{(Z)}_0.$

 Theorem \ref{theo} follows as a corollary of Theorem \ref{th1} if one can prove that the (EM) system satisfies the assumptions of Theorem \ref{th1}. An approximate solution satisfying Assumption \ref{ass_5} was constructed in the above section. The next sections are devoted to the verification of Assumptions \ref{ass_1} (hyperbolic structure, regularity of the eigenvalues and eigenprojectors), \ref{ass_3} (localization of the resonances) and \ref{ass_6} (transparency).


\begin{subsubsection}{Eigenvalues and eigenvectors} \label{3-2-2}

 We check in this section that the (EM) system satisfies Assumption \ref{ass_1}. The operator ${\cal A}$ defined at the beginning of section \ref{application} obviously belongs to $C^\infty {\cal M}^1.$ For all $\e, u ,\xi,$ ${\cal A}(\e,u,\xi)$ is a hermitian matrix, and satisfies a decomposition of the form \eqref{spectral-dec}. It remains to check hypothesis (i) to (iii) in Assumption \ref{ass_1}. 

 The eigenvalues arising in the spectral decomposition \eqref{spectral-dec} are the solutions, 
 $$ \om = \lam(\e, u, \xi), \quad \e >0, u \in \R^{14}, \xi \in \R^3,$$
 of the polynomial equation in $\om,$
 \begin{equation} \label{disp0} \mbox{det }(i \om + {\cal A}(\e, u, i \xi)) = 0.\end{equation}
 A look at the definition of ${\cal A}_1$ shows that the eigenvalues depend on $u$ only through the scalar terms, 
 $${\tt x}:= \e \theta_e v_e \cdot \xi, \quad  \mbox{and} \quad {\tt y}:= \e^2 v_i \cdot \xi,$$
 representing the electronic and ionic convections ($v_e$ and $v_i$ are coordinates of $u$).
 Equation \eqref{disp0} factorizes into a transverse, degree four equation,
 \begin{equation}\label{disp1}
   ( \om - {\tt x}) (\om - {\tt y})(\om^2 - 1 - |\xi|^2 - \frac{\e^2}{\theta_e^2}) = {\tt x} (\om- {\tt y}) + \frac{\e^2}{\theta_e^2} {\tt y}( \om - {\tt x}),
 \end{equation}
 and a longitudinal, degree five equation,
\begin{equation} \label{disp2} \begin{aligned} 
 \om( (\om - {\tt y})^2- \e^2 \alpha^2 |\xi|^2) ( (\om - {\tt x})^2 - 1 - \theta_e^2 |\xi|^2 ) & = - {\tt x}( (\om - {\tt y})^2 - \e^2 \a^2 |\xi|^2) \\ & + \frac{\e^2}{\theta_e^2} (\om - {\tt y})((\om - {\tt x})^2 - \theta_e^2 |\xi|^2).
 			       \end{aligned}\end{equation}

For all $\e, {\tt x}, {\tt y}, \xi,$ the Kernel of ${\cal A}(\e, u, \xi)$ has dimension one. It is generated by
 $$ e_0 := \big( \frac{\xi}{|\xi|}, 0, 0, 0, 0, 0 \big).$$
 The solutions of \eqref{disp1} have (algebraic and geometric) multiplicity two, while the solutions of \eqref{disp2} have multiplicity (algebraic and geometric) one. 
The solutions of \eqref{disp1} and \eqref{disp2} are algebraic functions of $\e, \xi, {\tt x}, {\tt y}.$ Evaluations of these functions $\om = \om(\e, ({\tt x}, {\tt y}), \xi)$ at $({\tt x}, {\tt y})= (\e \theta_e v_e \cdot \xi, \e^2 v_i \cdot \xi)$ give the non-zero eigenvalues of the (EM) system. 

 For $({\tt x}, {\tt y}) = (0,0),$ equation \eqref{disp1} and \eqref{disp2} simplify to
 \begin{equation} \label{disp1.1}  \om^2 (\om^2 - 1 - |\xi|^2 - \frac{\e^2}{\theta_e^2}) = 0,\end{equation}
 and
 \begin{equation} \label{disp1.2}  \om( (\om^2 - \e^2 \alpha^2 |\xi|^2) ( \om^2 - 1 - \theta_e^2 |\xi|^2 ) - \frac{\e^2}{\theta_e^2} ( \om^2 - \theta_e^2 |\xi|^2)) = 0.\end{equation}
    The solutions $\om = \om(\e,\xi)$ of \eqref{disp1.1}-\eqref{disp1.2} are represented on figure \ref{fig-char1}. 

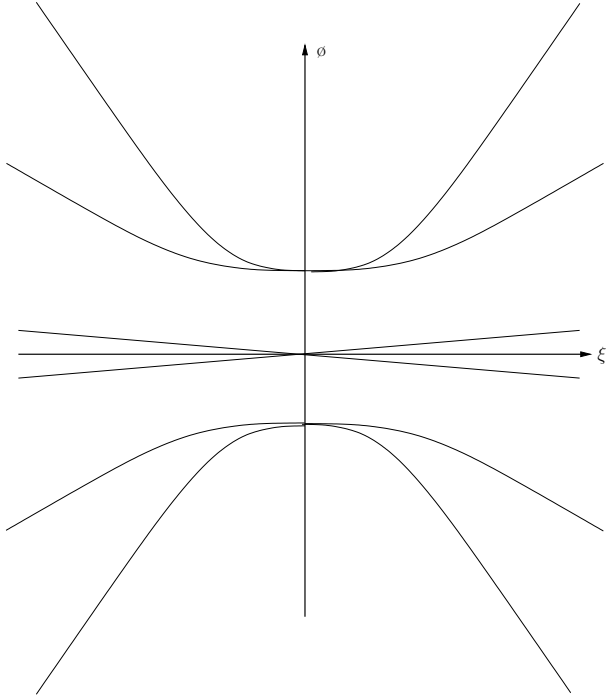
\begin{figure}[t]
\begin{center}
\scalebox{0.5}{\input{char1}}
\caption{The characteristic variety for the operator linearized around 0.}
\label{fig-char1}
\end{center}
\end{figure}

 The transverse modes at $(\e,0,\xi)$ are
 \begin{equation} \label{trans1} \lam_s(e_-) = 0, \quad \lam_s(e_+) = 0, \quad \lam_\pm = \pm \sqrt{1 + |\xi|^2 + \frac{\e^2}{\theta_e^2}}.\end{equation}
 The Klein-Gordon, longitudinal modes at $(\e,0,\xi)$ are
 $$ \mu_\pm = \pm \sqrt{1 + \theta_e^2 |\xi|^2} + \frac{\e}{2\theta_e^2 ( 1 + \theta_e^2 |\xi|^2)} + O(\e^2),$$
 and the acoustic, longitudinal modes at $(\e,0,\xi)$ are
  $$ \mu_s(e_-) = 0, \quad  \mu_s(e_+)_\pm = O(\e), \quad \mu_\pm = \pm \sqrt{1 + \theta_e^2 |\xi|^2} + O(\e),$$
  locally uniformly in $\xi.$ The Klein-Gordon modes $\mu_\pm$ have constant multiplicity, hence are analytical in $\e,\xi.$ Crossing of eigenvalues occur for the acoustic modes only at the zero frequency. Thus away from the zero frequency, the acoustic modes are analytical as well. It is easy to check that
  \begin{equation} \label{mu-1} 0 \leq |\mu_s(e^+)_\pm | \leq C \e |\xi|,\end{equation}
 uniformly in $\e, \xi,$ and that
 \begin{equation} \label{mu-1.1}
 \mu_s(e^+)_\pm = \pm \e |\xi| \big( \sqrt{1 + \a^2} + O(\e^2 + |\xi|^2) \big ),
 \end{equation}
uniformly in $\e, |\xi| \leq c_l.$ 
 Regularity at infinity can be directly checked using an exact description of the solutions of \eqref{disp1.2}: one finds that the longitudinal eigenvalues have the form
 \begin{equation} \label{mu-2}
 \mu_s(e_+)_\pm = \pm \e \a |\xi| + \frac{1}{|\xi|} F_\pm(\frac{1}{|\xi|}, \frac{\xi}{|\xi|}, \e), \quad \mu_\pm = \theta_e |\xi| + \frac{1}{|\xi|} G_\pm(\frac{1}{|\xi|}, \frac{\xi}{|\xi|}, \e),
 \end{equation}
where $F_\pm$ and $G_\pm$ are analytical in $(|\xi_0|,\infty) \times \S^{2} \times \R_+,$ for $|\xi_0| > 0.$

 Thus the eigenvalues at $({\tt x}, {\tt y}) = (0,0)$ satisfy \eqref{s1} (regularity in $\e,\xi$ for $\xi$ away from zero, decay at infinity) with $m=1,$ that is, condition (ii) in Assumption \ref{ass_1} is satisfied by the eigenvalues. 

For $\xi$ in a compact subset of $\R^3,$ the eigenvalues evaluated at $({\tt x}, {\tt y})= (\e \theta_e v_e \cdot \xi, \e^2 v_i \cdot \xi)$ are small perturbations of the eigenvalues at $({\tt x}, {\tt y})= (0,0).$ Thus the Klein-Gordon modes are separated from the acoustic modes. For large $\xi,$ the contribution of the convective terms to the eigenvalues is not negligible, but the acoustic modes are all $O(\e |\xi|),$ while the Klein-Gordon modes are $O(|\xi|).$ Thus condition (i) in Assumption \ref{ass_1} is satisfied.

The Klein-Gordon modes are single eigenvalues of \eqref{disp1} and \eqref{disp2}, for all $\xi.$ Hence the eigenvalues and the eigenprojectors corresponding to the Klein-Gordon modes are analytical

 For $\xi \neq 0,$ we let $\{ \xi_1, \xi_2\}$ be an orthonormal basis of $\{\xi\}^\perp.$
 
 At $(\e,0,\xi),$ $\xi \neq 0,$ the eigenvectors associated with the transverse eigenvalues are 
  \begin{eqnarray}  {e}_\pm & := & \big( \frac{|\xi| \xi_2}{\lam_\pm}\;, \xi_1\;, \frac{- i \xi_1}{\lam_\pm}\;, 0\;, \frac{i \e}{\theta_e}\frac{\xi_1}{\lam_\pm}\;, 0\big), \nonumber \\ 
  {e}'_\pm & := & \big( \frac{- |\xi| \xi_1}{\lam_\pm}\;, \xi_2\;, \frac{- i \xi_2}{\lam_\pm}\;, 0\;, \frac{i \e}{\theta_e}\frac{\xi_2}{\lam_\pm}\;, 0 \big),\nonumber \\
  {e}_s(e^-) & := & \frac{1}{\sqrt{1 + |\xi|^2}} \; \big( \xi_1, 0, i |\xi| \xi_2, 0,0,0 \big), \nonumber \\ 
 {e}'_s(e^-) & := & \frac{1}{\sqrt{1 + |\xi|^2}} \; \big( \xi_2, 0, i |\xi| \xi_1, 0,0,0 \big), \nonumber \\
 {e}_s(e^+) & : = & \frac{|\xi|}{|\xi| + \e} \big( \frac{- i \e \xi_2}{\theta_e |\xi|}, 0, 0, 0, \xi_1, 0 \big), \nonumber \\ 
  {e}'_s(e^+) & :=  & \frac{|\xi|}{|\xi| + \e} \big( \frac{i \e \xi_1}{\theta_e |\xi|}, 0, 0, 0, \xi_2, 0 \big), \nonumber 
  \end{eqnarray}

 At $(\e,0,\xi),$ $\xi \neq 0,$ the eigenvectors associated with the longitudinal eigenvalues are 
 $$ \begin{aligned} {f}_\pm  & := \Big( 0 \; ,\;  \frac{\mu_\pm^2 - \theta_e^2 |\xi|^2}{\mu_\pm} \frac{\xi}{|\xi|} \;,\;  -i  \frac{\xi}{|\xi|} \;,\; \frac{-i \theta_e |\xi|}{\mu_\pm} \; , & \frac{i \e}{\theta_e} \frac{\mu_\pm^2 - \theta_e^2 |\xi|^2}{\mu_\pm^2 - \e^2 \a^2 |\xi|^2} \frac{\xi}{|\xi|} \frac{\mu_\pm}{\mu_\pm} \; , \\ & &  \frac{i \a \e^2}{\theta_e} \frac{\mu_\pm^2 - \theta_e^2 |\xi|^2}{\mu_\pm^2 - \e^2 \a^2 |\eta|^2} \frac{|\xi|}{\mu_\pm} \Big), \end{aligned}$$
 $$  {f}_s(e^-) :=  \frac{1}{\sqrt{1 + \theta_e^2 |\xi|^2}} \big( 0 , - i \theta_e \xi, 0, 1, 0, -\frac{1}{\a} \big),$$
and
 $$ \begin{aligned} {f}_s(e^+)_\pm := \Big( 0 \; , & \; \frac{\theta_e}{\e} \tilde \mu_s(e^+)_\pm \frac{\xi}{|\xi|} \; ,    \; \frac{- i \theta_e}{\e}  \frac{\tilde \mu_s(e^+)_\pm \mu_s(e^+)_\pm}{\mu_s(e^+)_\pm^2 - \theta_e^2 |\xi|^2}  \frac{\xi}{|\xi|} \; , \\ & \; \frac{-i \theta_e^2}{\e} \frac{|\xi| \tilde \mu_s(e^+)}{\mu_s(e^+)_\pm^2 - \theta_e^2 |\xi|^2)} \frac{\xi}{|\xi|} \;,\; \frac{i \xi}{|\xi|} \; , \;  \frac{i \a \e |\xi|}{\mu_s(e^+)_\pm} \Big),\end{aligned}$$
 where $$ \tilde \mu_s(e^+)_\pm := \frac{\mu_s(e^+)_\pm^2 - \e^2 \a^2 |\xi|^2}{\mu_s(e^+)_\pm}.$$
From \eqref{mu-1.1}, one sees that $\tilde \mu_s(e^+)_\pm = O(\e |\xi|).$ This yields a description of $f_s(e^+)_\pm$ for small frequencies:
 $$ {f}_s(e^+)_\pm = \big( 0, O(|\xi|), O(\e), O(1), O(1), O(1) \big).$$
The corresponding orthogonal eigenprojectors ($\Pi_j, 1 \leq j \leq n$ with the notation introduced in Assumption \ref{ass_1}) are defined by
 $$ \un e \otimes \un e + \un e' \otimes \un e', \quad \un f \otimes \un f,$$
 where $e = e_{\pm}, e_{s}(e^{-})$ or $e_{s}(e^{+}),$ $f = f_{\pm}, f_{s}(e^{-})$ or $f_{s}(e^{+})_{\pm},$ and $\un e, \un f := \frac{e}{|e|}, \frac{f}{|f|}.$ With the above description of the eigenvalues for small and large $|\xi|,$ it is straightforward to check that all the eigenprojectors, evaluated at $(\e,0),$ are Fourier multipliers in the class $C^\infty {\cal M}^0.$ Thus condition (ii) in Assumption \ref{ass_1} is satisfied.


Condition (iii) follows from the above description of the acoustic modes, as they all have a prefactor $\e.$ 

We now turn to condition \eqref{reg-proj-totaux}. The total eigenprojectors, defined in \eqref{def-proj}, are 
 $$ \Pi_0 = \sum_{\pm} \un e_\pm \otimes \un e_\pm + \un e'_\pm \otimes \un e'_\pm +  \un f_\pm \otimes f_\pm,$$
 and
 $$ \begin{aligned} \Pi_s & = \sum_{\pm}\un e_s(e^\pm) \otimes \un e_s(e^\pm) + \un e'_s(e^\pm) \otimes \un e'_s(e^\pm) \\ & + \un f_s(e^-) \otimes \un f_s(e^-) + e_0 \otimes e_0 + \sum_\pm f_s(e^+)_\pm \otimes \un f_s(e^+)_\pm.\end{aligned}$$
 Given $(\e_0, u_0, \xi_0) \in [0,1] \times \R^{14} \times \R^3,$ for $(\e, u, \xi)$ in a neighborhood of $(\e_0, u_0, \xi_0),$
 \begin{equation} \label{p0} \Pi_0(\e, u,\xi) = \frac{1}{2 i \pi} \int_{\Gamma_{0+}} (z - {\cal A}(\e,u,\xi))^{-1} dz + \frac{1}{2 i \pi} \int_{\Gamma_{0-}} (z - {\cal A}(\e,u,\xi))^{-1} dz,\end{equation}
 and 
 \begin{equation} \label{ps} \Pi_s(\e,u,\xi) = \frac{1}{2 i \pi} \int_{\Gamma_s} (z - {\cal A}(\e,u,\xi))^{-1} dz,\end{equation}
 where $\Gamma_{0+}$ (resp. $\Gamma_{0-}$) is a contour enclosing all the positive (resp. negative) Klein-Gordon eigenvalues at $(\e_0, u_0, \xi_0),$ and no other eigenvalues of ${\cal A},$ and $\Gamma_s$ is a contour enclosing all the acoustic eigenvalues at $(\e_0, u_0, \xi_0),$ and no other eigenvalues of ${\cal A}.$ The above description of the eigenvalues show in particular that such contours do exist.   
 It follows from these representations that $\Pi_0$ and $\Pi_s$ are analytical in $\e, u, \xi.$ To investigate the behaviour for large $\xi,$ notice that ${\cal A}$ can be written, for $|\xi| \neq 0,$
 $$ {\cal A}(\e,u,\xi) = \un {\cal A}_0(\e) + \un {\cal A}_1(\e,u,\xi),$$
 where ${\cal A}_1$ is linear in $\xi,$ so that, using polar coordinates $\xi = (|\xi|, \theta) \in \R_+^* \times \S^{2},$
 $$ {\cal A}(\e,u,\xi) = |\xi| ( \frac{1}{|\xi|} \un {\cal A}_0(\e) + \un {\cal A}_1(\e,u,\theta)).$$
 Thus an eigenvector $e(\e, u, |\xi|, \theta)$ of ${\cal A}(\e, u, \xi),$ associated with an eigenvalue $\lam(\e,u,|\xi|, \theta),$ equals an eigenvector $\tilde e(\e, u, \frac{1}{|\xi|}, \theta)$ of 
 \begin{equation}\label{ondes-longues} \tilde {\cal A}(\e, u, |\xi|, \theta) := |\xi| \un {\cal A}_0(\e) + \un {\cal A}_1(\e,u,\theta),\end{equation}
 associated with the eigenvalue $\frac{1}{|\xi|} \lam(\e, u, \frac{1}{|\xi|}, \theta).$ The operator $\tilde {\cal A}$ is the long-wave operator associated with ${\cal A},$ introduced in \cite{T0} in the study of the short-wave limit.

 For small $|\xi|,$ the eigenvalues of $\tilde {\cal A}$ split into ``acoustic'' eigenvalues, of size $O(\e),$ and ``Klein-Gordon'' eigenvalues, of size $O(1).$ The total eigenprojectors,
 $$ \tilde \Pi_0 = \frac{1}{2 i \pi} \int_{\tilde \Gamma_{0+}} (z - \tilde {\cal A}(\e,u,|\xi|, \theta))^{-1} dz + \frac{1}{2 i \pi} \int_{\Gamma_{0-}} (z - \tilde {\cal A}(\e,u,|\xi|,\theta))^{-1} dz,$$
 and 
 $$ \tilde \Pi_s = \frac{1}{2 i \pi} \int_{\Gamma_s} (z - \tilde {\cal A}(\e,u,|\xi|, \theta))^{-1} dz,$$ 
 are analytical in $\e, u, |\xi|, \theta,$ for $|\xi|$ in a neighborhood of 0. There holds
 $$ \Pi_j(\e, u, \xi) = \tilde \Pi_j(\e, u, \frac{1}{|\xi|}, \frac{\xi}{|\xi|}), \quad j= 0, s,$$
 and thus the total eigenprojectors $\Pi_0, \Pi_s$ of ${\cal A}$ are analytical in $\e, u, \frac{1}{|\xi|}, \frac{\xi}{|\xi|},$ for $\xi$ in a neighborhood of $\infty.$ This implies that they satisfy decay estimates of the form \eqref{s1}, and \eqref{reg-proj-totaux} is proved. 
 
 To conclude this section, we now indicate how the Klein-Gordon eigenvalues and eigenvectors depend on ${\tt x}$ and ${\tt y}.$ These descriptions are needed in the evaluation of the interaction coefficients that enter Assumption \ref{ass_6}.
 
  The eigenvalues satisfy
$$ \d_{{\tt x}}\lam_{\pm}  = \frac{1}{2(1 + |\xi|^2)} + O(\e^2), \quad \d_{{\tt x}}\mu_{\pm}  = \frac{\theta_e^2 |\xi|^2}{2(1 + \theta_e^2 |\xi|^2)} + O(\e), $$
 at $({\tt x}, {\tt y}) = (0,0),$ locally uniformly in $\xi.$ 
 
  The eigenvectors are
  \begin{eqnarray}  {\un e}_\pm & := & \Big( \frac{|\eta| \eta_2}{\lam_\pm}\;, \eta_1\;, \frac{- i \eta_1}{\lam_\pm - {\tt x}}\;, 0\;, \frac{i \e}{\theta_e}\frac{\eta_1}{\lam_\pm- {\tt y}}\;, 0\Big), \nonumber \\ 
  {\un e}'_\pm & := & \Big( \frac{- |\eta| \eta_1}{\lam_\pm}\;, \eta_2\;, \frac{- i \eta_2}{\lam_\pm - {\tt x}}\;, 0\;, \frac{i \e}{\theta_e}\frac{\eta_2}{\lam_\pm - {\tt y}}\;, 0 \Big),\nonumber 
 \end{eqnarray} 
  and 
  $$ \begin{aligned} {\un f}_\pm  & :=  \Big( 0 \; , & \frac{(\mu_\pm - {\tt x})^2 - \theta_e^2 |\eta|^2}{\mu_\pm - {\tt x}} \frac{\eta}{|\eta|} \; , \;   -i  \frac{\eta}{|\eta|} \; , \;   \frac{-i \theta_e |\eta|}{\mu_\pm - {\tt x}} \; , \\ & & \frac{i \e}{\theta_e} \frac{(\mu_\pm - {\tt x})^2 - \theta_e^2 |\eta|^2}{(\mu_\pm - {\tt y})^2 - \e^2 \a^2 |\eta|^2} \frac{\eta}{|\eta|} \frac{\mu_\pm - {\tt y}}{\mu_\pm - {\tt x}} \; ,  \\  & & \frac{i \e^2}{\theta_e} \frac{(\mu_\pm - {\tt x})^2 - \theta_e^2 |\eta|^2}{(\mu_\pm - {\tt y})^2 - \e^2 \a^2 |\eta|^2} \frac{|\eta|}{\mu_\pm - {\tt x}} \Big). \end{aligned}$$
 
\end{subsubsection}


\begin{subsubsection}{Resonances}  \label{3-2-1}

 We use here the description of the eigenvalues given in section \ref{3-2-2}. The resonance equation \eqref{res-eq} for $j, k \leq n_0$ (Klein-Gordon/Klein-Gordon resonances) and $\e = 0$ is
 \begin{equation} \label{3fev3} \sqrt{1 + |\xi|^2} - \sqrt{1 + \theta_e^2 |\xi|^2} - 1 = 0 .\end{equation}
 If $\theta_e$ is small enough, we choose $c_m = 1$ and $C_m = 2.$ The left-hand side in \eqref{3fev3} is bounded away from zero for $|\xi| \notin [1, 2],$ and hypothesis (0-0) in Assumption \ref{ass_3} is satisfied. 
 
\begin{figure}[t]
\begin{center}
\scalebox{0.5}{\input{res-0-0}}
\caption{(0-0) resonances.}
\label{fig-0-0}
\end{center}
\end{figure}
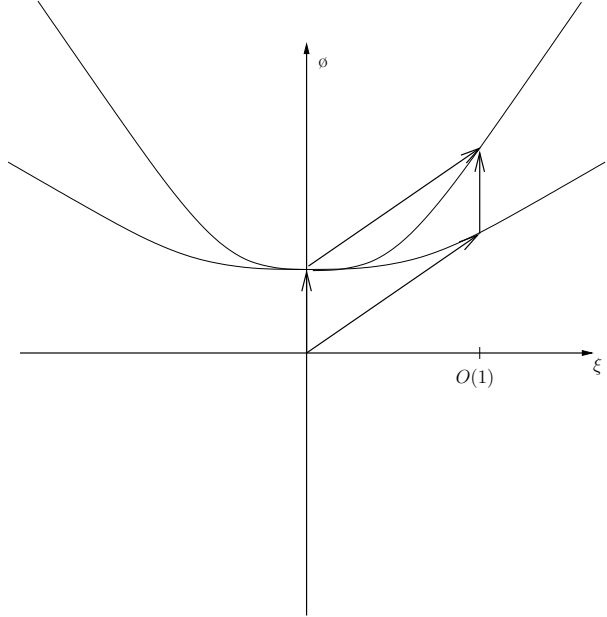

 The resonance equation \eqref{res-eq} for $j \leq n_0 < k$ (Klein-Gordon/acoustic resonance) is 
 \begin{equation} \label{3fev4} \sqrt{1 + \kappa^2 |\xi^2|} - 1 - \mu = 0.\end{equation}
  where $\kappa = 1$ or $\kappa = \theta_e,$ and $\mu = 0$ or $\mu = \mu_s(e^+)_\pm.$  If $\theta_e$ and $\e$ are small enough, we can choose $c_l = 1/2.$ Then the left-hand side in \eqref{3fev4} is bounded away from zero for $|\xi| \geq 1/2.$

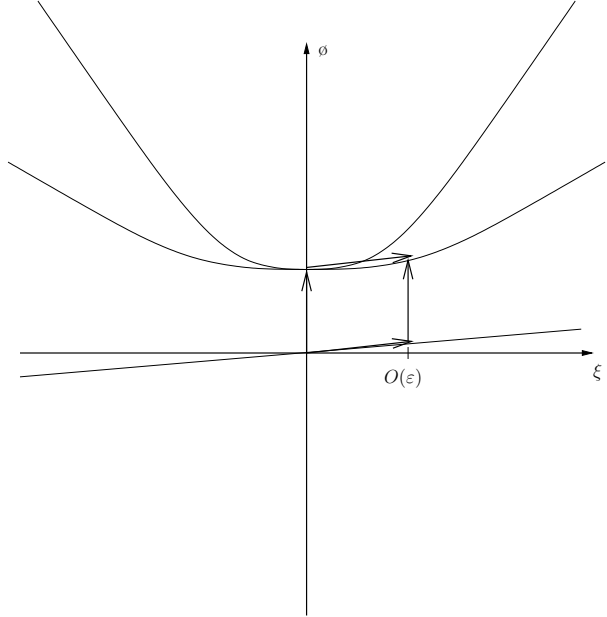
\begin{figure}[t]
\begin{center}
\scalebox{0.5}{\input{res-0-s}}
\caption{(0-s) resonances.}
\label{fig-0-s}
\end{center}
\end{figure}

 The resonance equation \eqref{res-eq-2} (secondary Klein-Gordon/acoustic resonances) is
  \begin{equation} \label{3fev5} \sqrt{1 + \kappa^2 |\xi|^2} - (p + p') = 0 ,\end{equation}
  where $p, p \in \{-1, 1\},$ $\kappa = 1$ or $\kappa = \theta_e.$ The left-hand side in \eqref{3fev5} is bounded away from 0 for $|\xi| \leq 1,$ and hypothesis (0-0-s) in Assumption \ref{ass_3} is satisfied.    

 The resonances are pictured on Figures \ref{fig-0-0} to \ref{fig-0-0-s}.

\begin{figure}[t]
\begin{center}
\scalebox{0.5}{\input{res-0-0-s}} 
\caption{(0-0-s) resonances.}
\label{fig-0-0-s}
\end{center}
\end{figure}
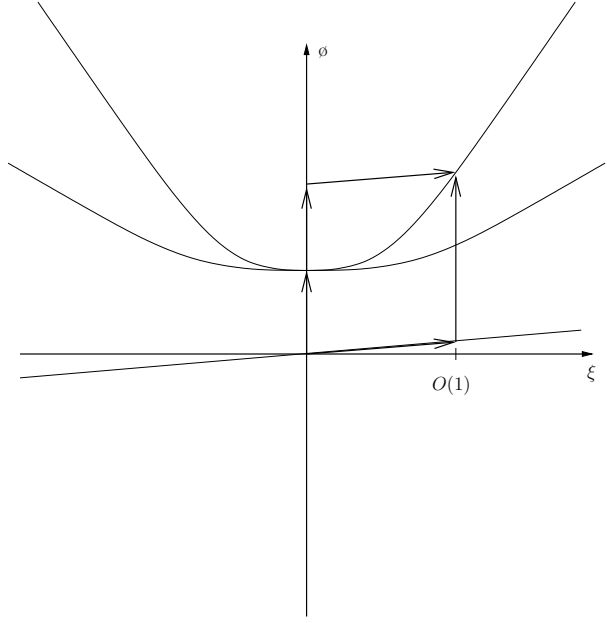

\end{subsubsection}



\begin{subsubsection}{Transparency} \label{3-2-4}

 We check in this section that Assumptions \ref{ass_6} and \ref{ass_sym} are satisfied.   
 
 With the definition of ${\cal A}$ given at the beginning of section \ref{application}, one has ${\cal A}^{(2)} = 0.$ It is straightforward to check that $\Pi_s(0,0) {\cal D}(u_a) \Pi_0(0,0),$ where ${\cal D}$ is defined by \eqref{d}, has size $O(\e),$ and \eqref{interaction-d} is satisfied.

 Up to terms of size $O(\e^2),$ the source term defined in \eqref{b} is 
 $$ {\cal B}(u_a) u = \left(\begin{array}{c} 0 \\ n_e {\bf v}_{e} + \e {\bf n}_{e1} v_e \\ - \theta_e ( {\bf v}_{e} \times B + ({\bf v}_e \cdot \xi) v_e + \e {\bf B}_1 \times v_e + \e v_e \cdot \nabla {\bf v}_e) \\ -\theta_e ({\bf v}_e \cdot \xi) n_e \\ 0 \\ 0 \end{array}\right).$$
 It depends on the coordinates ${\bf B},$ ${\bf v}_e,$ ${\bf n}_e,$ of the approximate solution $u_a,$ evaluated at $\theta = \om t/\e^2.$
 
   One computes,
 $$ e_+^* {\cal B}(u_a) e_s(e^-) = \frac{i \theta_e}{\lam_+^2} \frac{{\bf v}_{e0} \cdot {\xi}}{|\xi|} + O(\e),$$
 and also
 $$ f_+^* {\cal B}(u_a) f_s(e^-) = \frac{\mu_+^2 - \theta_e^2 |\xi|^2}{\mu_+^2} \frac{{\bf v}_{e0} \cdot \xi}{|\xi|} + O(\e).$$
 In particular, these interaction coefficients are not small for small frequencies. This shows that the non-transparency condition \eqref{non-transp} is satisfied. 

 We now check that the transparency condition of Assumption \ref{ass_6} is satisfied. 

The symbol $\rho$ defined in \eqref{rho} is $ \rho = \rho^{(0)} + O(\e),$ where
 $$ \begin{aligned} \rho^{(0)} := \Pi_s \sum_{|\a| = 1} \Big( & \d_\xi^\a (\Pi_s {\cal A}^{(0)}) \d_v \Pi_0 \cdot \d_x^\a u_a + \d_\xi^\a \Pi_s {\cal A}^{(1)}(\d_x^\a u_a) \\ & + \Pi_s {\cal A}^{(0)} \d_\xi^\a \Pi_0 \d_v \Pi_0 \cdot \d_x^\a u_a - \d_\xi^\a \Pi_s {\cal B}(0, \d_x^\a u_a) \Big) \Pi_0,\end{aligned}$$
 where, unless otherwise noted, the symbols and their derivatives are evaluated at $(\e,v) = (0,0).$ 

 Assumption \ref{ass_6} is a transparency assumption for the interaction coefficient $\d_u B^{\tt r}(\e,0) \cdot u_a.$ Direct computations, using the description of the eigenvalues and the eigenvectors given above, yield the bound
 \begin{equation} \label{int-coeff}
 | \Pi_k(\e,0) (\d_u B^{\tt r}(\e,0) \cdot u_a) \Pi_j(\e,0) | \leq C_B ( |\xi|^2 + \e |\xi|),
 \end{equation}
 for $j \leq n_0 < k,$ uniformly in $\e, t, x$ and $|\xi| \leq c_l.$ Now \eqref{int-coeff} and the above description of the phases imply \eqref{interaction0}, as follows. Let $\delta > 0$ be given. 
 \begin{itemize}
 \item If $|\xi| \leq \delta \e,$ then \eqref{int-coeff} directly implies \eqref{interaction0}, with $C = C_B \delta (1 + \delta).$   
  \item If $|\xi| > \delta \e,$ then for $0 < \e < \e_0,$ $|\Phi_{j,k,p}(\e)| > C_0(\delta) |\xi|^2,$ with
   $$ C_0(\delta) :=  \frac{3\theta_e^2}{8} - \frac{C(\a,\e_0)}{\delta},$$
 for some nondecreasing function $C(\a,\e_0).$ If $\delta$ is chosen to be large enough, then $C_0(\delta) > 0,$ and \eqref{interaction0} is satisfied with $$C = C_B(\frac{1}{\delta} + 1) C_0(\delta)^{-1}.$$   
  \end{itemize}

\bigskip

Assumption \ref{ass_sym} is a symmetrizability condition for the interaction coefficients
 $$ B_0 = \Pi_0 {\cal B}(u_a) \Pi_0, \quad B_s = \Pi_s {\cal B} \Pi_s.$$
 As these interaction coefficients enter the equation with a prefactor $\frac{1}{\e},$ it is sufficient to consider the leading term in $\e$ in ${\cal B}(u_a).$ In particular, the contribution of ${\cal A}^{(1)}$ in $B_0$ and $B_s$ is only $O(\e).$ We can write indeed
 \begin{equation} \label{a1} \ompi_\e({\cal A}^{(1)}(u)) u_a = \e \ompi_1({\cal A}^{(1)}(u)) u_a,\end{equation}
 because ${\cal A}^{(1)}$ is linear in $\xi,$ and the $H^s_\e$ norm of \eqref{a1} is bounded by
 $$ \e C(\| u_a \|_{1,s+d_0}) \| u \|_{\e,s}.$$
 The auto-interaction coefficients are all purely imaginary:
 $$ e^* {\cal B}(u_a) e \in i \R,$$
 where $e$ is any eigenvector of the Euler-Maxwell equations.
 The other interaction coefficients between the Klein-Gordon modes are
 $$ \begin{aligned} 
  e_\pm^* {\cal B}(u_a) e'_\pm & = 0, & e_\pm^* {\cal B}(u_a) f_\pm & = - i \theta_e \frac{|\xi|}{\mu_\pm} {\bf v}_e \cdot \xi_1,\\
  (e'_\pm)^* {\cal B}(u_a) e_\pm & = 0,  & (e'_\pm)^* {\cal B}(u_a) f_\pm & = - i \theta_e \frac{|\xi|}{\mu_\pm} {\bf v}_e \cdot \xi_2,\\
 f_\pm^* {\cal B}(u_a) e_\pm & =  - i \theta_e \frac{|\xi|}{\lam_\pm} {\bf v}_e \cdot \xi_1, \quad & f_\pm^* {\cal B}(u_a) e'_+ & = - i \theta_e \frac{|\xi|}{\lam_\pm} {\bf v}_e \cdot \xi_2,\end{aligned}$$
 where for all $\xi \neq 0,$ $\{ \xi_1, \xi_2\}$ is an orthonormal basis of $\{\xi\}^\perp.$
  Let $S_0 = S_0^{(+)} + S_0^{(-)},$ where
 $$ S_0^{(\pm)} := e_\pm \otimes e_\pm + e'_\pm \otimes e'_\pm + \frac{\lam_\pm}{\mu_\pm} f_\pm \otimes f_\pm.$$
 Then $S_0$ is a symmetrizer, in the sense that 
 $$ \frac{1}{\e} (S_0 B_0 + (S_0 B_0)^*) \in C^\infty {\cal M}^0.$$
 Finally, we turn to the interaction coefficients between acoustic modes. Remark first that the divergence-free condition for the magnetic field,
 \begin{equation} \label{div-free} \nabla \cdot B = 0,\end{equation}
 is equivalent to
 $$ \ompi_\e(e_0) u = 0.$$
 Condition \eqref{div-free} is propagated by the equations, that is, the solution belongs to the orthogonal of the image of $\ompi_\e(e_0)$ in $L^2(\R^3)$ if the initial datum does. Thus we can overlook the interaction coefficient involving $e_0.$ At first order in $\e,$ the interaction coefficients with the other acoustic modes all vanish, except for
 \begin{eqnarray} f_s(e_-)^* {\cal B}(u_a) f_s(e_+)_\pm & = & \frac{{-\bf v}_e \cdot \xi}{\sqrt{1 + \theta_e^2 |\xi|^2}} \frac{\mu_s(e_+)_\pm}{\e \theta_e} \frac{- \theta_e^2 |\xi|^2}{\mu_s(e_+)_\pm^2 - \theta_e^2 |\xi|^2}  \nonumber,\\
 (f_s(e_+)_\pm)^* {\cal B}(u_a)  f_s(e_-) & = & \frac{{\bf v}_e \cdot \xi}{\sqrt{1 + \theta_e^2 |\xi|^2}} \frac{\mu_s(e_+)_\pm}{\e \theta_e}  \nonumber \end{eqnarray}
 It follows from the description of the eigenvalues given in section \ref{3-2-2} that, for $\e$ small enough, 
 $$ \frac{1}{2} \leq \frac{- \theta_e^2 |\xi|^2}{\mu_s(e_+)_\pm^2 - \theta_e^2 |\xi|^2} \leq 2.$$ 
 uniformly in $\xi \in \R^3.$ Thus we can take $\g = 2$ in Assumption \ref{ass_sym}. Let 
 $$ \begin{aligned}
 S_s := & e_0 \otimes e_0 + e_s(e_-) \otimes e_s(e_-) + e_s(e_+) \otimes e_s(e_+) + f_s(e_-) \otimes f_s(e_-) \nonumber\\
 & + \frac{- \theta_e^2 |\xi|^2}{\mu_s(e_+)_+^2 - \theta_e^2 |\xi|^2} f_s(e_+)_+ \otimes  f_s(e_+)_+ \\ & + \frac{- \theta_e^2 |\xi|^2}{\mu_s(e_+)_-^2 - \theta_e^2 |\xi|^2}  f_s(e_-)_+ \otimes  f_s(e_+)_-.\end{aligned}$$
 Then
 $$ \frac{1}{\e} (S_s B_s + (S_s B_s)^*) \in C^\infty {\cal M}^0.$$
 Because $S_0$ and $S_s$ are diagonal in a basis of eigenvectors of ${\cal A},$ the matrices $S_0 A_0$ and $S_s A_s$ are hermitian. Finally,
 $$ S := \left(\begin{array}{cc} S_0 & 0 \\ 0 & S_s \end{array}\right),$$
 defines a symmetrizer that satisfies Assumption \ref{ass_6}. 

 \end{subsubsection}

\end{subsection}

\end{section}




\end{document}

%% file: troncs.tex
\begin{picture}(0,0)%
\epsfig{file=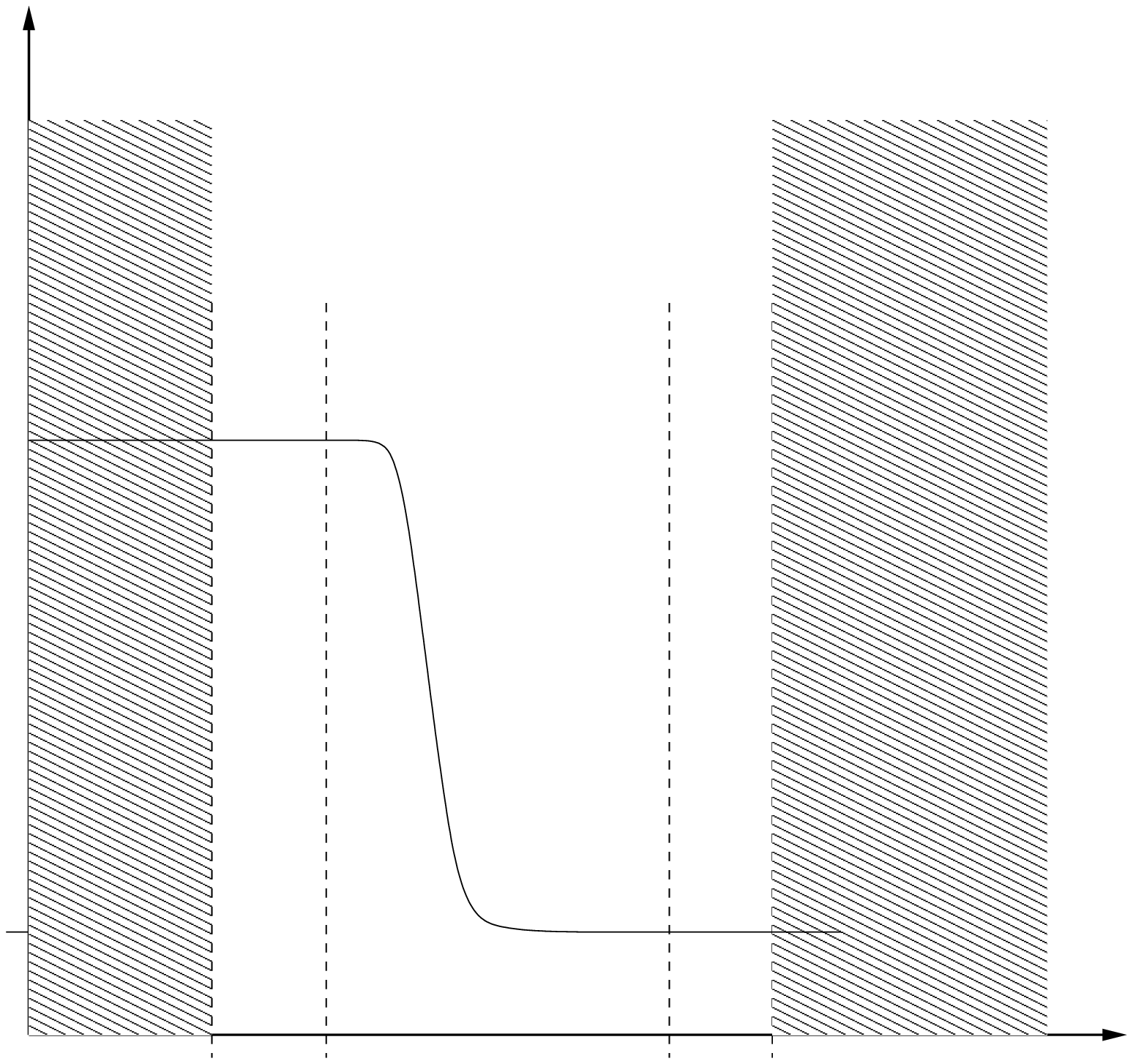}%
\end{picture}%
\setlength{\unitlength}{3947sp}%
\begingroup\makeatletter\ifx\SetFigFont\undefined%
\gdef\SetFigFont#1#2#3#4#5{%
  \reset@font\fontsize{#1}{#2pt}%
  \fontfamily{#3}\fontseries{#4}\fontshape{#5}%
  \selectfont}%
\fi\endgroup%
\begin{picture}(7822,7407)(1801,-6799)
\put(7276, 14){\makebox(0,0)[lb]{\smash{{\SetFigFont{12}{14.4}{\rmdefault}{\mddefault}{\updefault}(0-0-s) resonances}}}}
\put(1801,-2536){\makebox(0,0)[lb]{\smash{{\SetFigFont{14}{16.8}{\rmdefault}{\mddefault}{\updefault}$1$}}}}
\put(4951,-2761){\makebox(0,0)[lb]{\smash{{\SetFigFont{14}{16.8}{\rmdefault}{\mddefault}{\updefault}$\chi_\e$}}}}
\put(1801,-5761){\makebox(0,0)[lb]{\smash{{\SetFigFont{14}{16.8}{\rmdefault}{\mddefault}{\updefault}$\e$}}}}
\put(3301,-6736){\makebox(0,0)[lb]{\smash{{\SetFigFont{14}{16.8}{\rmdefault}{\mddefault}{\updefault}$c_l$}}}}
\put(4126,-6736){\makebox(0,0)[lb]{\smash{{\SetFigFont{14}{16.8}{\rmdefault}{\mddefault}{\updefault}$c_0$}}}}
\put(6376,-6736){\makebox(0,0)[lb]{\smash{{\SetFigFont{14}{16.8}{\rmdefault}{\mddefault}{\updefault}$c_1$}}}}
\put(7051,-6736){\makebox(0,0)[lb]{\smash{{\SetFigFont{14}{16.8}{\rmdefault}{\mddefault}{\updefault}$c_m$}}}}
\put(9301,-6736){\makebox(0,0)[lb]{\smash{{\SetFigFont{14}{16.8}{\rmdefault}{\mddefault}{\updefault}$|\xi|$}}}}
\put(2551, 14){\makebox(0,0)[lb]{\smash{{\SetFigFont{12}{14.4}{\rmdefault}{\mddefault}{\updefault}(0-s) resonances}}}}
\put(7276,464){\makebox(0,0)[lb]{\smash{{\SetFigFont{12}{14.4}{\rmdefault}{\mddefault}{\updefault}(0-s) resonances}}}}
\end{picture}%

%% file: char1.tex
\begin{picture}(0,0)%
\epsfig{file=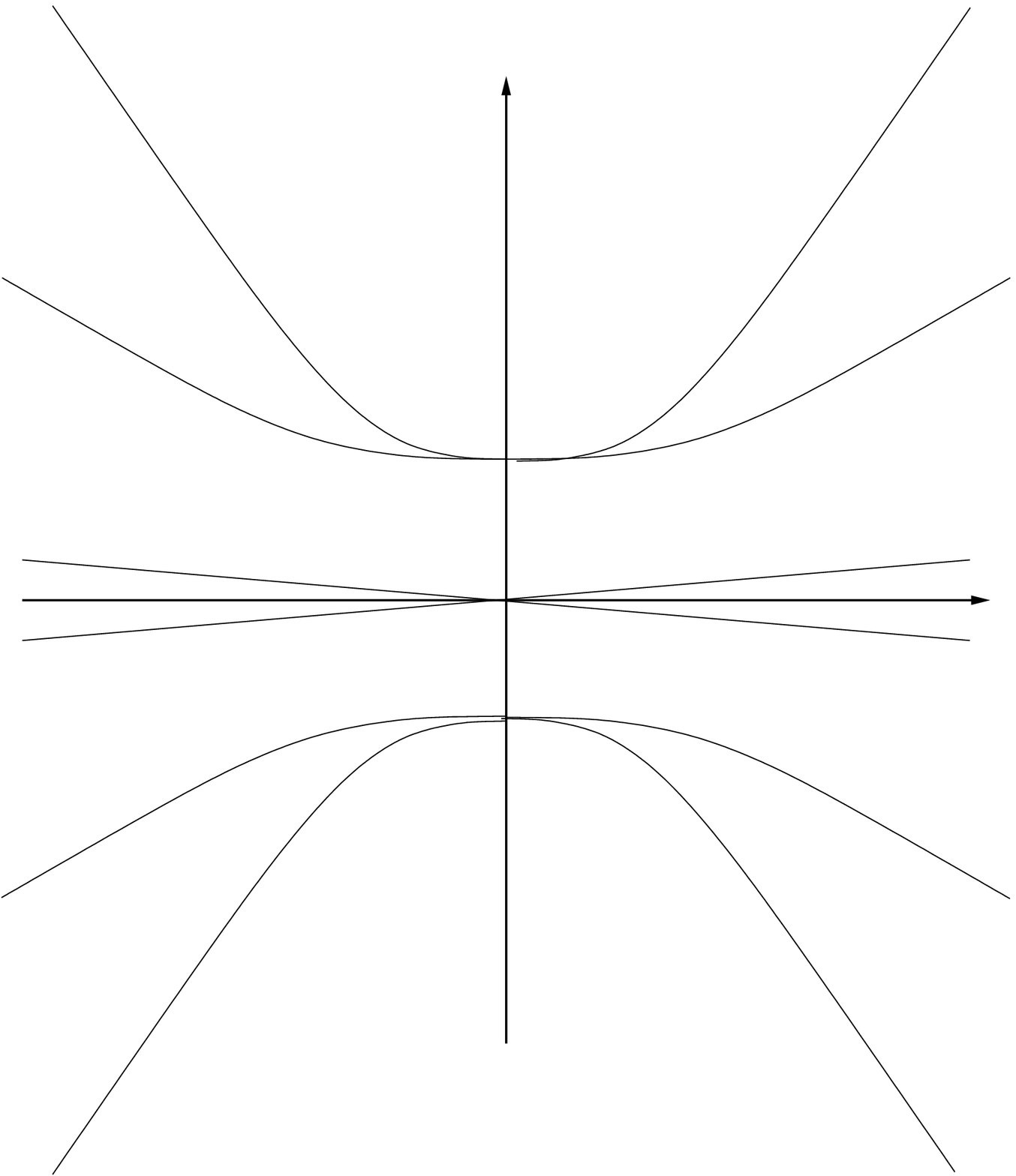}%
\end{picture}%
\setlength{\unitlength}{3947sp}%
\begingroup\makeatletter\ifx\SetFigFont\undefined%
\gdef\SetFigFont#1#2#3#4#5{%
  \reset@font\fontsize{#1}{#2pt}%
  \fontfamily{#3}\fontseries{#4}\fontshape{#5}%
  \selectfont}%
\fi\endgroup%
\begin{picture}(7529,8724)(2234,-8248)
\put(9676,-4036){\makebox(0,0)[lb]{\smash{{\SetFigFont{14}{16.8}{\rmdefault}{\mddefault}{\updefault}$\xi$}}}}
\put(6151,-211){\makebox(0,0)[lb]{\smash{{\SetFigFont{14}{16.8}{\rmdefault}{\mddefault}{\updefault}$\o$}}}}
\end{picture}%

%% file: res-0-0.tex
\begin{picture}(0,0)%
\epsfig{file=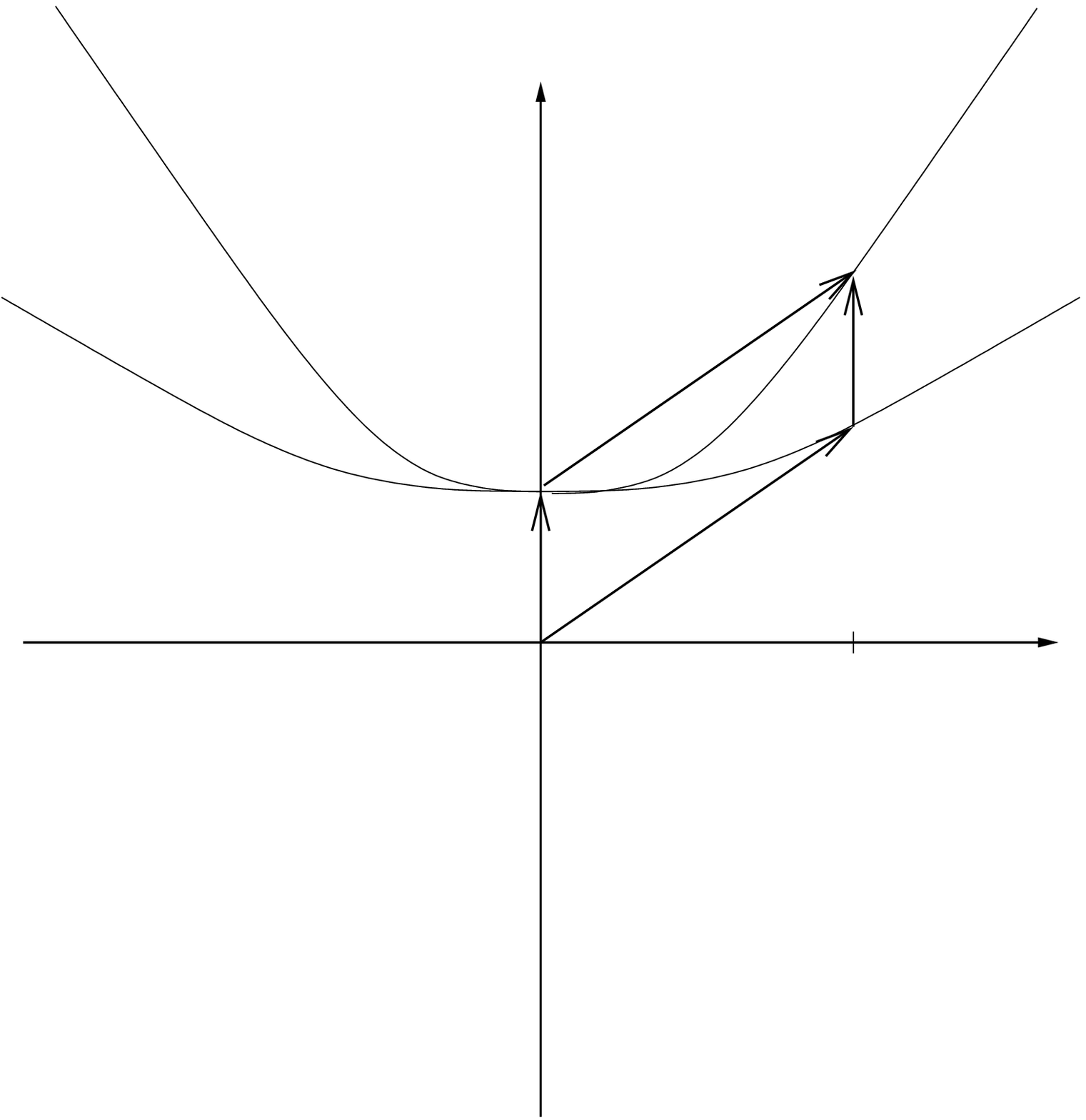}%
\end{picture}%
\setlength{\unitlength}{3947sp}%
\begingroup\makeatletter\ifx\SetFigFont\undefined%
\gdef\SetFigFont#1#2#3#4#5{%
  \reset@font\fontsize{#1}{#2pt}%
  \fontfamily{#3}\fontseries{#4}\fontshape{#5}%
  \selectfont}%
\fi\endgroup%
\begin{picture}(7524,7759)(2239,-7283)
\put(9601,-4186){\makebox(0,0)[lb]{\smash{{\SetFigFont{14}{16.8}{\rmdefault}{\mddefault}{\updefault}$\xi$}}}}
\put(7876,-4336){\makebox(0,0)[lb]{\smash{{\SetFigFont{14}{16.8}{\rmdefault}{\mddefault}{\updefault}$O(1)$}}}}
\put(6151,-361){\makebox(0,0)[lb]{\smash{{\SetFigFont{14}{16.8}{\rmdefault}{\mddefault}{\updefault}$\o$}}}}
\end{picture}%

%% file: res-0-s.tex
\begin{picture}(0,0)%
\epsfig{file=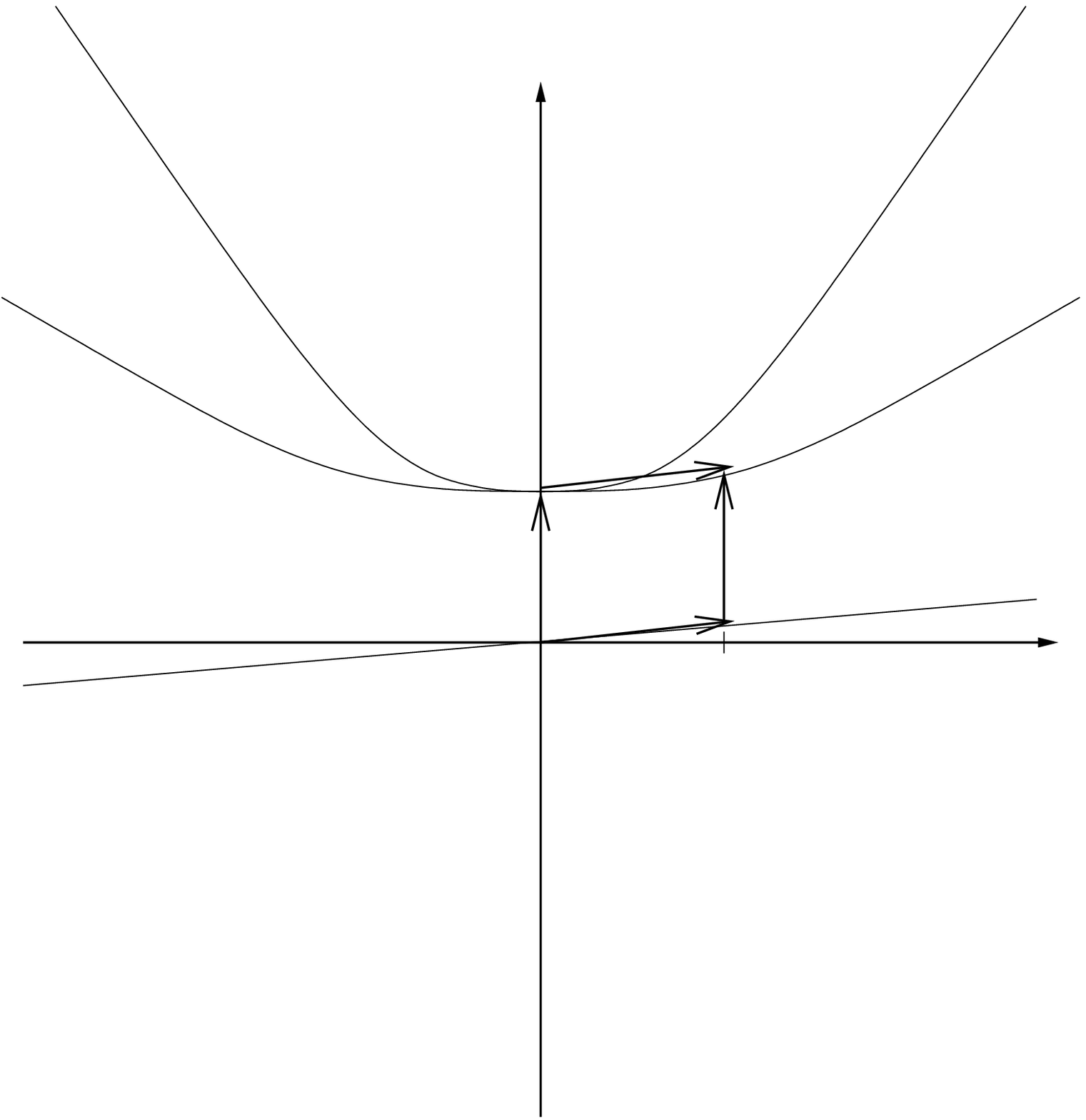}%
\end{picture}%
\setlength{\unitlength}{3947sp}%
\begingroup\makeatletter\ifx\SetFigFont\undefined%
\gdef\SetFigFont#1#2#3#4#5{%
  \reset@font\fontsize{#1}{#2pt}%
  \fontfamily{#3}\fontseries{#4}\fontshape{#5}%
  \selectfont}%
\fi\endgroup%
\begin{picture}(7524,7759)(2239,-7283)
\put(9601,-4261){\makebox(0,0)[lb]{\smash{{\SetFigFont{14}{16.8}{\rmdefault}{\mddefault}{\updefault}$\xi$}}}}
\put(6151,-211){\makebox(0,0)[lb]{\smash{{\SetFigFont{14}{16.8}{\rmdefault}{\mddefault}{\updefault}$\o$}}}}
\put(6976,-4336){\makebox(0,0)[lb]{\smash{{\SetFigFont{14}{16.8}{\rmdefault}{\mddefault}{\updefault}$O(\e)$}}}}
\end{picture}%

%% file: res-0-0-s.tex
\begin{picture}(0,0)%
\epsfig{file=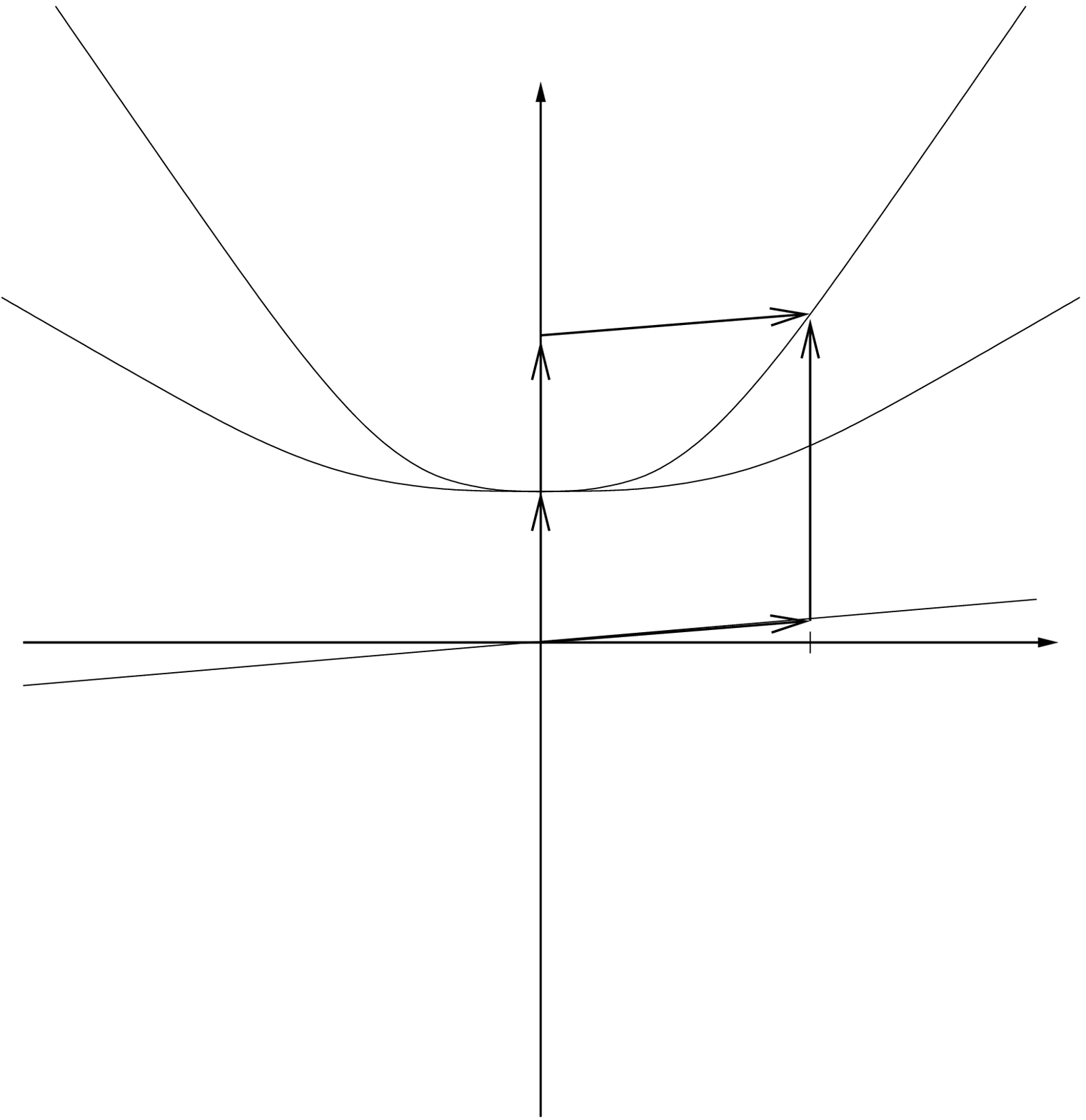}%
\end{picture}%
\setlength{\unitlength}{3947sp}%
\begingroup\makeatletter\ifx\SetFigFont\undefined%
\gdef\SetFigFont#1#2#3#4#5{%
  \reset@font\fontsize{#1}{#2pt}%
  \fontfamily{#3}\fontseries{#4}\fontshape{#5}%
  \selectfont}%
\fi\endgroup%
\begin{picture}(7524,7759)(2239,-7283)
\put(9526,-4261){\makebox(0,0)[lb]{\smash{{\SetFigFont{14}{16.8}{\rmdefault}{\mddefault}{\updefault}$\xi$}}}}
\put(6151,-211){\makebox(0,0)[lb]{\smash{{\SetFigFont{14}{16.8}{\rmdefault}{\mddefault}{\updefault}$\o$}}}}
\put(7576,-4411){\makebox(0,0)[lb]{\smash{{\SetFigFont{14}{16.8}{\rmdefault}{\mddefault}{\updefault}$O(1)$}}}}
\end{picture}%